\documentclass[sii,nonatbib]{ipart}

\usepackage{subfigure}
\usepackage{mathtools}
\usepackage{mathdots}
\usepackage{mathrsfs}
\usepackage{mathabx}
\usepackage[dvipsnames,table]{xcolor}

\newtheorem{theorem}{Theorem}[section]
\newtheorem{lemma}{Lemma}[section]
\newtheorem{proposition}{Proposition}[section]
\newtheorem{corollary}{Corollary}[section]

\newtheorem{remark}{Remark}[section]

\startlocaldefs
\newcommand{\be}{\mathbf{e}}
\newcommand{\bp}{\mathbf{p}}
\newcommand{\bq}{\mathbf{q}}
\newcommand{\br}{\mathbf{r}}
\newcommand{\bA}{\mathbf{A}}
\newcommand{\bB}{\mathbf{B}}
\newcommand{\bS}{\mathbf{S}}
\newcommand{\bQ}{\mathbf{Q}}
\newcommand{\bP}{\mathbf{P}}
\newcommand{\bL}{\mathbf{L}}
\newcommand{\bR}{\mathbf{R}}
\newcommand{\bW}{\mathbf{W}}
\newcommand{\bI}{\mathbf{I}}
\newcommand{\bX}{\mathbf{X}}
\newcommand{\bU}{\mathbf{U}}
\newcommand{\bV}{\mathbf{V}}
\newcommand{\bZ}{\mathbf{Z}}
\newcommand{\bSigma}{\mathbf{\Sigma}}

\newcommand{\bw}{\mathbf{w}}
\DeclareMathOperator{\spanm}{span}
\DeclareMathOperator{\rank}{rank}
\DeclareMathOperator{\vecm}{vec}
\DeclareMathOperator{\trace}{Trace}
\DeclareMathOperator*{\argmin}{arg\,min}
\newcommand{\diag}{\operatorname{diag}}

\DeclareMathOperator{\sto}{subject\,to}
\newcommand{\kron}{\otimes}

\newcommand{\NN}{\mathbb{N}}
\newcommand{\CC}{\mathbb{C}}
\newcommand{\len}{{N}}

\endlocaldefs

\pubyear{2022}
\volume{0}
\issue{0}
\firstpage{1}
\lastpage{1}

\begin{document}

\begin{frontmatter}

\title{Hankel low-rank approximation and completion in time series analysis and forecasting: a brief review\protect}

\begin{aug}

\author{\inits{J.}\fnms{Jonathan} \snm{Gillard}\thanksref{t2}\ead[label=e1]{gillardjw@cardiff.ac.uk}}
\address{School of Mathematics, Cardiff University, UK\\\printead{e1}}
\and\author{\inits{K.}\fnms{Konstantin} \snm{Usevich}\ead[label=e2]{konstantin.usevich@univ-lorraine.fr}}
\address{Universit\'{e} de Lorraine, CNRS, CRAN, Nancy, France\\\printead{e2}}
\thankstext{t2}{Corresponding author.}
\end{aug}

\begin{abstract}
In this paper we offer a review and bibliography of work on Hankel low-rank approximation and completion, with particular emphasis on how this methodology can be used for time series analysis and forecasting. We begin by describing possible formulations of the problem and offer commentary on related topics and challenges in obtaining globally optimal solutions. Key theorems are provided, and the paper closes with some expository examples.
\end{abstract}




\end{frontmatter}

\section{Introduction}

Low-rank techniques have become one of the great success stories in computational mathematics,
linear algebra, engineering and statistics over the last decade. They now constitute the dominant
paradigm for solving large-scale problems from computational fluid dynamics, through computational
meteorology, to molecular modelling. They are also used in inverse problems, where restricting the
rank acts as a kind of regularisation to make ill-posed problems stable and tractable. Arguably one of the most exciting developments in statistics in recent years has been the exploitation of (approximate) low-rank structure of observed data. This is the idea that underpins the success of compressed sensing and matrix completion.

In statistics, approximating a given matrix by one of small rank, is closely related to principal component analysis and variations of it \cite{shen2008sparse}, factor analysis \cite{srebro2003weighted} and total least squares \cite{markovsky2007overview}. Such low-rank approximations are now essential for the application of kernel methods to large-scale learning problems \cite{kulis2006learning}. Their motivation and need is well described in the paper \cite{udell2019big}, where it is asserted that matrices of low-rank are pervasive in data science and that any sufficiently large matrix with small spectral norm can be well approximated by a low-rank matrix.

There has been much recent work in the application and development of low-rank approximation methods to solve typical problems of time series analysis, namely those of modelling and forecasting. Generally speaking the main idea of the methodology is to embed an observed time series into a Hankel matrix, the rank of which is taken to be a measure of the complexity of the time series. There are rich and historical results connecting properties of a Hankel matrix and its embedded original time series which directly mean that finding low-rank approximations of the Hankel matrix has sensible meaning in yielding a time-series which approximates the original one, but is of smaller complexity.

In this paper we offer a review and bibliography of work on Hankel low-rank approximation and completion, with particular emphasis on how this methodology can be used for time series analysis and forecasting. The structure of the paper is as follows. Section \ref{sec:hlra} introduces the problem of Hankel low-rank approximation, describing possible choices of norm to define the optimization problem, as well as its possible formulation via a vector or matrix form. Unstructured low-rank approximation and its solution via the classic Eckart-Young-Mirsky-Schmidt theorem is described as well as related problems. The important relationship between low-rank Hankel matrices and linear recurrence sequences is provided in Section~\ref{sec:lin_rec}. Algorithms to approximate a solution of the Hankel low-rank approximation problem are described in Section~\ref{sec:algs}, and the problem of forecasting is framed as one of low-rank matrix completion in Section~\ref{sec:conv_rel} where use of the nuclear norm as a convex relaxation of the rank constraint is explained. The paper closes with some numerical examples in Section~\ref{sec:ex}.

\section{Hankel low-rank approximation}\label{sec:hlra}
\subsection{Problem statement}
Let
\begin{equation}\label{eq:p_vector}
\bp=(p_{1},p_{2},\ldots,p_{\len}),
\end{equation}
be a time series of length $N>1$, taking real or complex values.
We mostly focus on real-valued time series in this paper since they more readily appear in statistical application.

Most low-rank based methods of time series analysis embed the vector $\bp$ into a Hankel matrix.  Let $L$, $K$ be integers such that $L+K-1 = N$. For simplicity of notation we will assume $L\le K$ throughout. The $L \times K$ Hankel matrix $\mathcal{S}(\bp)$ parameterized by the vector $\bp$ is defined as
\begin{equation}\label{eq:struct_hankel_general}
\mathcal{S}(\bp) = \mathcal{H}_{L}(\bp)  =
\begin{pmatrix}
p_1 & p_2 & p_3 & \cdots & \cdots & p_K \\
p_2 & p_3 & \iddots & \cdots& \cdots & \vdots \\
p_3 & \iddots & \iddots & \iddots &\iddots & \vdots \\
\vdots &  \iddots &\iddots & \iddots&  \iddots & p_{N-1} \\
p_L & \cdots & \cdots & \cdots & p_{N-1} & p_N \\
\end{pmatrix}.
\end{equation}
The parameter $L$ is specified by the user, and it referred to as the window length in techniques such as singular spectrum analysis (SSA) \cite{golyandina2020singular}.
Its choice depends on the low-rank approximation (LRA) problem under consideration.
The most recent book on singular spectrum analysis is \cite{golyandina2020singular} and advice on the choice of $L$ for SSA is given in \cite{atikur2013note} for this technique, with some more general discussion available in \cite{golyandina2010choice}.
Commonly $L$ is chosen so that the resulting Hankel matrix is (close to) square \cite{markovsky2019book}.
Another possible strategy is to choose the window length as small as possible in order to make the matrix ``fat''.

Let the given vector of $N$ observations (time series) be
\[
\bp_{0}=(p_{0,1},p_{0,2},\ldots,p_{0,N})
\]
and $\bX_0 =  S(\bp_0)$ be the associated Hankel structured matrix.
Low-rank based methods of approximating/forecasting a given time series aim to obtain a  so-called structured low-rank approximation (SLRA) of the matrix  $\bX_0$.
The problem of SLRA in the literature is typically formulated in two ways: vector and matrix form, which informally speaking refers to the domain in which the SLRA optimization problem is defined.
In fact, these formulations can be shown to be equivalent, but both of them are used in the literature.
We begin our paper with the vector form, which seems to be more suited to time series analysis.

\paragraph*{Vector (parameter) form}
For a given time series $\bp_{0}$ and fixed rank $r$, the SLRA problem is to find the closest time series $\widehat{\bp}$, such that the corresponding structured matrix $\mathcal{S}(\widehat{\bp})$ is low-rank:
\begin{equation}\label{eq:SLRA}
\min_{\widehat{\bp}} \|\widehat{\bp}-\bp_{0} \| \quad \sto \quad \rank{\mathcal{S}(\widehat{\bp})} \le r,
\end{equation}
where $\|\cdot\|$ is  a vector norm (or, in general, an extended semi-norm).
The SLRA problem thus can be interpreted as an approximation of a time series by a simpler time series having some structure.
The class of time series corresponding to low-rank Hankel matrices is described in section~\ref{sec:lin_rec}.

The norm $\|\cdot\|$ measures the distance between the time series and its approximation.
Typical choices of norm include:
\begin{itemize}
\item (Euclidean) $2$-norm $\|\cdot\|_2$, which puts equal weights to all observations;
\item elementwise weighted extended semi-norm
\[
\|\bp\|_{\bw} =  \sqrt{\sum\limits_{k} w_{k} p_{k}^2},
\]
defined by a nonnegative vector of weights $w_k \in [0;+\infty]$; more precisely,
\begin{itemize}
\item $w_k \in (0,+\infty)$ quantifies importance of  observations;
\item $w_k = \infty$ corresponds to fixed values (i.e., an additional constraint $\widehat{p}_k = p_{0,k}$ in \eqref{eq:SLRA});
\item $w_k = 0$ implies that the value $ p_{0,k}$ is not important, and therefore can encode missing values;
\end{itemize}
in particular, if all $w_k \in (0,+\infty)$, then $\|\bp\|_{\bw}$ is a norm.
\item  General weighted norm
\[
\|p\|_{\bW} := \sqrt{ p^{\top} \bW  p} ,
\]
where $\bW$ is a symmetric positive semidefinite matrix. In particular,
\begin{itemize}
\item the choice $\bW = \diag(w_1, \ldots, w_N)$ corresponds to the elementwise weighted extended norm ($w_k \in (0,+\infty)$);
\item a general matrix $\bW$ can be taken as an inverse covariance matrix of the noise (if an additive noise model is assumed).
\end{itemize}
\end{itemize}
Other types of norms can be considered (see, e.g., examples in the matrix formulation described next), but the SLRA problem becomes much more difficult in that case.

\paragraph*{Matrix form}
It is the matrix form that was used in many of the original papers and classical algorithms on SLRA \cite{chu2003slra}.
For a given  matrix $\bX_0 = \mathcal{S}(\bp_0) \in \mathscr{S}$ belonging to the set of structured matrices $\mathscr{S} = \{  \mathcal{S}(\bp) \,|\, \bp \in \mathbb{R}^{N}\}$ and a chosen rank $r$ we seek its best rank-$r$ structured approximation, i.e.
\begin{equation}\label{eq:SLRA_mat_form}
\min_{\widehat{\bX}} \|\widehat{\bX} - \bX_0\| \quad \sto \quad \rank{\widehat{\bX}} \le r \text{ and }  \widehat{\bX} \in \mathscr{S},
\end{equation}
where $\|\cdot\|$ is some matrix norm or, in general,  an extended semi-norm.
Examples  of  common choices of norms are:
\begin{itemize}
\item the classic Frobenius norm
\[
\| \bX \|_{F} =  \sqrt{\sum\limits_{i,j} X_{ij}^2};
\]
\item the elementwise weighted Frobenius norm
\[
\| \bX \|_{\bZ}  =  \sqrt{ \sum\limits_{i,j} Z_{ij} X_{ij}^2},
\]
where $\bZ$ is a matrix of  nonnegative weights encoding the importance of the elements in the approximation.
As in the vector formulation described above, the weights can take the values $\infty$ and $0$, as to encode fixed or missing values respectively, which makes $\|\cdot\|_{\bZ}$ an (extended) semi-norm.

\item so-called (Q,R)-norm \cite{gillard2016weighted}:
\[
\| \bX \|_{\bQ,\bR} = \sqrt{\trace (\bQ\bX\bR\bX^{T})},
\]
where $\bQ$ and $\bR$ are  positive definite matrices.
\item Generalized weighted Frobenius norm, defined by a symmetric positive definite matrix $\bW_{m} \in \mathbb{R}^{LK\times LK}$,
\begin{equation}\label{eq:general_weight_matrix}
\| \bX \|_{\bW_{m}} = \sqrt{ (\vecm{\bX})^{T}  \bW_{m} \vecm{\bX}}.
\end{equation}
In fact, the general weighted Frobenius norm generalizes the previous two norms, with
\begin{itemize}
\item  $\bW_{m} = \diag(\vecm{\bZ})$ for the elementwise weighted norm;
\item $\bW_{m} = \bR \kron \bQ$ for the (Q,R)-norm, with $\kron$ being the Kronecker product;
\item in particular, a norm that is both an elementwise and a (Q,R)-norm corresponds to a rank-one matrix $\bZ$ and diagonal matrices $\bQ$ and $\bR$:
\[
\bZ = \bq \br^{T}, \quad \bQ = \diag(\bq) \quad\mbox{and}\quad \bR =   \diag(\br)\, .
\]
\end{itemize}

\item Other norms can be considered, e.g., the spectral norm, but this results in a much more difficult problem and very few results and practical algorithms are available \cite{knirsch2021optimal,Antoulas1997,rump2003structured,floater2021best}. The spectral norm is the dual of the nuclear norm \cite{fazel2001rank}, and so is often studied when considering properties of nuclear norm minimization. The spectral norm underpins the work of Adamyan, Arov and Krein which has
come to be known as AAK theory \cite{adamjan1971analytic}: a series of results connecting the theory of complex functions on the unit circle to infinite Hankel matrices. Optimal approximations of Hankel matrices in the spectral norm is  gaining traction in developing modern model reduction methods necessary for recurrent neural networks \cite{balle2021optimal}.
\end{itemize}


\paragraph*{General remarks}
\begin{itemize}
\item For the Frobenius norm the unstructured LRA (without the structure constraint) is solved by the singular value decomposition (see the following subsection).

\item The SLRA problem is formulated exactly in the same manner as in  for general (affine) matrix structures; in this paper, however, we focus on the Hankel SLRA, but discuss related structures in section~\ref{sec:other_structures}.

\item
In the general case, the  SLRA problem is a difficult NP-hard global optimization problem \ \cite{gillis2011low}. For Hankel SLRA (as shown in \cite[Corollary 3.9]{ottaviani2014exact} for the case of Frobenius norm), the number of  stationary points grows polynomially in $\len$ and exponentially in $r$. The minimizers possess narrow regions of attraction and the effect of observing noisy observations further compounds the difficulty, dampening and shifting the location of the global minimum \cite{gillard2013optimization}.

\item The two formulations (matrix and vector), in general, are equivalent in the case of the general weighted semi-norm \cite{ishteva2014factorization}. We review this correspondence in  section~\ref{sec:parameterizations}.

\item All the definitions and problem formulations given thus far are naturally extended to the complex case, whereupon we replace transpose by the Hermitian transpose, and squares (e.g., $p_i^2$) by the squares of the absolute values (e.g., $|p_i|^2$ in the definitions of the norms.
\end{itemize}

\subsection{Unstructured low-rank approximation and global solutions}
Before discussing properties of  SLRA, we recall results on unstructured low-rank approximation problems.
The first result is the classic Eckart-Young-Mirsky-Schmidt theorem  (see e.g., \cite{golub2012matrix}) which treats the case of unitarily invariant norms, and the   Frobenius norm in particular.

\begin{theorem}\label{thm:svd}
Let $\bA \in \mathbb{R}^{L \times K}$ be a given matrix ($L \leq K$) with singular value decomposition (SVD) given by $\bA=\bU \bSigma \bV^T$, where $\bU \in \mathbb{R}^{L \times L}$, $\bV \in \mathbb{R}^{K \times K}$ and $\bSigma \in \mathbb{R}^{L \times K}$,  with elements $\bSigma_{i,i}=\sigma_{i}$ being the singular values $\sigma_{1} \geq \sigma_{2} \geq \cdots \geq \sigma_{L} \ge 0$. Let $\bSigma^{(r)} \in \mathbb{R}^{L \times K}$ with elements
$$
\bSigma^{(r)}_{i,j}=\left\{
                \begin{array}{ll}
                  \sigma_{i}, & 1 \leq i \leq r \textrm{ and } j=i \\
                  0, & \hbox{otherwise.}
                \end{array}
              \right.
$$
Then $$
\bA^{(r)}:=\bU \bSigma^{(r)} \bV^{T}=\argmin_{\bB: \rank(\bB)\leq r}\| \bA-\bB \|_{\diamond}^{2}\, ,
$$
where $\| \cdot \|_{\diamond}$ is an arbitrary unitarily invariant norm.
\end{theorem}

\begin{remark}
The best  low-rank approximation is not uniquely defined in the case $\sigma_r = 
\sigma_{r+1}$, which is due to nonuniqueness of the singular vectors corresponding to repeating singular values. In such a case, all the possible SVDs give the set of all best rank-$r$ approximations, according to Theorem~\ref{thm:svd}.
\end{remark}

An extension of the above theorem to the (Q,R)-norms is as follows.
\begin{lemma}\label{lem:posdef}
Let $\bA \in \mathbb{R}^{L \times K}$, $\bB \in \mathbb{R}^{L \times K}$, $\bQ$ and $\bR$ be positive definite matrices of dimension $L \times L$ and $K \times K$ respectively. Then
$$
\| \bA-\bB\|_{\bQ,\bR}^2=\| \bQ^{1/2}(\bA-\bB)\bR^{1/2} \|_{F}^{2}\, .
$$
\end{lemma}
 A proof of Lemma \ref{lem:posdef} and a discussion of its importance in relating problems of time-series analysis and low-rank approximation is given in \cite[Theorem 4.12]{markovsky2019book} (see also \cite{gillard2016weighted} for a detailed study). Note that as the matrices $\bQ$ and $\bR$ are assumed to be positive definite, then $\bQ^{1/2}$ and $\bR^{1/2}$ exist, as do their inverses.

\begin{theorem}
Let $\bA \in \mathbb{R}^{L \times K}$ be a given matrix ($L \leq K$) with SVD of $\tilde{A}:= \bQ^{1/2} \bA\bR^{1/2} = \bU {\Sigma}\bV$ as described in Theorem \ref{thm:svd}, and $\bQ$ and $\bR$ be given positive definite matrices of dimension  $L \times L$ and $K \times K$ respectively. Then
\begin{equation*}
\bA^{(r)}_{\bQ,\bR}:=(\bQ^{1/2})^{-1}\bU {\Sigma}^{(r)}\bV (\bR^{1/2})^{-1}=\argmin_{\bB: \rank(\bB) \leq r} \|\bA-\bB \|^{2}_{\bQ,\bR}
\end{equation*}
 \end{theorem}

These results do not consider structured approximations, in that there are no constraints on the elements of the matrix approximation. However these results are necessary for SLRA algorithms to obtain either: an initial point to begin iterations of some numerical algorithm and/or to provide projections to the `closest' matrix of rank $\leq r$. Indeed we will use these results to introduce some algorithms which aim to approximate a solution of \eqref{eq:SLRA} in Section \ref{sec:algs}.

\subsection{Related problems}
In some literature, instead of low-rank matrix approximation, the problem is formulated as one of rank minimization.
For example, in the vector formulation, for a given approximation error $\varepsilon$,  we wish to minimize the rank
\begin{equation}\label{eq:rankmin}
\min_{\widehat{\bp}} \rank  \,\mathcal{S}(\bp)\,\,\,\,  {\rm subject\,\, to}\,\,\,\|\widehat{\bp}-\bp_{0}\| \leq \varepsilon.
\end{equation}
This formulation is dual to that of SLRA \eqref{eq:SLRA}, as both problems explore the same Pareto optimality front (see e.g., \cite{markovsky2008structured} or \cite{markovsky2019book} for a discussion on this topic).
In fact, if we are able to solve one of the problems we can find a solution of the other one.
One particular drawback of the formulation in \eqref{eq:rankmin} is that the solution is in most cases non-unique, and this is one reason why we prefer the SLRA formulation.

Another possible formulation is the ``regularized'' form,
\begin{equation}\label{eq:rankmin_reg}
\min_{\widehat{\bp}} \rank  \,\mathcal{S}(\bp) + \lambda \|\widehat{\bp}-\bp_{0}\|^2,
\end{equation}
but this formulation is much less intuitive and makes it difficult to find the regularization parameter.
Yet another related problem is trace regression \cite{rohde2011estimation} (or matrix sensing), but we leave it out of the scope of this paper.

\subsection{Parameterizations and norms}\label{sec:parameterizations}
The Hankel matrix structure belongs to the class of affine matrix structures \cite{markovsky2019book}, i.e. the matrix structures which can be parameterized as:
\begin{equation}\label{eq:linstr}
\mathcal{S}(\bp)=\bS_{0}+\sum_{i=1}^{N}p_{i}\bS_{i}
\end{equation}
where $\bS_{i}$, $i \in \{ 0,1,\ldots, N\} \in \mathbb{R}^{L \times K}$ are given linearly independent basis matrices. In particular, for the Hankel matrix structure \eqref{eq:struct_hankel}, the  basis matrices in \eqref{eq:linstr} are given as $\bS_{0} = 0$,
\begin{equation}
\begin{split}
&\bS_1 =\begin{psmallmatrix}
1    & 0     & \cdots & 0 & 0 \\
0     & 0     & \iddots & 0 & 0     \\
\vdots  & \iddots  & \iddots & \iddots  &  \vdots \\
0 & 0 & \iddots &  0 & 0\\
0 & 0  &  \cdots & 0 & 0
\end{psmallmatrix},
\bS_2 =\begin{psmallmatrix}
0    & 1     & \cdots & 0 & 0 \\
1     & 0     & \iddots & 0 & 0     \\
\vdots  & \iddots  & \iddots & \iddots  &  \vdots \\
0 & 0 & \iddots &  0 & 0\\
0 & 0  &  \cdots & 0 & 0
\end{psmallmatrix},\ldots,\\
&\bS_{N-1} =\begin{psmallmatrix}
0    & 0     & \cdots & 0 & 0 \\
0     & 0     & \iddots & 0 & 0     \\
\vdots  & \iddots  & \iddots & \iddots  &  \vdots \\
0 & 0 & \iddots &  0 & 1\\
0 & 0  &  \cdots & 1 & 0
\end{psmallmatrix},
\bS_{N} =\begin{psmallmatrix}
0    & 0     & \cdots & 0 & 0 \\
0     & 0     & \iddots & 0 & 0     \\
\vdots  & \iddots  & \iddots & \iddots  &  \vdots \\
0 & 0 & \iddots &  0 & 0\\
0 & 0  &  \cdots & 0 & 1
\end{psmallmatrix}.
\end{split}
\end{equation}
The optimization problems \eqref{eq:SLRA} and \eqref{eq:SLRA_mat_form} are equivalent for particular choices of vector and matrix norm. For example
$$
\|\mathcal{S}(\bp)-\mathcal{S}(\bp_{0}) \|_{F}^{2}=\sum_{j=1}^{N}\kappa_{j}(p_{j}-p_{0,j})^2
$$
where the weights
\begin{equation}\label{eq:froweights}
\kappa_{j}=\left\{
  \begin{array}{ll}
    j, & \text{if } 1 \leq j < L, \\
    L, & \text{if } L \leq j \leq K, \\
    N-j+1, &\text{if } K < j \leq N.
  \end{array}
\right.
\end{equation}

This remark has a large implication on the use of SLRA for time series applications. Should the Frobenius norm be used in \eqref{eq:SLRA}, then the weights as given in \eqref{eq:froweights} are (possibly unintentionally) bestowed upon the time series $\bp_{0}$. Work to suggest solutions to rectify this problem is covered in \cite{zvonarev2017iterative}, whereby given desired weights $w_{j}$ methods to obtain matrices $\bQ$ and $\bR$ such that
$$
\| \mathcal{S}(\bp)-\mathcal{S}(\bp_{0})\|_{\bQ,\bR}^2 \approx \sum_{j=1}^{N}w_{j}(p_{j}-p_{0,j})^2\, ,
$$
are described. For certain cases of elementwise weights $w_{j}$ it is possible to find exact equivalencies between $\| \mathcal{S}(\bp)-\mathcal{S}(\bp_{0})\|_{\bQ,\bR}^2$ and $\sum_{j=1}^{N}w_{j}(p_{j}-p_{0,j})^2$, but in general only approximations are possible.
However, if the general weight matrices are allowed, as in \eqref{eq:general_weight_matrix}, then the vector and matrix formulations  \eqref{eq:SLRA} and \eqref{eq:SLRA_mat_form} can be shown to be equivalent \cite{ishteva2014factorization}.


\section{Low-rank Hankel matrices and linear recurrence sequences}\label{sec:lin_rec}
Rank-deficient Hankel matrices are well-studied  \cite{heinig1984algebraic}, in particular, they are tightly linked to linear recurrence sequences.
Such sequences, can be characterized by a so-called canonical representation of sums of products of polynomial, exponential and cosine functions, which is a very useful class of time series.
In this section, we provide a summary of these connections, together with references to the relevant literature.

\subsection{Rank properties of Hankel matrices and time series of finite rank}\label{sec:rankHankel}
The algebraic theory of Hankel matrices goes back to Sylvester \cite{sylvester1886binary}, who studied them in the context of decompositions of binary forms, and Prony \cite{prony1975}.
The most complete reference is the book of Heing \& Rost \cite[Ch.5 and 8]{heinig1984algebraic}, however, this reference was not known very well until recently for historical reasons.
While not pretending to provide a comprehensive treatment of the subject, we briefly review the key points that are particularly useful in the context of time series analysis, following the details provided in
\cite{usevich2010signal} and \cite{gillard2018structured}.

\begin{figure}
  \centering
  \includegraphics[width=0.35\textwidth]{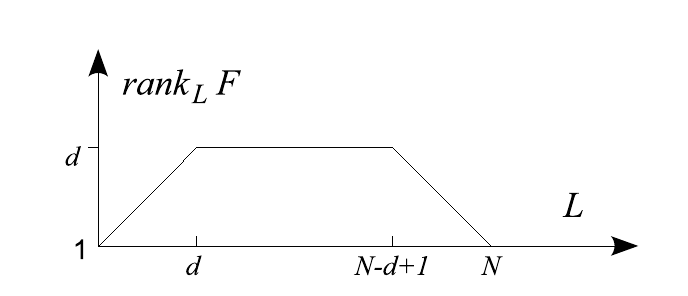}
  \caption{Rank of the Hankel matrix $\mathcal{H}_{L}$ for varying $L$}\label{fig:ranks}
\end{figure}

The key property is to that of how the rank varies with respect to the window length.
Assume that we fix the length of the time series $N$ and we vary the window length $L$ from $1$ to $N$ (allowing  in this case that $L$ may be bigger than $K$).
\begin{lemma}[{\cite[Corollary 5.2]{heinig1984algebraic}}]\label{lem:hankel_rank}
For any sequence $\bp \in \mathbb{C}^{N}$ there exists a number $d \le \lfloor \frac{N+1}{2} \rfloor$ such that
\[
\rank  \mathcal{H}_{L}(\bp) = \min(L,K,d).
\]
\end{lemma}
Figure~\ref{fig:ranks} illustrates the behaviour of ranks as described in Lemma~\ref{lem:hankel_rank}.
As it can be seen from Fig.~\ref{fig:ranks}, the rank of the matrix is maximal, excepting a ``plateau'' at $d$.
The case $d = \frac{N+1}{2}$ (for time series of odd length) is special, as in this case the Hankel matrix is full rank for any value of $L$ (and is square nonsingular for $L = \frac{N+1}{2}$ in particular).
For other cases ($d < \frac{N+1}{2}$), the matrices have a very particular structure; we will call such $\bp$ \emph{time series of finite rank}, as in the SSA literature \cite{golyandina2001analysis}.

There is a following important consequence of Lemma~\ref{lem:hankel_rank}: if a Hankel matrix is rank-deficient (i.e., $\rank  \mathcal{H}_{L}(\bp) \le r < \min(L,K)$, and this is a time series of finite rank),  then $\rank  \mathcal{H}_{L}(\bp)  = \rank  \mathcal{H}_{r+1}(\bp)$, where $\mathcal{H}_{r+1}(\bp)$ is a Hankel matrix of built from the same vector but with $r+1$ rows.
This leads to the following corollary.
\begin{corollary}\label{cor:shape_hankel_fat}
In the vector formulation 
\eqref{eq:SLRA}, the SLRA problem with $r < L$ and $\mathcal{S} = \mathcal{H}_{L}$ is equivalent to the SLRA problem for the rank $r$ and structure $\mathcal{S} = \mathcal{H}_{r+1}$, for the same choice of the (semi-)norm $\|\cdot\|$ (i.e., the search space and the optimal solutions of both problems coincide).
\end{corollary}

Thanks to the Corollary~\ref{cor:shape_hankel_fat} mentioned above, some authors (e.g., in the work of Markovsky and coauthors, \cite{markovsky2019book}) argue that it makes no sense to consider  Hankel SLRA with other window length than $L=r+1$ (and therefore ignore other window lengths).
This, however, leads to a confusion in the literature as it leaves an impression that the approach of \cite{markovsky2019book} cannot handle other cases than rank reduction by $1$.

\subsection{Linear recurrence relations (LRR)}
Low-rank Hankel matrices are tightly linked to linear recurrence relations.
We say that a (real- or complex-valued) time series $\bp \in \mathbb{C}^{N}$ satisfies a linear recurrence relation (LRR) of order $r$ if
there exists a non-zero vector $\boldsymbol{\theta} = (\theta_0, \ldots, \theta_r) \in \mathbb{C}^{r+1}\setminus\{0\}$ such that
\begin{equation}\label{eq:lin_rec_nonminimal}
\theta_0 p_k  + \ldots + \theta_r p_{k+r} = 0, \quad \mbox{for all } k = 1,..., N-r.
\end{equation}
For real-valued $\bp$, the coefficients of the linear recurrence relation can be also be chosen to be real-valued.

Next, we remark that the linear recurrence property is related to rank-deficiency of a Hankel matrix.
\begin{remark}\label{rem:lrr_and_rank}
In the matrix form, we can rewrite \eqref{eq:lin_rec_nonminimal} as
\begin{equation}\label{eq:kernel_rank_constraint}
\begin{pmatrix}\theta_0 &  \cdots & \theta_r \end{pmatrix}  \mathcal{H}_{r+1} (\bp) = 0.
\end{equation}
Here, as in Lemma~\ref{lem:hankel_rank}, the Hankel matrix may have more rows than columns.
This implies that if a Hankel matrix is rank-deficient, then the underlying time series satisfies an LRR.
\end{remark}
Note that, in general, $\theta_0$ and $\theta_r$ may be zero.
If the last coefficient is non-zero, then the linear recurrence \eqref{eq:lin_rec_nonminimal} can be represented in a slightly different form, which is more common in time series analysis (e.g., \cite{golyandina2001analysis}).
\begin{remark}
If\footnote{The case $\theta_r \neq 0$ is called ``continuable time series'' in the SSA theory \cite{golyandina2001analysis}. Indeed, from Remark~\ref{rem:lrr_and_rank}, the case $\theta_r \neq 0$ corresponds to the case when $\be_{L} \not\in \spanm( \mathcal{H}_{L} )$, which is a necessary and sufficient condition of continuability in \cite{golyandina2001analysis}.} $\theta_r \neq 0$, then we can define $a_k = -\theta_k/\theta_r$, and we can rewrite \eqref{eq:lin_rec_nonminimal} as
\begin{equation}\label{eq:LRR}
p_{j}=a_{r-1}p_{j-1}+a_{r-2}p_{j-2}+\cdots+a_{0}p_{j-r}.
\end{equation}
\end{remark}
The LRR \eqref{eq:LRR}  (linear recurrent formula in the terminology of \cite{golyandina2001analysis}) gives a way to recursively compute the next values of the time series based on the previous ones.
This serves as a base for time series forecasting, explained in section~\ref{eq:slrmc_solution}.

\subsection{Linear recurrence relations and sums of exponentials}
In this subsection, we explain the connection between time series of finite rank and sums of (modulated) complex exponentials.
This connection is well-known for infinite sequences (or infinite time series), see, for example, \cite{hall1967combinatorial}, but in this and subsequent sections, we show what happens for finite-length time series as well.

The key tool is the \emph{characteristic polynomial} of the LRR, which is
\[
\theta(z) = \theta_r z^r +\theta_{r-1}z^{r-1}+\cdots+\theta_{1}z+\theta_{0}\, .
\]
Assuming that the first and last coefficients are nonzero, the characteristic polynomial can be factorized thanks to the fundamental theorem of algebra,
\begin{equation}\label{eq:char_poly_nonminimal}
\theta(z) = (z-\lambda_{1})^{\nu_{1}}(z-\lambda_{2})^{\nu_{2}}\cdots(z-\lambda_{s})^{\nu_{s}},
\end{equation}
with $\nu_{1} + \ldots + \nu_{s} = r$.
Moreover, the form of the time series can be determined from its roots, as shown by the following lemma.
\begin{theorem}[{\cite[Thm. 3.1.1]{hall1967combinatorial}}]\label{thm:complex}
A complex-valued time series\footnote{In \cite{hall1967combinatorial} Theorem~\ref{thm:complex} was formulated for  infinite-length time series, but it holds for finite time series by restriction.} $\bp$ satisfies a linear recurrence with $\theta_0 \neq 0$ and $\theta_s \neq 0$ if and only if it can be expressed as
\begin{equation}\label{eq:complex}
p_{k}=\sum_{j=1}^{s}P_{j}(k)\lambda_{j}^{k}, \,\,\, k=1,2,\ldots
\end{equation}
where $P_{j}(k)$ are complex polynomials of degrees at most $\nu_{j}-1$ and $\lambda_{j} \in \mathbb{C}\setminus \{ 0\}$ are the roots of the polynomial in \eqref{eq:char_poly_nonminimal}.
The coefficients of the polynomials in \eqref{eq:complex} are determined by the first $s$ values of the time series.
\end{theorem}

Roughly speaking, Theorem~\ref{thm:complex} implies that the Hankel SLRA problem becomes a problem of approximating a given time series by time series of the form \eqref{eq:complex}.
Indeed, thanks to Remark~\ref{rem:lrr_and_rank},  the Hankel matrix is rank-deficient if and only if satisfies a linear recurrence relation, and therefore in most cases its form is given by Theorem~\ref{thm:complex}. (The general case will be discussed at the end of this subsection.)
In the case when the characteristic polynomial has simple roots ($s=r$ and $\nu_j = 1$ in \eqref{eq:char_poly_nonminimal}), we have that
\begin{equation}\label{eq:complex_simple}
p_{k}=\sum_{j=1}^{r}c_k\lambda_{j}^{k}, \,\,\, k=1,2,\ldots
\end{equation}
which is known as a sum-of-exponential model in signal processing, see also subsection~\ref{sec:other_structures}.

In time series analysis, the real-valued time series are more common.
Therefore, we formulate a corollary of Theorem~\ref{thm:complex}.

\begin{corollary}\label{cor:real_lin_rec}
A real-valued time series satisfies a linear recurrence \eqref{eq:lin_rec_nonminimal} with real coefficients $\boldsymbol{\theta} \in \mathbb{R}^{r+1}$, $\theta_0,\theta_r \neq 0$ if and only if it can be written  as the following sum \cite{lemmerling2001analysis}:
\begin{equation}\label{eq:damped}
p_{k}=\sum_{j=1}^{q}P_j(k) \exp(d_{k}k)\sin(2\pi\omega_{j}k+\phi_{j}),\,\,\, k=1,\ldots,N
\end{equation}
where $\{P_j(k)\}$ are real polynomials (which can be just constant amplitudes), $\{d_{l}\}$ are the dampings, $\{\omega_{l}\}$ are the frequencies and $\{\phi_{l}\}$ are the phases, to use terminology common in signal processing.
\end{corollary}
Corollary~\ref{cor:real_lin_rec} follows from the fact that if the coefficients of the linear recurrence are real, then the non-real  roots of the characteristic polynomial  must come in complex conjugate pairs.
The number of terms $q$ is equal to the number of real roots of $\theta(z)$ plus the number of complex conjugate pairs (i.e., for the case of $q=s/2$ in \eqref{eq:complex_simple} if there are no real roots).
More details on Corollary~\ref{cor:real_lin_rec} can be found, e.g., in  \cite{golyandina2001analysis, usevich2010signal, gillard2018structured,zvonarev2021low}.

\subsection{Minimal linear recurrences and complete characterisation of finite-length case}
While Theorem~\ref{thm:complex} and Corollary~\ref{cor:real_lin_rec}  gives us a form of the time series for several cases, it does not give any information whether this representation is unique.
Another drawback is that it does not give a complete characterization, as it does not treat the cases when $\theta_0 = 0$ and $\theta_r = 0$.

Let us first discuss the uniqueness issue.
In fact, we can define the linear recurrence relation of minimal length \cite{649704}, which we can find as a corollary of Lemma~\ref{lem:hankel_rank} and  Remark~\ref{rem:lrr_and_rank}
\begin{lemma}\label{eq:min_linrec}
For a time series of finite rank $d$, the minimal length linear recurrence relation  is of  length $d$ with coefficients $(q_0, \ldots, q_{d})$ . It exists and is unique up to multiplication by a constant factor, and can be found from the left kernel of the corresponding Hankel matrix
\[
(q_0, \ldots, q_d)  \mathcal{H}_{d+1} (\bp) = 0.
\]
\end{lemma}
Indeed, a shorter linear recurrence is not possible (which would imply $\mathcal{H}_d$ to be rank-deficient),
and the left kernel of $\mathcal{H}_{d+1} (\bp)$ is generated by a single vector, which is defined uniquely up to normalization and sign.
Such a vector $(q_0, \ldots, q_d) $ is related to the first characteristic polynomial in the terminology of \cite[]{heinig1984algebraic}.
For a range of lengths, all LRRs can be obtained from $\bq$, as explained in section~\ref{sec:left_kernel}.

Moreover, the time series  $\bp$ possesses a unique minimal representation similar to  \eqref{eq:complex_simple}.
In what follows, we formulate the general case in the following theorem (as to allow the coefficients $q_0$ and $q_d$ to be zero).
Let $\lambda_1, \ldots,\lambda_s \in  \mathbb{C} \setminus \{0\}$,  $\nu_1,\ldots,\nu_s$ be some positive integers, and $\nu_0, \nu_{\infty} \in \NN$ be nonnegative integers such that 
\begin{equation}\label{eq:multiplicities_sum}
\nu_0 + \nu_1 + \cdots + \nu_{s} + \nu_{\infty} = d.
\end{equation}
Consider a signal $\bp \in \CC^{N}$ such that
\begin{equation}\label{eq:canonical}
p_k = \underbrace{\sum\limits_{j=1}^{\nu_0} a_j \delta_{j-1}(k) + \sum\limits_{l=1}^{\nu_{\infty}} b_l \delta_{N-l}(k)}_{\text{transient terms (beginning and end)}} + \sum\limits_{j=1}^{s} P_j(k) \lambda_j^{k},
\end{equation}
where
\begin{enumerate}
\item each $P_j(k)$ is a polynomial of degree $\nu_k-1$;
\item $\delta_x(k)$ is the (Kronecker) delta-function:
\[
\delta_x(k) =
\begin{cases}
1, & x = k, \\
0, & x \neq k;
\end{cases}
\]
\item and $a_{\nu_0} \neq 0$,  $b_{\nu_{\infty}} \neq 0$.
\end{enumerate}
Finally, let $(q_0, \ldots, q_{d})$ be coefficients of the characteristic polynomial of the minimal linear recurrence
\begin{equation}\label{eq:char_poly_general}
q(z) = q_0 + q_1 z +  \cdots + q_{d} z^d = z^{\nu_0}(z-\lambda_1)^{\nu_1} \cdots(z-\lambda_s)^{\nu_s};
\end{equation}
note $\nu_0$ denotes the multiplicity of the root $0$ (if present) and $\nu_{\infty}$ encodes the number\footnote{This corresponds to the multiplicity of $\infty$ root in the terminology of \cite{heinig1984algebraic}.} of zero leading coefficients ($q_d$, $q_{d-1}$, ...).
The $q(z)$ is precisely the characteristic polynomial of the minimal linear recurrence relation $\bq$ in Lemma~\ref{eq:min_linrec}.

\begin{theorem}[Special case of {\cite[Thm. 8.2]{heinig1984algebraic}}]\label{thm:finite_rank_general}
The time series $\bp$ is of finite rank $d$ if and only if the signal $\bp$ has the form \eqref{eq:canonical}.
Moreover, the representation \eqref{eq:canonical} with $d$ satisfying \eqref{eq:multiplicities_sum}  (called the \emph{canonical representation} in \cite{heinig1984algebraic}), is unique.
\end{theorem}
For real-valued time series, a corollary of Theorem~\ref{thm:finite_rank_general} can be derived giving a complete characterisation of real values time series of finite rank in terms of their canonical representation,
in a similar manner to Corollary~\ref{cor:real_lin_rec}.

\subsection{Non-minimal linear recurrences and left kernel of Hankel matrices}\label{sec:left_kernel}
In this section, we give some details about the way the other linear recurrences \eqref{eq:lin_rec_nonminimal} (for $r > d$) can be obtained from the minimal one.
The following proposition gives a complete characterisation for a wide range of values of $r$. 
\begin{proposition}[{Special case of \cite[Theorem 5.1]{heinig1984algebraic}}]
Let $\bp$ be a complex-valued time series of finite rank $d$ and $\bq \in \CC^{d+1}$ be the coefficients of its minimal LRR. Then for any $r \in \{d,\ldots, N-d\}$ all linear recurrence relations $\begin{pmatrix}\theta_0 &  \cdots & \theta_r \end{pmatrix}$ of length $r$ satisfied by $\bp$ must have the form
\begin{equation}\label{eq:low-order-lrf}
\begin{pmatrix}\theta_0 \\  \vdots \\ \theta_r \end{pmatrix} = \underbrace{\begin{pmatrix}
q_{0} &    0    &     \cdots & 0        \\ 
q_{1} & q_{0} & \ddots &\vdots  \\
\vdots     & q_1 & \ddots &   0         \\
q_{d} &   \vdots  &\ddots   & q_{0}\\
0 & q_{d}  & & q_{1}   \\
\vdots & \ddots & \ddots  & \vdots  \\
0 & \cdots  &  0      & q_{d}     
\end{pmatrix} }_{\mathbf{M}_{r-d}(\mathbf{q})}
\begin{pmatrix}b_0  \\ b_1\\  \vdots \\ b_{r-d} \end{pmatrix},
\end{equation}
where $b_k$ are some complex coefficients (in the case of real-valued $\bp$, the coefficients $b_k$ can be chosen to be real as well).
\end{proposition}
In fact, \eqref{eq:low-order-lrf} can be conveniently interpreted in terms of the  {characteristic polynomial} of the linear recurrence.
Indeed, the equation \eqref{eq:low-order-lrf} is equivalent to 
\[
\theta(z)  = q(z) c(z);
\]
this is why $\mathbf{M}_{r-d}(\mathbf{q})$ is often called a multiplication matrix, because it represents  multiplication by the polynomial $q(z)$.
Thanks to Remark~\ref{rem:lrr_and_rank}, such linear recurrences correspond to all the vectors in the left kernel of the matrix $\mathcal{H}_{r+1} (\bp)$ (i.e., the set of all vectors such that $\boldsymbol{\theta}^{\top}\mathcal{H}_{r+1} (\bp) = 0$). 
Thus, for $d \le r \le  N-d$, the span of columns of the matrix ${\mathbf{M}_{r-d}(\mathbf{q})}$ in \eqref{eq:low-order-lrf} coincides with the left kernel of the matrix $\mathcal{H}_{r+1} (\bp)$ (which is a linear subspace of $\mathbb{C}^{r+1}$ of dimension $r-d+1$).

Note that for $N-d < r \le N-1$, the situation is more complicated: a second  generating linear recurrence comes into play, as explained in \cite[Theorem 5.1, Definition 5.4]{heinig1984algebraic}. 
This is due to the fact that as $r > N-d$, the rank defect of $\mathcal{H}_{r+1} (\bp)$ grows twice faster as $r$ increases (see Fig.~\ref{fig:ranks}).
Note that in the case of time series of infinite length, the second polynomial does not appear,  and all possible linear recurrences are given by \eqref{eq:low-order-lrf}, see e.g.~\cite{hall1967combinatorial}.

Finally, we make a remark that there is a representation similar to \eqref{eq:low-order-lrf} for the tangent space to the manifold of time series of finite rank $d$ (of dimension $2d$).
In fact, the tangent space is spanned by the matrix ${\mathbf{M}_{N-1-2d}(\mathbf{q^2})}$, where $\mathbf{q^2}$ is the vector of coefficients of squared polynomial  $q^2(z)$ \cite[Theorem 2.3]{zvonarev2021low}.
Similar results are known in the field of algebraic geometry in the study of so-called secant varieties (see e.g. \cite[Theorem 3.9]{iarobbino1999power} or \cite[Theorem 2.1]{brambilla2008alexander}).

\subsection{Other uses of Hankel SLRA and other matrix structures}\label{sec:other_structures}
As already remarked, the SLRA problem \eqref{eq:SLRA} can be formulated in  the same manner for an arbitrary affine structure.
We list below some common problems and structures.
\paragraph*{Hankel structure in signal processing}
Fitting and analysing the sum-of-exponentials model is very common in signal processing and so there has been much attention on Hankel SLRA to obtain solutions in problems such as: spectral estimation \cite{razavilar1996spectral,stoica2005spectral}, sparse signal recovery \cite{Markovsky:2018}, and denoising particular types of signal \cite{bresler1986exact,condat2015cadzow}.
Other uses of the sum-of-exponentials model include: finite rate of innovation  \cite{dragotti2007fri}, shape-from-moments \cite{elad2004shape}, and superresolution  \cite{candes2013super} problems.
Note that in spectral estimation it is customary to restrict to unit-norm exponents $|\lambda_k|=1$, which requires an additional constraint on positive semidefiniteness of the corresponding Toeplitz matrix \cite{candes2013super}. 
Some literature exists on developing approximate Prony methods (see \cite{zhang2019optimal} and references therein), which are closely linked to Hankel SLRA.
The Hankel structure is a special case of the block-Hankel structure which is discussed later in this subsection.

Yet another and classical use of the model \eqref{eq:complex_simple} in pure mathematics appeared in the field of algebraic geometry and tensor decompositions.
In fact, \eqref{eq:complex_simple}  corresponds exactly to the decomposition of a $2\times \cdots \times 2$ symmetric tensor  into a sum of rank-one symmetric tensors \cite{sylvester1886binary,iarobbino1999power,comon1996sums}, see also \cite{ottaviani2014exact} (equivalently, as a decomposition of a binary form into sums of powers of linear forms, the so-called Waring problem).
This is why many aspects of the theory of Hankel matrices trace back to classic works in algebraic geometry and are in parallel studied in the works on tensor decompositions.
Low-rank approximation of Hankel matrices is also a useful tool for function approximation \cite{beylkin2005exponential,gonnet2013robust}.

\paragraph*{Block-Hankel and mosaic Hankel structure}
This structure often appears when we wish to model a multivariate time series, i.e., a collection of $M$ time series $\bp^{(1)}, \ldots, \bp^{(M)}$.
A typical way is to stack the Hankel matrices corresponding to different time series  into one structure
\begin{equation}\label{eq:striped_hankel}
\begin{bmatrix} \mathcal{H}_{L}(\bp^{(1)}) &  \cdots & \mathcal{H}_{L}(\bp^{(M)})  \end{bmatrix}.
\end{equation}
A low-rank approximation of such series is performed when we assume some common structure of time series \cite{papy2006common},
which is used, for example, in multivariate SSA \cite{golyandina2015jss}.
Another use of such matrices occurs when the different time series are inputs and outputs of    a dynamical systems, which allows for solving system identification problems, thanks to behavioural systems theory \cite{markovsky2006siam,markovsky2008structured,markovsky2019book,willems1986linear}.

If the time series are of the same length, then the matrix \eqref{eq:striped_hankel} can be reshaped to have a block-Hankel structure.
In the general case (of different lengths), \eqref{eq:striped_hankel} is a special case of so-called mosaic Hankel structure \cite{slra-software,heinig1995mosaic}.
The block-Hankel matrix structure has the form
 which in general has the following block form:
 \begin{equation}\label{eq:block_hankel}
\begin{pmatrix}
\mathbf{H}_1 &\mathbf{H}_2 &\mathbf{H}_3 & \cdots & \mathbf{H}_K \\
\mathbf{H}_2 & \mathbf{H}_3  &\mathbf{H}_4 & \cdots &\mathbf{H}_{K+1} \\
\vdots &  \vdots &\iddots & \iddots & \vdots \\
\mathbf{H}_L & \mathbf{H}_{L+1}  &\mathbf{H}_{L+2}  & \cdots  & \mathbf{H}_N \\
\end{pmatrix},
\end{equation}
where $\mathbf{H}_k$ are $p\times q$ Hankel matrices.
Low-rank approximations of block-Hankel matrices are widely used in linear system theory, where it is related to the problems of realization and model reduction for linear systems \cite{kung1978identification,ho1966,markovsky2008structured}.
There is also a link between block-Hankel matrices and the short-time Fourier transform \cite{usevich2018stft}.

The algebraic theory for finite-size block-Hankel matrices (similar to that given in Sections~\ref{sec:rankHankel}-\ref{sec:left_kernel})  is also linked to linear recurrences, however, several linear recurrences are needed in this case \cite{heinig1992block}, which makes the explicit parameterization more complicated \cite{feldmann1999parameterization}.

\paragraph*{Multilevel Hankel matrices}
Multilevel Hankel matrices (also called quasi-Hankel matrices  \cite{mourrain2000multivariate}) appear in problems dealing with higher-dimensional data arrays (e.g., images), where the model of sum-of-exponentials \cite{mourrain2018polynomial,yang2016vandermonde} is often applicable.
A special case (for two-dimensional data) concerns Hankel-block-Hankel matrices, which have the structure as in \eqref{eq:block_hankel}, but with block also having Hankel structure.
Some application areas include: magnetic resonance imaging \cite{shin2014mri},
seismic data analysis \cite{jia2016geophysics,oropeza2011seismic}, two-dimensional extensions of SSA \cite{golyandina20102dssa,golyandina2015jss,shlemov2014shaped}, and symmetric tensor decompositions in the general case  \cite{brachat2010tensors}.
A complete algebraic theory for finite-size  multilevel Hankel matrices is not available; the characterization is available in the infinite  case, where the multivariate linear recurrences are linked to polynomial ideals \cite{mourrain2018polynomial,kurakin1995linear}.

\paragraph*{Sylvester matrices}
The Sylvester matrix structure is composed of  (two or more) multiplication matrices (i.e., the matrix in \eqref{eq:low-order-lrf}) stacked together.
The rank of the Sylvester matrix is related to the degree of the greatest common divisor of the polynomials whose co-efficients populate the elements of the Sylvester matrix. Consequently, Sylvester low-rank approximation appears in the problem of finding approximate common divisors of polynomials \cite{chin1998agcd,usevich2017agcd}, and in some signal processing applications, e.g., blind deconvolution \cite{fazzi2018agcd}.


\section{Algorithms}\label{sec:algs}
\subsection{Alternating projections and  subspace-based approaches}
\subsubsection{Cadzow's iterations and Singular Spectrum Analysis}
One of the simplest (and earliest) algorithms to describe are the so-called Cadzow iterations \cite{cadzow1988signal} which can broadly be described as an alternating projections algorithm, which alternates approximations between the closest low-rank matrix and the closest Hankel matrices. Each of these projections are described next.

\paragraph*{Projection to the closest low-rank matrix}
Suppose we are given $\bA=\mathcal{S}(\bp_{0})$. Denote
\begin{equation*}
\pi_{(\bQ,\bR)}^{(r)}(\bA)=\argmin_{\bB:\rank(\bB) \leq r}\|\bA-\bB \|^{2}_{\bQ,\bR}\, .
\end{equation*}
The matrix $\pi_{(\bQ,\bR)}^{(r)}(\bA)$ is obtained from the SVD of $\bA$, as described in Theorem \ref{thm:svd}.

\paragraph*{Projection to the closest Hankel matrix}
For a given matrix $\bA=\mathcal{S}(\bp_{0})$ we denote
\begin{equation}\label{eq:projH}
\pi_{(\bQ,\bR)}^{\mathcal{H}}(\bA)=\argmin_{\bB:\bB \in \mathscr{S}} \|\bA-\bB \|^{2}_{\bQ,\bR}, 
\end{equation}
where $\mathscr{S}$ is the set of structured (Hankel) matrices.
In the special case when $\bQ$ and $\bR$ are identity matrices then \eqref{eq:projH} has an explicit solution. In this case, let $\bA^\mathcal{H}=\pi_{(\bI,\bI)}^{\mathcal{H}}(\bA)$ and the projection is given by averaging over the antidiagonals
$$
\bA^\mathcal{H}_{l,k}=\kappa_{l+k-1}^{-1}\sum_{l+k=l'+k'}\bA_{l',k'}\,
$$
with $\kappa_{j}$ as defined in \eqref{eq:froweights}.
The advantage of this case (Frobenius norm) that the projection can be efficiently implemented with the Fast Fourier Transform, which brings the complexity of the Cadzow iterations down to $O(r N\log N)$ (see e.g. \cite{korobeynikov2010computation,golyandina2015jss}); in addition, the whole Hankel matrix need not be stored.

In the general case,  
thanks to Lemma~\ref{lem:posdef},  the projection can be found as $\pi_{(\bQ,\bR)}^{\mathcal{H}}(\bA) =\mathcal{S}(\widetilde{\bp})$ with
\begin{equation}\label{eq:projH_modified}
 \widetilde{\bp} = \argmin_{\bp} \|\widetilde{\mathcal{S}}(\bp) - \bQ^{1/2}\bA\bR^{1/2} \|_F,
\end{equation}
where $\widetilde{\mathcal{S}}(\bp) = \bQ^{1/2}\mathcal{S}(\bp)\bR^{1/2}$ is a modified matrix structure, which is still linear.
The minimizer of \eqref{eq:projH_modified} has an explicit solution as shown in \cite[Lemma 1]{ishteva2014factorization}. \color{black}

\paragraph*{Alternating projections}
Using the projections $\pi_{(\bQ,\bR)}^{(r)}$ and $\pi_{(\bQ,\bR)}^{\mathcal{H}}$ then we can construct the following alternating projections algorithm.

Set $\bA^{(0)}=\bA=\mathcal{S}(\bp_{0})$. For $i=1,2,\ldots$
\begin{equation}\label{eq:cad}
\bA^{(i)}=\pi_{(\bQ,\bR)}^{\mathcal{H}}\left(  \pi_{(\bQ,\bR)}^{(r)}(\bA^{(i-1)})  \right)
\end{equation}
If $\bQ$ and $\bR$ are identity matrices then iterations of \eqref{eq:cad} are known as Cadzow iterations, with the special case of one iteration corresponding to the basic version of  singular spectrum analysis (SSA). 
Numerical studies of Cadzow iterations and SSA are provided in \cite{gillard2010cadzow}. Versions and the relative merits of SSA and Cadzow iterations where observations are allocated different weights are described in \cite{zvonarev2017iterative}.

\subsubsection{Convergence of Cadzow iterations and  improvements}
As seen from the form of the Cadzow iterations, the matrices $\mathbf{A}^{(i)}$ get closer to both the set $\mathscr{S}$ of structured matrices and the set $\mathscr{M}_{\le r}$ of the manifold of fixed rank; but they are not, in general guaranteed to converge to a point in the intersection of $\mathscr{S}$ and $\mathscr{M}_{\le r}$.
However, if a point in the intersection is well-behaved (satisfies some regularity conditions, see \cite{lewis2008ap} and \cite{andersson2013alternating}), then the Cadzow iterations converge linearly; as argued in \cite[\S 7]{andersson2013alternating}, for Hankel matrices almost all matrices in the intersection are well-behaved.
Also, for rank-$1$ matrices, it was shown in \cite{knirsch2021optimal} that the Cadzow iterations always converge.
For rank-$r$ matrices, as shown \cite{gillard2013optimization,zvonarev2017iterative} there always exists a convergent subsequence of Cadzow iterations.

Another issue is that the Cadzow iterations do not minimize the cost function \eqref{eq:SLRA_mat_form} of interest (i.e., it is not guaranteed that the point in the intersection is the closest to $\mathcal{S}(\bp_{0})$); a classic example of convergence to a suboptimal solution is given in \cite{de1993structured}.
An easy and cheap scalar correction is described in \cite{gillard2015stochastic} (see also \cite{zvonarev2017iterative}) is as follows. Suppose $\bZ$ is an approximation of $\bA=\mathcal{S}(\bp_{0})$. Then, if we take $c$ as the solution
\begin{equation}\label{eq:localcorr}
c=\argmin_{b \in \mathbb{R}}\|\bA-b\bZ \|_{\bQ,\bR}^{2}
\end{equation}
then this  defines a better approximation in the sense
$
\|c\bZ-\mathcal{S}(\bp_{0}) \|_{\bQ,\bR}^{2} \leq \|\bZ-\mathcal{S}(\bp_{0}) \|_{\bQ,\bR}^{2}.
$
Note that such a correction preserves the rank and the structure of the matrix (e.g., Hankel).
In case of the Frobenius norm ($\bQ,\bR = \bI$), the correction can be simply computed as
\[
c=
\frac{\trace( \bZ^{T}\mathcal{S}(\bp_{0}))}{\trace (\bZ^{T}\bZ)}.    
\]
The correction can be applied at each step or at the end of the algorithm.

A further improvement was proposed in \cite{schost2016quadratic}, where instead of  $\pi_{(\bQ,\bR)}^{(r)}(\bA^{(i-1)})$ on the space $\mathscr{S}$ of structured matrices, a Newton-like step was suggested: the next iteration $\bA^{(i-1)}$ is computed as a projection of $\bA^{(0)}$ on the intersection of $\mathscr{S}$ and  the tangent space to the manifold of $\mathscr{M}_{r}$ of fixed rank matrices at $\pi_{(\bQ,\bR)}^{(r)}(\bA^{(i-1)})$.
It was shown that this algorithm is quadratically convergent, but the cost of each iteration becomes higher.
Another recently proposed modification of Cadzow \cite{condat2015cadzow}  is based on mimicking proximal splitting methods from convex optimization. 
Yet another method which can be seen as an modification of Cadzow, is the alternating direction method of multipliers proposed in \cite{andersson2014admm}; it requires however storing the dual variables, which increases the storage complexity compared to Cadzow iterations.

\subsubsection{Multistart/global optimization}\label{sect:apbr}
With suitable modifications, \eqref{eq:cad} can be used as the basis to construct algorithms known to (at least linearly) converge to the global optima. Possible options include running iterations of \eqref{eq:cad} from several starting points (akin to random multistart), as a multi-stage algorithm and also as an evolutionary method.

The main algorithm studied within the family constructed by the authors in \cite{gillard2015stochastic} is the so-called multistart APBR (alternating projections with backtracking and randomization) where theoretical conditions are also given for its convergence to the global optima. Within this algorithm, backtracking to the initial data is beneficial to ensure that the approximation remains close to it, whilst adding randomisation assists in moving the iterations away from possible local minima.

The multistart APBR is briefly described as follows. Let $\bX$ denote a random Hankel matrix which corresponds to a realization of a white noise Gaussian process ${\boldsymbol \xi}=(\xi_{1},\xi_{2},\ldots, \xi_{N})$, that is, $\bX=\mathcal{S}({\boldsymbol \xi})=\mathcal{H}_{L}({\boldsymbol \xi})$. We assume each $\xi$ are independent Gaussian random variables with zero mean and some variance to be set by the user. 

Set $\bA^{(0)}=\bA=\mathcal{S}(\bp_{0})$,  In multistart APBR, we run $M$ independent trajectories starting at random Hankel matrices
$
\bA^{(0,j)}=(1-s_0)\bA^{(0)}+s_0\bX,
$
with some $0 \leq s_{0} \leq 1$, where for trajectory $j=1,\ldots, M$,
\begin{equation}\label{eq:apbr}
\bA^{(i+1,j)}=c_{i,j}\left[(1-\delta_{i})\pi^{\mathcal{H}}\left(  \pi^{(r)}(\bA^{(i,j)})  \right)+\delta_{i}\bA^{(0)}+\sigma_{i}\bX\right],
\end{equation}
where $\pi^{\mathcal{H}} = \pi^{\mathcal{H}}_{(\bQ,\bR)}$, $\pi^{(r)} = \pi^{(r)}_{(\bQ,\bR)}$
and $c_{i,j}$ is the correction defined in \eqref{eq:localcorr} for the matrix in brackets. Each trajectory is either run until convergence or for a pre-specified number of iterations.  If $s_0=\delta_i=\sigma_i=0$ then the iterations \eqref{eq:apbr} co-incide with those of \eqref{eq:cad} albeit with some local improvement as already described.   If $s_0>0$ then the $j$-th trajectory of multistart APBR starts at a random matrix in the neighbourhood of $\bA^{(0)}$, with the width of this neighbourhood controlled by the parameter $s_0$. If $\sigma_{i}>0$ then there is a `random mutation' at iteration $i$. When $\delta_{i}>0$ the $i$-th iteration `backtracks' towards $\bA^{(0)}$  At some point it is recommended to set $\delta_i=\sigma_i=0$ and directly have iterations of the form \eqref{eq:cad} to accelerate convergence to the set of rank-$r$ Hankel matrices. Recommendations as to the selection of the parameters $\delta_{i}$ and $\sigma_{i}$ are also offered in \cite{gillard2015stochastic}.

\subsection{Local optimization}
There is a significant body of work adressing the  problem \eqref{eq:SLRA} focused on developing computationally efficient methods for computing a locally optimal solution.
Two types of parameterization are typically used: kernel and image representation, where the latter is related to optimizing the manifold of rank-$r$ matrices.
We provide below a summary of these approaches, with a focus on Hankel SLRA
(see also the introduction of \cite{ishteva2014factorization} for more details on kernel vs. image representation).

We note that there is no guarantee that a globally optimal solution is found, and the solution will depend on the initial point selected.
Many existing algorithms are highly sensitive to the choice of the initial point, and due to the existence of several local minima often do not move significantly from it \cite{gillard2015stochastic}.
However, these methods can be improved by choosing different (random) starting points.

\subsubsection{Kernel representation and variable projection}
Consider the Hankel low-rank approximation problem \eqref{eq:SLRA} in the vector formulation with $\mathcal{S} =  \mathcal{H}_{r+1}$ (so that $L = r+1$), which aims at reducing the rank by one\footnote{Note that  by virtue of Corollary~\ref{cor:shape_hankel_fat}, any Hankel SLRA problem with $L > r$ can be reduced to the case $L = r+1$}.
Then thanks to Remark~\ref{rem:lrr_and_rank}, we can replace the rank constraint with a constraint on the kernel of the Hankel matrix, which leads to an equivalent formulation \begin{equation}\label{eq:slra_kernel}
\begin{split}
&\min_{\boldsymbol{\theta} \in \mathbb{R}^{r+1}, \boldsymbol{\theta} \neq 0, \widehat{\bp} \in \mathbb{R}^{r+1}} \|\widehat{\bp}-\bp_{0}\|^2 \\
&{\rm subject\,\, to} \quad  \boldsymbol{\theta}^{\top}\mathcal{H}_{r+1} (\widehat{\bp}) = 0,
\end{split}
\end{equation}
where $\boldsymbol{\theta}^{\top}\mathcal{H}_{r+1} (\bp) = 0$ is a shortcut for the constraint  \eqref{eq:kernel_rank_constraint},
and  the norm is squared to make the objective smooth.
The reformulation \eqref{eq:slra_kernel} is called the \emph{kernel representation}, because the search space is augmented  by an additional vector parameterizing the kernel.
The kernel representation is also closely related to desingularization in  algebraic geometry \cite{khrulkov2018desingularization}.
Also note the kernel representation was described here only for rank reduction of the given matrix by $1$,  rank reduction by more than one is also possible (which leads to minimization of a cost function on a Grassmann manifold, see \cite{slra-software}).

In the kernel representation, we are free to choose a scaling of the  vector $\boldsymbol{\theta}$, as it does not change the cost function.
Earlier approaches  addressed Hankel SLRA from a  structured total least squares \cite{abatzoglou1991tsp} viewpoint (see \cite{markovsky2010bibliography} or \cite{markovsky2006siam} and for a survey),
which corresponds to a particular choice of $\theta_r = 1$ (fixing one of the coordinates to $1$, which parameterizes a part of the search space).
 This was the approach of De Moor \cite{de1993structured}, who developed a non-linear extension of the SVD named the Riemannian SVD, which aims to solve a Lagrange multipliers formulation of \eqref{eq:SLRA}, and also  Rosen et al. \cite{rosen1996total}, who proposed a linear approximation to \eqref{eq:slra_kernel} and then use Newton's method.
A modification of this approach for rank reduction more than by $1$ was  described in Van Huffel \cite{van1996formulation}, but this involves Kronecker products which dramatically inflate the dimension of the involved matrices.
Structured total least squares can be also viewed as the problem of estimating an errors-in-variables regression model which has constraints on the matrix structure.

One of the most successful approaches is the variable projection approach,  initially proposed in  \cite{markovsky2006siam}, which was developed in a series of papers  \cite{slra-software,slra-missdata}). It uses the fact that for a fixed $\boldsymbol{\theta}$, the minimum of \eqref{eq:slra_kernel} (denoted as $f(\boldsymbol{\theta})$) has a closed form solution and thus the SLRA problem is equivalent to the minimization of $f(\boldsymbol{\theta})$ (and thus the variable $\widehat{\bp}$ was eliminated).
The cost function $f(\boldsymbol{\theta})$ can be minimized by general purpose constrained optimization methods, such as Gauss-Newton and Levenberg-Marquardt.
For elementwise weights, the complexity of evaluating $f(\boldsymbol{\theta})$ and its gradient is $O(rN)$ and the complexity of evaluating an approximation of the Hessian is $O(r^2N)$  \cite{usevich2014varpro}.
However, this approach is more difficult to implement efficiently in case of missing data \cite{slra-missdata} (unlike other approaches, presented in the previous and the following subsections).

There are some further improvements of the variable projection approach and other methods that operate in the kernel representation. For example, in \cite{zvonarev2021low,zvonarev2021fast}, a modified Gauss-Newton iteration was proposed, that makes use of the parameterization of the tangent space and a fast projection onto it.
A convex relaxation based on the kernel representation was proposed in \cite{cifuentes2021convex} that globally solves the SLRA if  $\mathcal{S}(\bp_{0})$ is  in a small neighborhood of a low-rank and structured matrix.
Another recent approach  using kernel representation is based on an iterative reweighting scheme \cite{zhang2019optimal}. An augmented Lagrangian approach in kernel representation was proposed in \cite[Ch. 5]{borsdorf2012structured}.

\subsubsection{Image representation and other methods}
The image representation approach (according to the categorization in \cite{ishteva2014factorization,markovsky2019book}), consists of using the fact that any rank-$r$ matrix $\bX \in \mathscr{M}_r$ can be factorized as $\bX = \bP\bL$, where $\bP \in \mathbb{R}^{L\times r}$ and $\bL \in \mathbb{R}^{r\times L}$ and hence the SLRA problem can be posed as optimization over $\bP$ and $\bL$.
The difficulty in using the image representation for the SLRA problem is that an additional constraint $\bX \in \mathscr{S}$ needs to be handled.
In \cite{ishteva2014factorization}, it was proposed to minimize the following objective:
\[
\min_{\bP,\bL} \|\bP\bL - \bX_0\|_{\mathbf{W}_{m}}^2 + \lambda \| \bP\bL - \pi^{\mathcal{H}} (\bP\bL)\|_F^2,  
\]
where the second term penalizes the distance to structured matrices and $\lambda$ is increased adaptively, from small to large.
For fixed $\lambda$, the objective function is minimized using the alternating least squares strategy (for fixed $\bP$, the problem becomes a least squares problem in $\bL$, and vice versa).
An extension of this approach based on Lagrangian formulation was proposed in \cite{hage2015robust}.
Another extension based on proximal methods was proposed in \cite{liu2020simax} in order to handle multiple rank constraints.
In general, other methods for optimization on the set  (manifold) of rank-$r$ matrices can be applied to such a reformulation of SLRA.

Recently, in a series of papers \cite{guglielmi2017ode,fazzi2018agcd,fazzi2021ode} it was shown that a solution to SLRA can be found as a solution of a gradient system for eigenvalues or singular values. Therefore, it was suggested to solve SLRA by integrating the corresponding ordinary differential equation.
While the approach works very well for some matrix structures (e.g., Sylvester \cite{guglielmi2017ode,fazzi2018agcd}), its behaviour in the Hankel case is not completely clear and ad-hoc adjustments to the algorithm are needed \cite{fazzi2021ode}.

A class of convex relaxation methods have been proposed, whereby a common strategy is to replace the rank constraint with the nuclear norm \cite{fazel2002matrix,fazel2013hankel,blomberg2015regpaths} (these approaches will be mentioned in the next section in the context of time series forecasting).
Tighter relaxations were also proposed based on convex envelopes \cite{caprani2019icassp,andersson2017convex,alshabili2019icassp},\cite[\S VI]{grussler2018convex} which take into account the particular cost function under consideration.
Finally, several iteratively reweighted least squares methods based on surrogates for rank were proposed, which possess local convergence and theoretical guarantees \cite[Ch.3]{kummerle2019irls}, \cite{kummerle2019phd}, and allow for a fast implementation \cite{ongie2017irls}.

\section{Convex relaxations for forecasting}\label{sec:conv_rel}

\subsection{Forecasting as matrix completion}
For forecasting, we introduce some additional notation for simplicity of the presentation.
We denote $N = n+m$, where $n$ is be the length of the given time series that we wish to forecast and  $m$ is the number of observations to be forecast. We use the notation $\bp_{(1:n)}=(p_{1},p_{2},\ldots,p_{n})$ for the first `known' $n$ elements of $\bp$. Two obvious important cases are when $m=0$ (where there are no observations to forecast) and $m>1$. The Hankel matrix with $\bp$ embedded takes the form
\begin{equation}\label{eq:struct_hankel}
\mathcal{S}(\bp)  =
\left(
\begin{array}{cccccc}
\cellcolor{lightgray}p_1  & \cellcolor{lightgray}p_2&  \cellcolor{lightgray}\cdots & \cellcolor{lightgray} \cdots  & \cellcolor{lightgray} \cdots  & \cellcolor{lightgray}p_K \\
\cellcolor{lightgray}p_2 & \cellcolor{lightgray}p_3& \cellcolor{lightgray} \cdots  &  \cellcolor{lightgray}\iddots & \cellcolor{lightgray} \cdots & \cellcolor{lightgray} \vdots \\
\cellcolor{lightgray}\vdots &\cellcolor{lightgray}  \cdots &\cellcolor{lightgray}\iddots &  \cellcolor{lightgray}\iddots & \cellcolor{lightgray} \cdots &  \cellcolor{lightgray} p_{n} \\
\cellcolor{lightgray}\vdots & \cellcolor{lightgray} \iddots &\cellcolor{lightgray}\iddots &  \cellcolor{lightgray}\iddots & \cellcolor{lightgray} \iddots & p_{n+1} \\
\cellcolor{lightgray}\vdots &   \cellcolor{lightgray} \iddots &   \cellcolor{lightgray}\iddots&   \cellcolor{lightgray}\iddots&  \iddots & \vdots \\
\cellcolor{lightgray}p_L &   \cellcolor{lightgray}\cdots &  \cellcolor{lightgray}p_{n}  & p_{n+1}&   \cdots & p_{n+m}
\end{array}
\right).
\end{equation}
In \eqref{eq:struct_hankel}, the grey-shaded values are ``known'' and others are ``missing''.

\subsection{Nuclear norm relaxation}
Often the nuclear norm relaxation of the rank constraint is used, since the nuclear norm is a convex envelope of the rank function. Formally, for a matrix $\mathbf{X} \in \mathbb{C}^{L\times K}$ its nuclear norm is defined as
\[
\|\bX\|_* = \sum_{k=1}^{\min(L,K)} |\sigma_k (\bX)|,
\]
where $\sigma_k(\mathbf{X})$ are the singular values of $\bX$. This idea has a rich recent history and has been discussed by several eminent authors, see \cite{Candes2009,candes2010matrix,recht2010guaranteed} for example.

There is an analogy between the use of the nuclear norm in matrix completion problems and the use of the $\ell_1$ norm for sparse approximation (see \cite{rish2014sparse} for an overview of this topic) and its use has underpinned several successful works describing how to impute the missing values of a matrix with astonishing accuracy. More generally nuclear norm relaxation has proved a useful tool in:
spectral estimation \cite{akccay2014subspace}, recommender systems \cite{kang2016top}, 
system identification \cite{liu2010interior} and several other application areas.

In what follows we will describe the low-rank matrix completion problem  and the remainder of this review describes work which followed the initial proposal of Butcher and Gillard \cite{butcher2016simple} who proposed that a time series be embedded into a Hankel matrix, with the data to be forecasted being missing elements stored in the bottom right-hand corner of this matrix.

As well as low-rank approximation offering the potential for many new algorithms for forecasting, the main theoretical question is when the nuclear norm relaxation solves the original low-rank matrix completion problem \cite{markovsky2012sysid}. Existing work mainly assumes that the position of the missing elements is random \cite{chen2014robust}, and often that the known entries are non-random. In general, unstructured matrices are studied, yielding many of the theoretical results in famous works such as \cite{candes2010matrix} redundant for our problem of Hankel structure with missing values in a particular section of the matrix.

\subsection{Exact low-rank matrix completion}
Let $\bp_{0}=(p_{0,1},p_{0,2},\ldots,p_{0,n})$ be a given vector of observations (a time series).
For a given matrix structure \eqref{eq:linstr}, the exact Structured Low-Rank Matrix Completion (SLRMC) problem is posed as
\begin{equation}\label{eq:slrmc}
\tilde{\bp}=\argmin_{\bp \in \mathbb{R}^{(n+m)}} \rank  \,\mathcal{S}(\bp)\,\,\,\,  {\rm subject\,\, to}\,\,\,\bp_{(1:n)}=\bp_{0}. \,
\end{equation}
The implicit low-rank assumption of the Hankel matrix corresponds to the class of time series of so-called finite rank, which are described in \cite{usevich2010signal}. Here the notation $\bp_{(1:n)}$ represents the first $n$ (known) elements of $\bp$ and so the constraint $\bp_{(1:n)}=\bp_{0}$ ensures that the known elements are fixed in the approximation. The optimization problem \eqref{eq:slrmc} thus acts on $m$ variables, corresponding to the $m$ missing values.

A convex relaxation of \eqref{eq:slrmc} is obtained by replacing the rank with the nuclear norm:
\begin{equation}\label{eq:nmorig}
\hat{\bp}=\argmin_{\bp \in \mathbb{R}^{(n+m)}} \| \mathcal{S}(\bp) \|_{\ast}\,\,\,\,  {\rm subject\,\, to}\,\,\,\bp_{(1:n)}=\bp_{0} \, .
\end{equation}
The intuition behind this relaxation is the same as for using the $\ell_1$ norm in compressed sensing: the nuclear norm is expected to force all but a few singular values to be zero hopefully rendering solutions of \eqref{eq:nmorig} at least close to \eqref{eq:slrmc}.


\subsubsection{Solution of the exact matrix completion}\label{eq:slrmc_solution}
For the Hankel matrix case, the solution of \eqref{eq:slrmc} is known. Before describing the approximate case we recall the class of so-called time series of finite rank and hence give the solution of \eqref{eq:slrmc}. The performance of the nuclear norm relaxation is important in this exact case and is crucial to understand the behaviour of forecasting in the approximate case to be introduced later. 

For time series of finite rank, the solution to \eqref{eq:slrmc} with Hankel structure \eqref{eq:struct_hankel} can be found by the so-called minimal rank extension of Hankel matrices, solved by \cite{iohvidov1982toeplitz} for square matrices and \cite{heinig1984algebraic} for rectangular matrices. We summarise their results in the following theorem.

\begin{theorem}\label{thm:exact_completion}
Consider the Hankel matrix structure \eqref{eq:struct_hankel} and let $r \leq \min(L-1,K-1,\frac{n}{2})$.
For any complex-valued time series $(p_{0,1},p_{0,2},\ldots,p_{0,n},\ldots)$ of finite rank $r$ and vector $\bp_{0}=(p_{0,1},p_{0,2},\ldots,p_{0,n})$ then the unique solution of \eqref{eq:slrmc} is given by
$$
\tilde{\bp}=(p_{0,1},p_{0,2},\ldots,p_{0,n},p_{0,n+1},\ldots,p_{0,n+m})
$$
where $p_{0,n+1},\ldots,p_{0,n+m}$ are found from the linear recurrence formula \eqref{eq:LRR}.
\end{theorem}

A proof is available in, for example, \cite{usevich2016hankel}.

\subsubsection{Known results on the performance of the nuclear norm}
Recent results are described in \cite{usevich2016hankel,gillard2018structured} with the first described by \cite{dai2015nuclear}. We summarise some of the known results here.

\begin{theorem}
Let $\bp_{0}=(p_{0,1},p_{0,2}, \ldots, p_{0,n})$  be a complex-valued vector given as
\[
p_{0,k} = c \lambda^{k},  \quad  k=1,\ldots,n,
\]
where $\lambda \in \mathbb{C}$.
\begin{itemize}
\item If $|\lambda| \le 1$, the solution of \eqref{eq:slrmc}, i.e.
\[
p_{0,k} = c \lambda^{k},  \quad  k=n+1,\ldots,n+m
\]
is also a solution of \eqref{eq:nmorig}, in particular, if $|\lambda| < 1$, then the solution of \eqref{eq:nmorig} is unique.

\item If $|\lambda| > 1$ (and $L=K=n=m+1$), then the unique solution of \eqref{eq:nmorig} is given by
\[
p_{0,n+k} = c \frac{\lambda^{n}}{\overline{\lambda}^{k}},  \quad  k=1,\ldots,m.
\]
\end{itemize}
\end{theorem}

A version of Theorem 3 was originally given in \cite{dai2015nuclear}, and has been later refined in \cite[Thm. 6]{usevich2016hankel} and \cite[Prop. 4.2 and 4.5]{gillard2018structured}. In this case the exponential needs to be damped for the solution of \eqref{eq:slrmc} to equal that of \eqref{eq:nmorig}. Despite concerning what might be a simple case of a ``rank 1" Hankel matrix it yields many clues as to when the nuclear norm relaxation will coincide with the original problem when there is no convex relaxation, and this is described in the next theorem, originally stated in \cite{usevich2016hankel} and extended to rectangular matrices in  \cite{gillard2018structured}.

\begin{theorem}\label{thm:nn_rankr}
Fix $L,K,n,m$ in \eqref{eq:struct_hankel} as in Theorem~\ref{thm:exact_completion}, such that $r \le \frac{n}{2}$. Then there exists a number $0 < \rho_{max,r,m} < 1$ such that
for any complex-valued time series of finite rank $r$  with
\[
|\lambda_j| < \rho_{max,r,m},
\]
the solution of \eqref{eq:nmorig} is unique and coincides with  the solution as described in Theorem~\ref{thm:exact_completion}.
\end{theorem}

\subsection{Approximate matrix completion}\label{sec:approx_mc}
The approximate rank minimization can be posed as follows. Given $\bp_{0}$,  vector of weights $\bw \in \mathbb{R}^{n}$, $m\geq 0$ and $\tau \geq 0$ find
\begin{equation}\label{eq:slrmc_approx}
\tilde{\bp}=\argmin_{\bp \in \mathbb{R}^{(n+m)}} \rank  \,\mathcal{S}(\bp)\,\,  {\rm subject\,\, to}\,\|\bp_{(1:n)}-\bp_{0}\|_{W} \leq \tau \,,
\end{equation}
which is a special case of the rank minimization problem \eqref{eq:rankmin}.
The parameter $\tau$ controls the precision of approximation and two extreme cases can be distinguished:
\begin{itemize}
\item if $\tau = 0$, then  \eqref{eq:slrmc_approx} is equivalent to the exact matrix completion problem  \eqref{eq:slrmc}, i.e. there is no approximation;
\item if $m=0$ in \eqref{eq:slrmc_approx}, then $\tilde{\bp}$ is an approximation to the given vector $\bp_{0}$ with no forecast.
\end{itemize}
Unlike  \eqref{eq:slrmc}, the problem \eqref{eq:slrmc_approx} does not have a known solution.

In \cite{gillard2018structured}, the following relaxation of \eqref{eq:slrmc_approx} using the nuclear norm was introduced and studied:
\begin{equation}\label{eq:nm}
{\bp}_{\ast}=\argmin_{\bp \in \mathbb{R}^{(n+m)}}  \| \mathcal{S}(\bp) \|_{\ast}  \,\,\, {\rm subject\,\, to}\,\,\,\|\bp_{(1:n)}-\bp_{0}\|_{W} \leq \tau \, .
\end{equation}

There are alternative ways to extend the problem \eqref{eq:nmorig}  to approximate versions. For example, consider the following equivalent formulations:
\begin{align}
& \min_{\bp \in \mathbb{R}^{(n+m)}} \|\bp_{(1:n)}-\bp_{0}\|_{W}   \,\,\, {\rm subject\,\, to}\,\,\,\| \mathcal{S}(\bp) \|_{\ast}  \leq \delta, \label{eq:regularisation2}\\
& \min_{\bp \in \mathbb{R}^{(n+m)}} \|\bp_{(1:n)}-\bp_{0}\|_{W}   + \gamma \| \mathcal{S}(\bp) \|_{\ast}, \label{eq:regularisation3}
\end{align}
where $\delta$ and $\gamma$ are regularisation parameters for each of the formulations.

Problems \eqref{eq:nm}, \eqref{eq:regularisation2} and \eqref{eq:regularisation3} are equivalent in the following sense: for any value of $\tau$, there exist $\delta$ and $\gamma$ such that the solutions to \eqref{eq:nm}, \eqref{eq:regularisation2} and \eqref{eq:regularisation3} coincide.
However, the relation between ``equivalent'' $\tau$, $\delta$ and $\gamma$ is not known (a priori).

\section{Numerical examples}\label{sec:ex}

\subsection{Example 1: Structured low-rank approximations}
Here we provide some commentary on the methods described earlier, via a simple example,  to find a structured low-rank approximation of a Hankel matrix in which a time series has been embedded.

Assume we have observed a time-series ($n=20$, $m=0$) of the form $\bp=(p_{1},p_{2},\ldots,p_{20})$ where $p_{j}=s_{j}+\varepsilon_{j}$. Here $s_{j}=\exp(\lambda j)\sin(2\pi\omega j)$ and $\varepsilon_{j}$ is a noise term. We take $L~=~10$ and compare $r=3$ structured low-rank approximations to $\mathcal{S}(\bp)$ using Cadzow iterations \eqref{eq:cad}, multistart APBR (as described in \cite{gillard2015stochastic} and in Section \ref{sect:apbr}), and the variable projections algorithm of Markovsky using the software provided in \url{https://slra.github.io/software-slra.html} and described in \cite{slra-software}, which we will abbreviate as VP. We take $\lambda=0.05$ and $\omega=0.2$. We consider different forms of the noise term $\varepsilon_{j}$.

\begin{figure}[ht!]
\centering
    \subfigure[$\sigma=0.3$ ]{\includegraphics[width=0.16\textwidth]{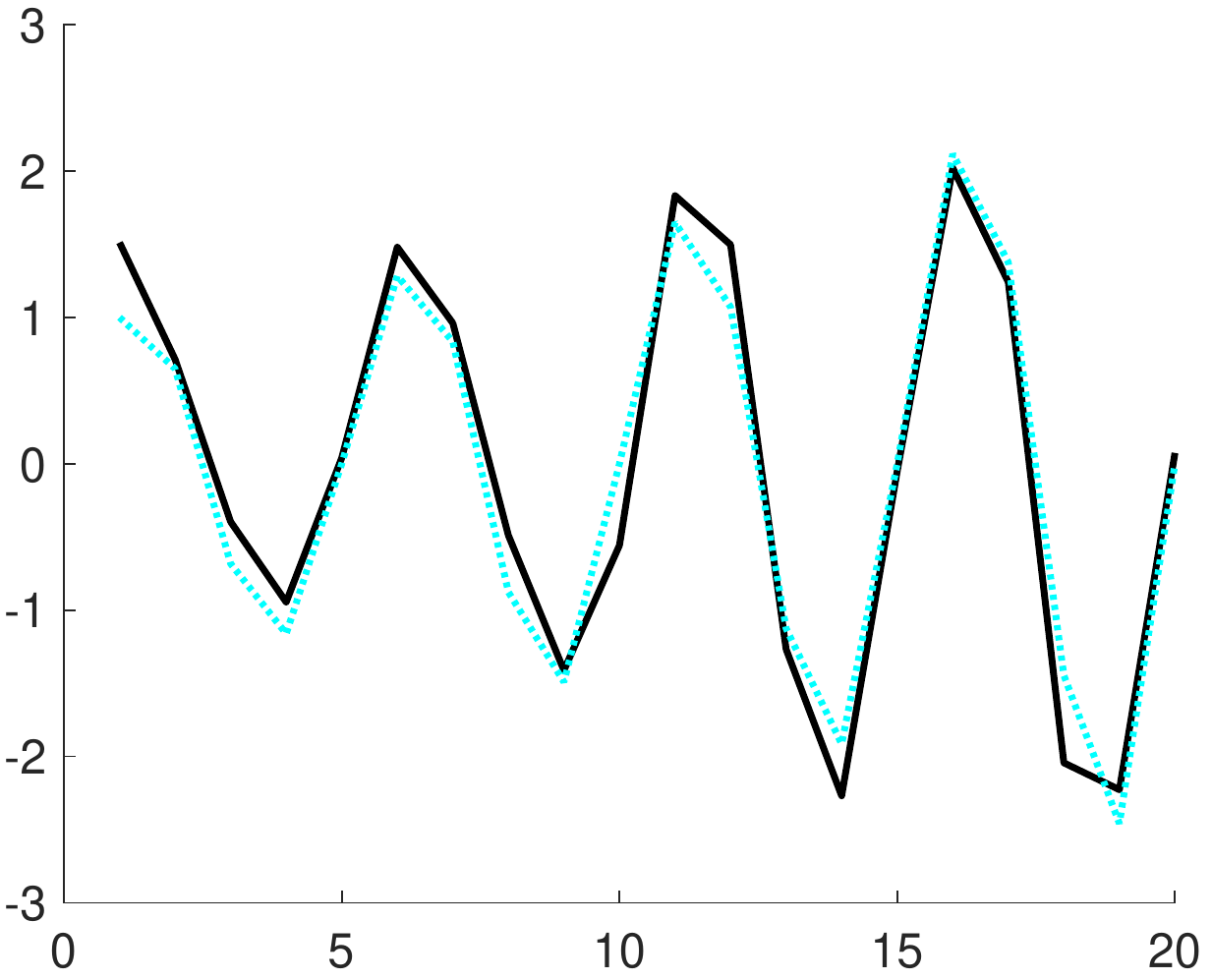}}
    \subfigure[$\sigma=0.6$ ]{\includegraphics[width=0.16\textwidth]{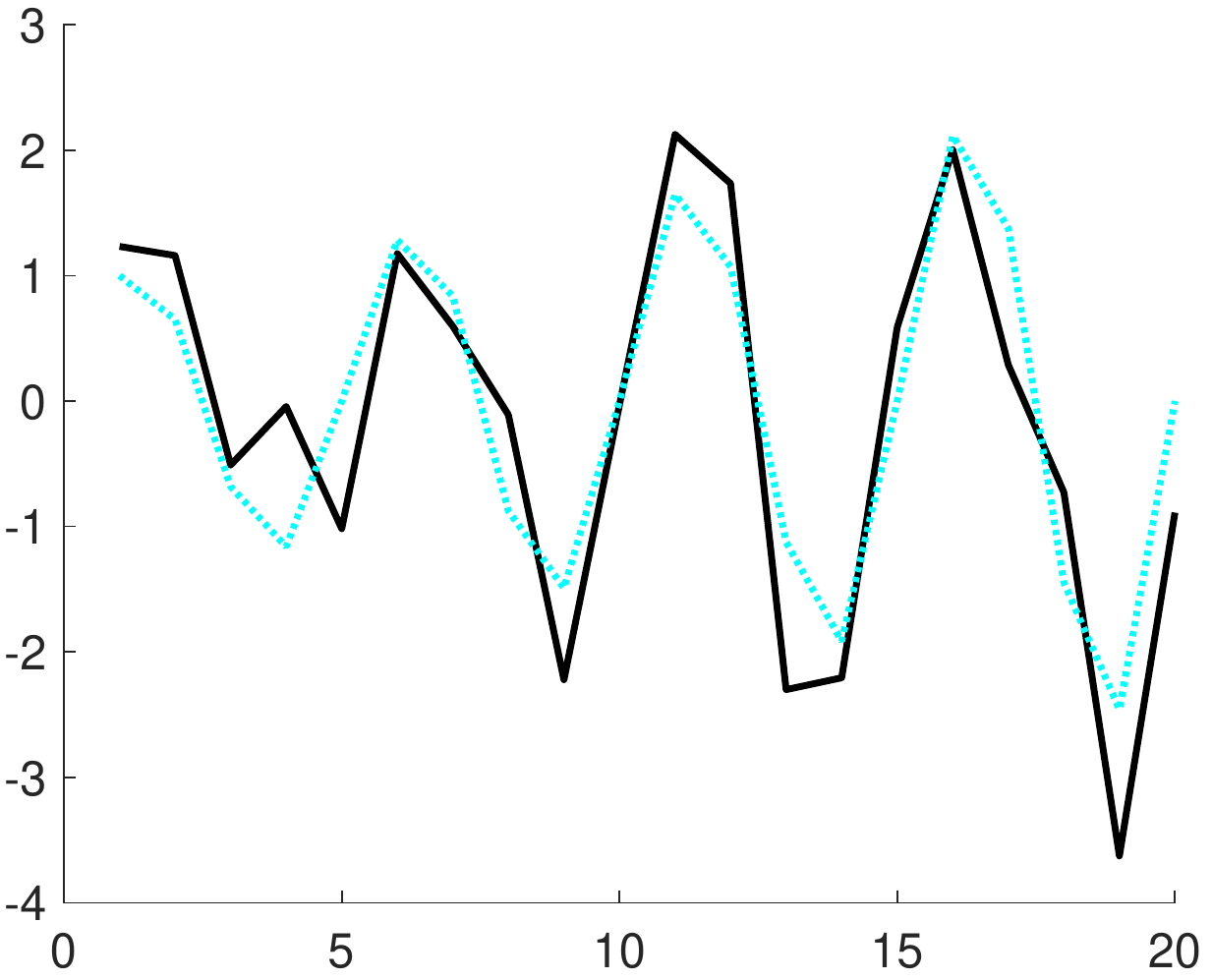}}
    \subfigure[$\sigma=0.9$ ]{\includegraphics[width=0.16\textwidth]{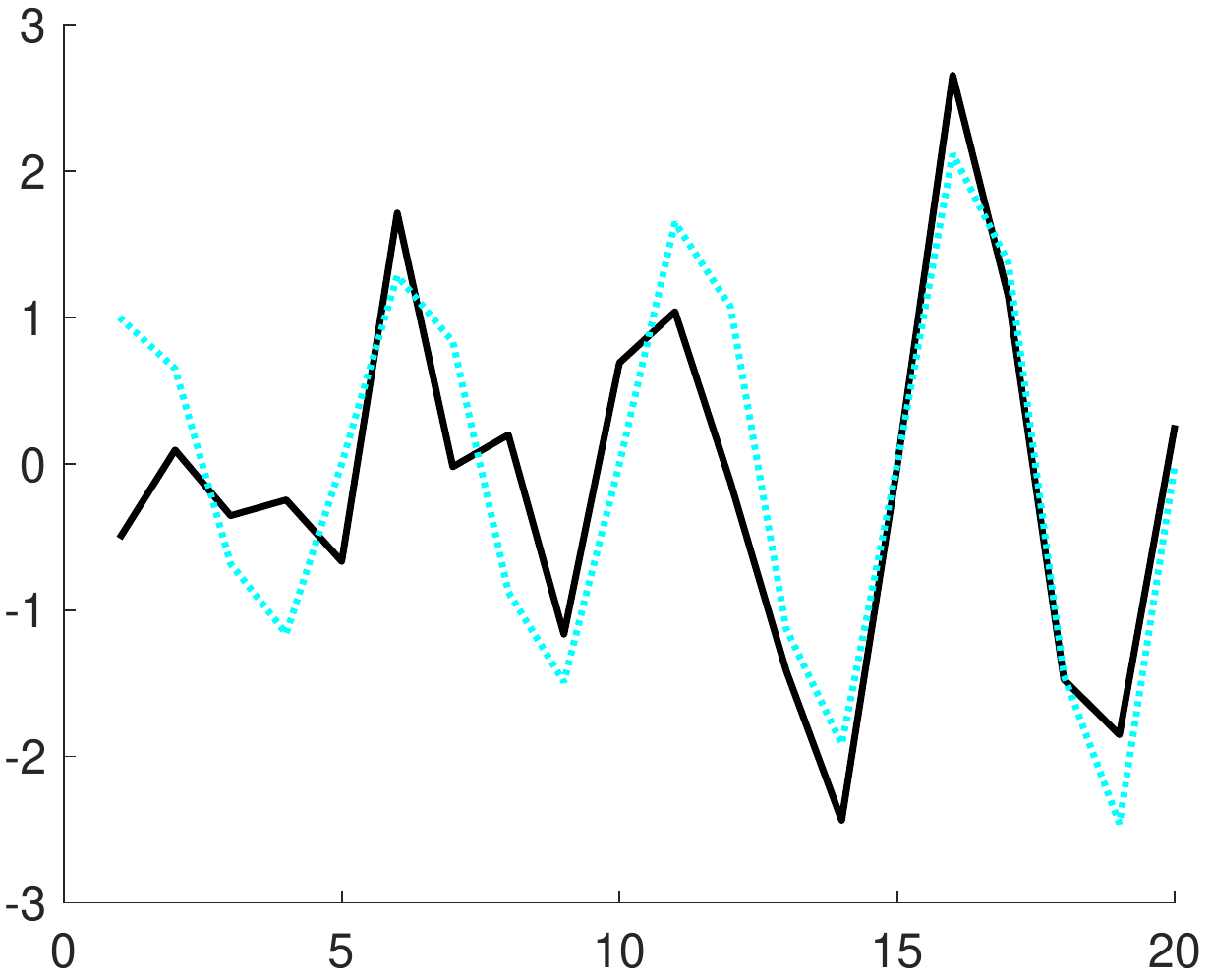}}
  \caption{Plots of $\bp$ for different $\sigma$, plot of time-series without noise in light colour, plot of time-series with noise in black.}
\label{fig:typ}
\end{figure}

\paragraph*{Case 1: White noise, $\varepsilon_{j} \sim N[0,\sigma^2]$}
We first consider the case of so-called white noise, that is  $\varepsilon_{j} \sim N[0,\sigma^2]$ for $j=1,2,\ldots,20$, with some $\sigma^2$, and $cov[\varepsilon_{i},\varepsilon_{j}]=0$ for $i \neq j$.  Plots of typical time-series with different values of $\sigma^2$ are given in Figure \ref{fig:typ}.

\begin{figure}[ht!]
\centering
    \subfigure[$\sigma=0.3$ ]{\includegraphics[width=0.16\textwidth]{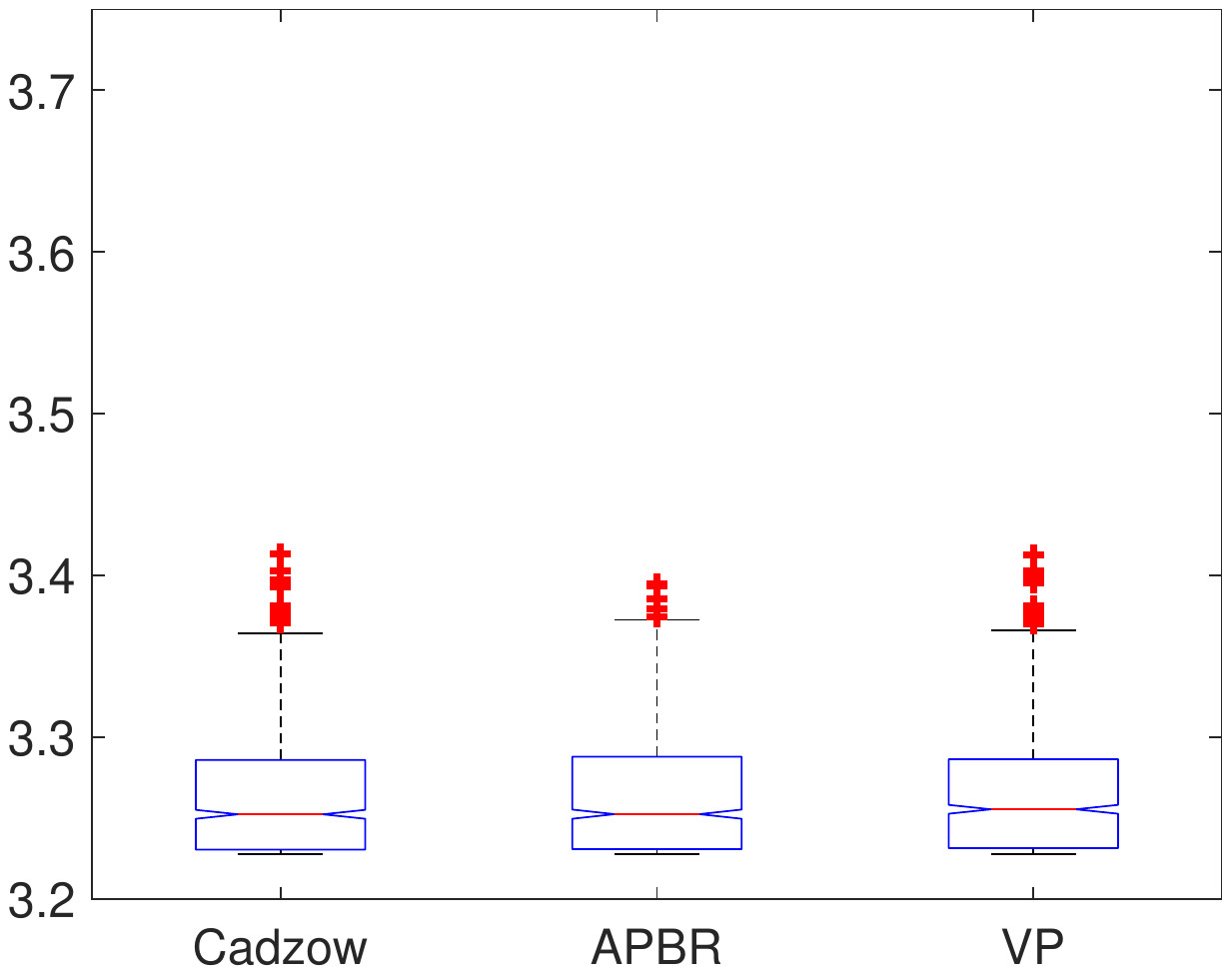}}
    \subfigure[$\sigma=0.6$ ]{\includegraphics[width=0.16\textwidth]{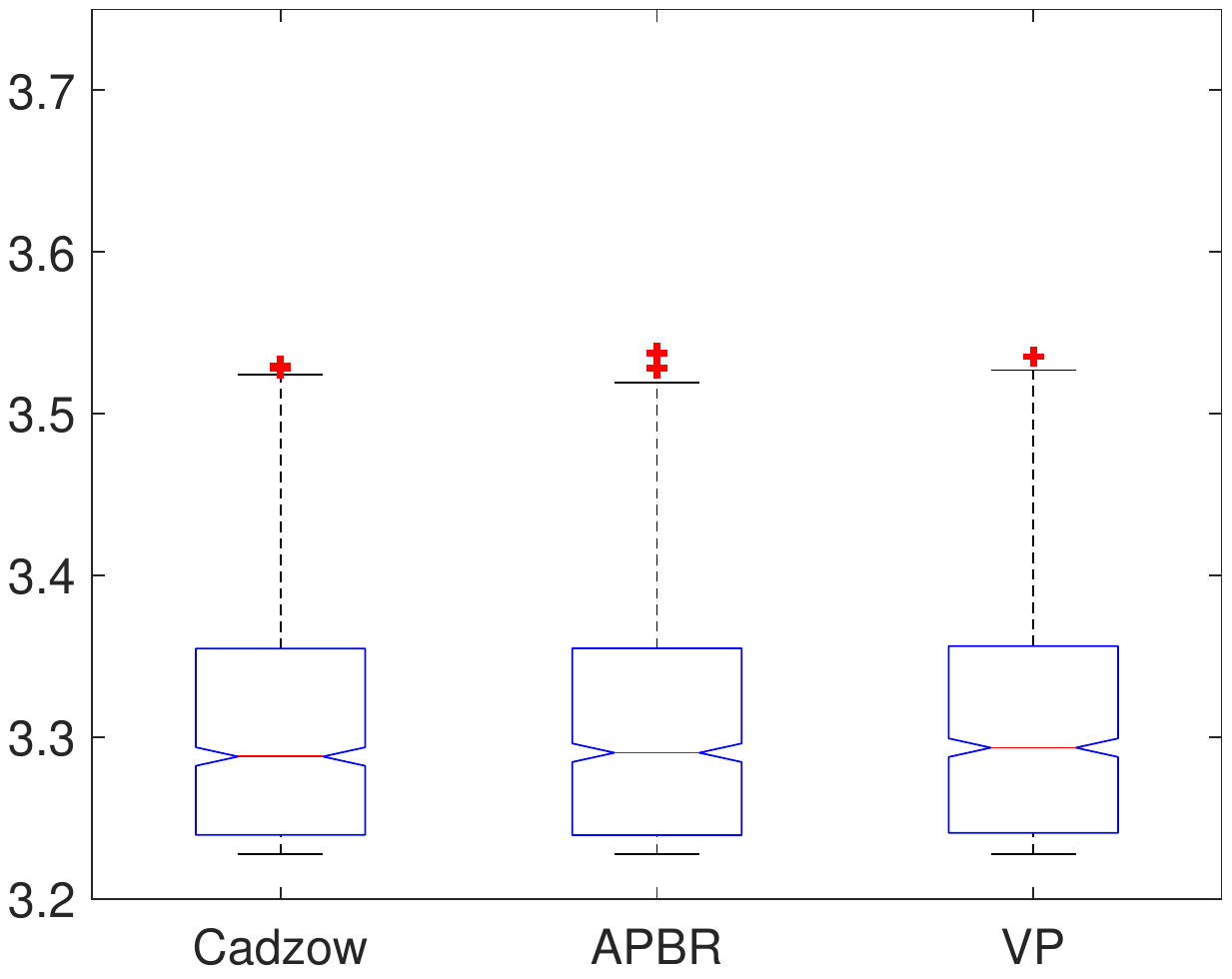}}
    \subfigure[$\sigma=0.9$ ]{\includegraphics[width=0.16\textwidth]{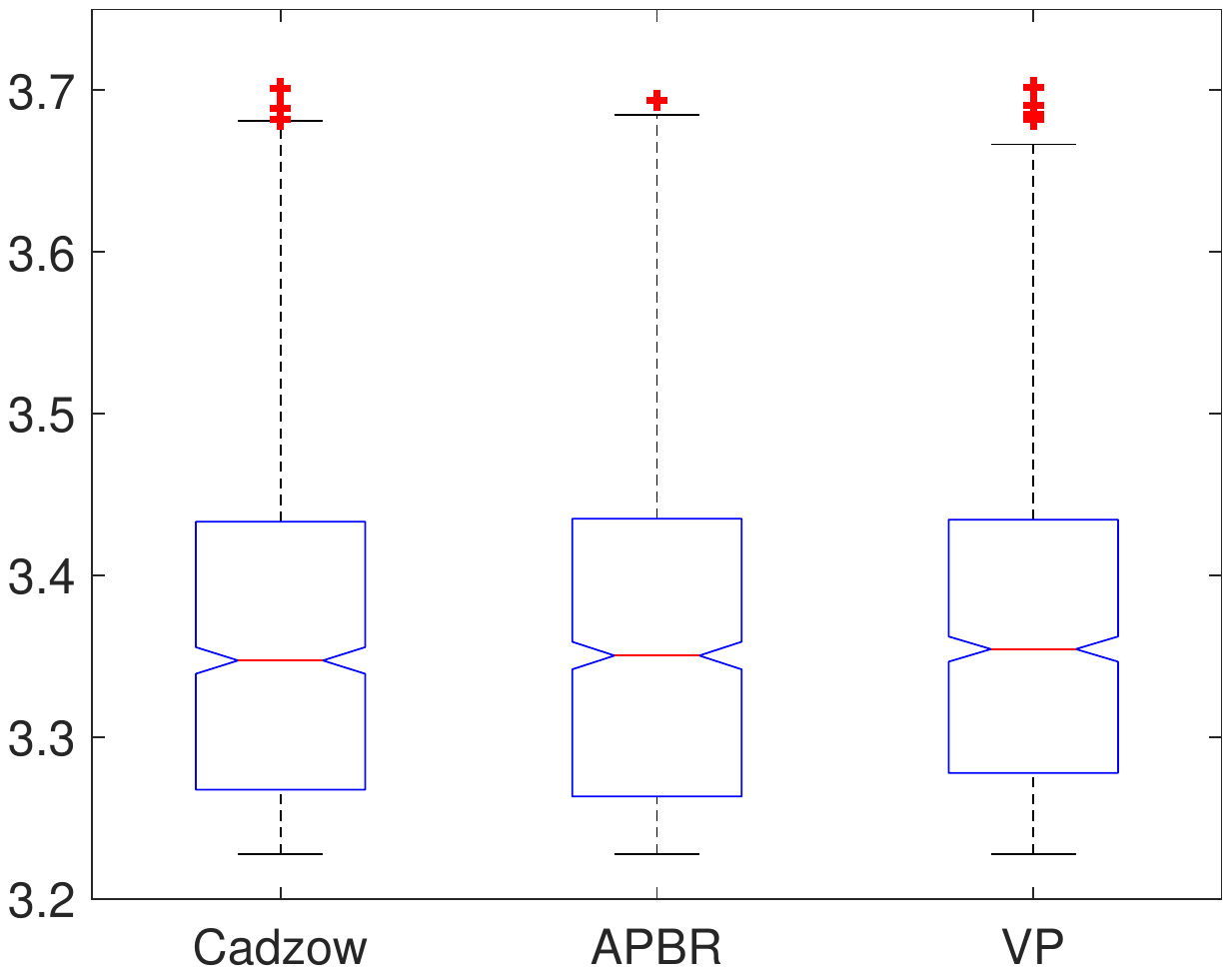}}
  \caption{Boxplots of $\log \|\bp-\bp_{approx} \|_{2}$ taken over 1000 simulations for the case of white noise.}
\label{fig:norm}
\end{figure}

We vary $\sigma$ and below report the norm $\|\bp-\bp_{approx} \|_{2}$ taken over 1000 simulations, where $\bp_{approx}$ is the solution to \eqref{eq:SLRA} as found by the three methods described above. Figure \ref{fig:norm} contains boxplots for each of these methods for different $\sigma$.

\paragraph*{Case 2: Deterministic `noise'}
We now consider the case $\varepsilon_{j}=c (-1)^j$ for some $c>0$. Table \ref{tab:c} reports the norm $\|\bp-\bp_{approx} \|_{2}$ obtained from the three methods considered.
\begin{table}[h!]
\caption{The norm $\|\bp-\bp_{approx} \|_{2}$ obtained from three methods, for different $c$}
\begin{tabular}{|l|c|c|c|c|}
  \hline
            & 0.2        & 0.4        & 0.6        & 0.8        \\ \hline
  Cadzow    & 25.6373    & 26.6507    & 28.2001    & 30.2832    \\
  APBR      & 25.4312    & 26.4329    & 28.1950    & 30.1874    \\
  VP        & 25.6373    & 26.6507    & 28.2001    & 30.2832    \\
  \hline
\end{tabular}
\label{tab:c}
\end{table}

\paragraph*{Case 3: Red noise, $\varepsilon_{j} = \alpha \varepsilon_{j-1}+\eta_{j}$ }
We now consider the case $\varepsilon_{j} = \sigma\left(\alpha \varepsilon_{j-1}+\eta_{j}\right)$, where $\eta_{j} \sim N[0,1-\alpha^2]$ for $j=1,2,\ldots,20$, $\alpha=0.5$,  $cov[\eta_{i},\eta_{j}]=0$ for $i \neq j$, and $\sigma > 0$. This is red noise, more precisely, an autoregressive process of order 1. We vary $\sigma$ and below report the norm $\|\bp-\bp_{approx} \|_{2}$ taken over 1000 simulations, where $\bp_{approx}$ is the solution to \eqref{eq:SLRA} as found by the three methods described above. Figure \ref{fig:norm} contains boxplots for each of these methods for different~$\sigma$.

\begin{figure}[ht!]
\centering
    \subfigure[$\sigma=0.3$ ]{\includegraphics[width=0.16\textwidth]{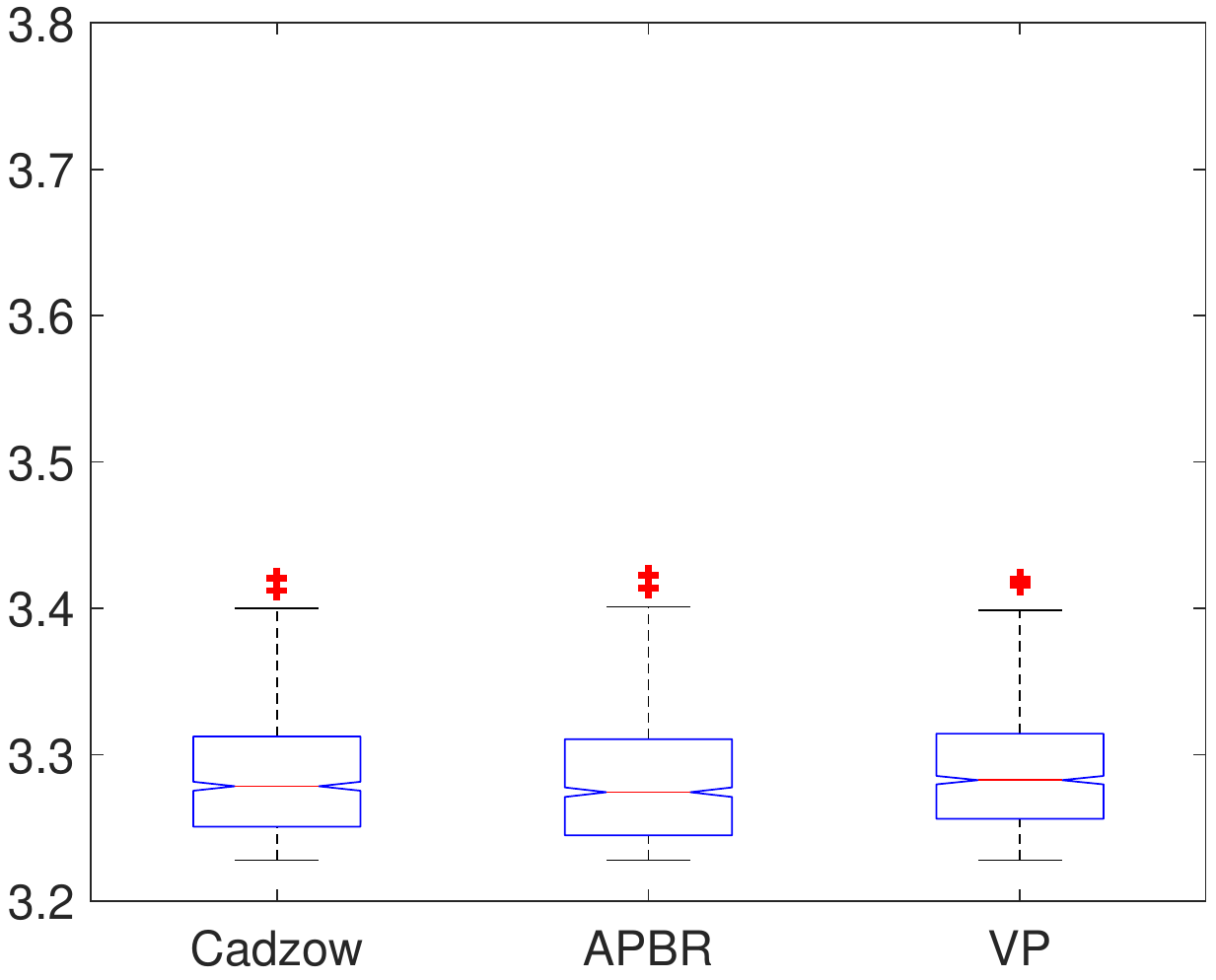}}
    \subfigure[$\sigma=0.6$ ]{\includegraphics[width=0.16\textwidth]{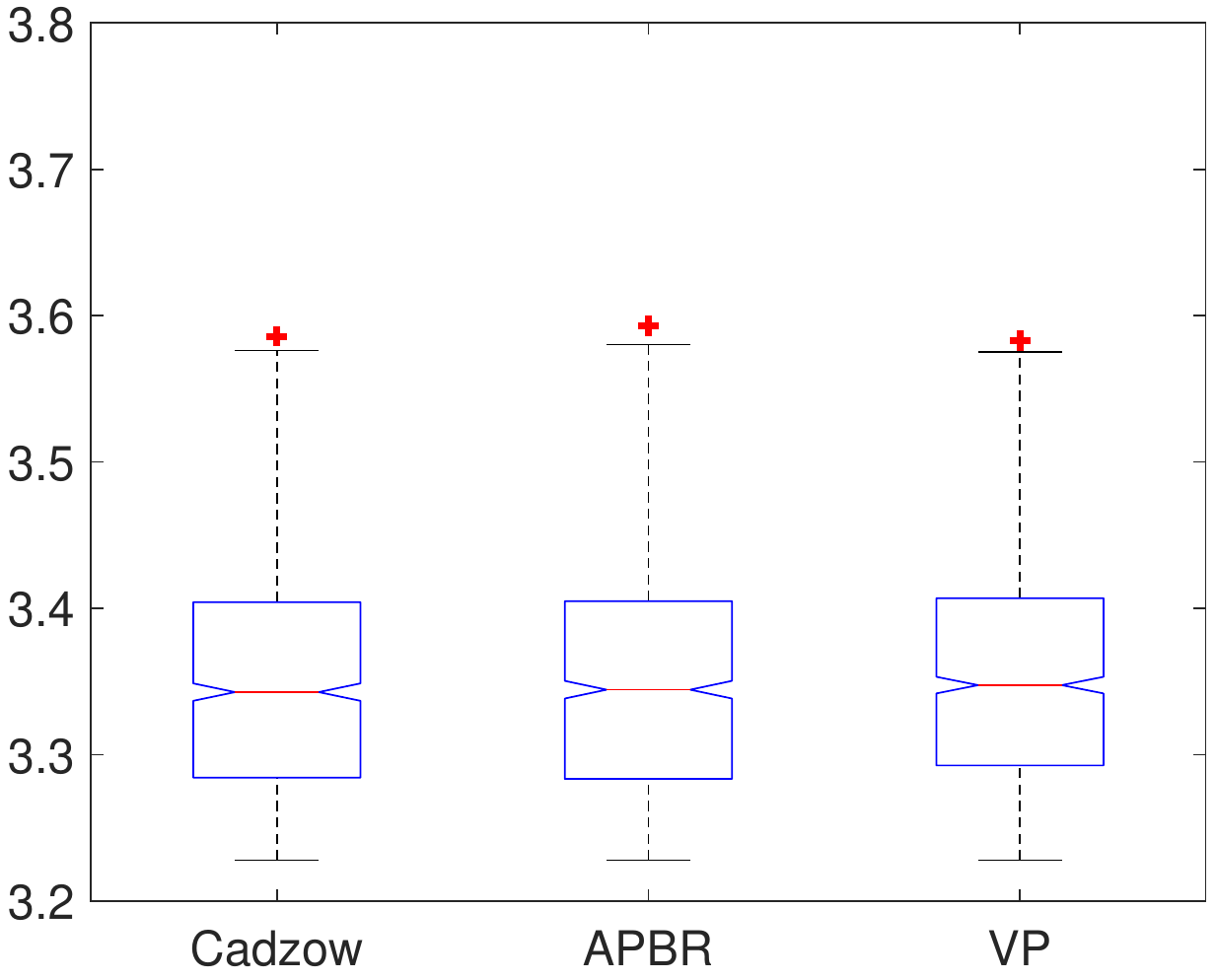}}
    \subfigure[$\sigma=0.9$ ]{\includegraphics[width=0.16\textwidth]{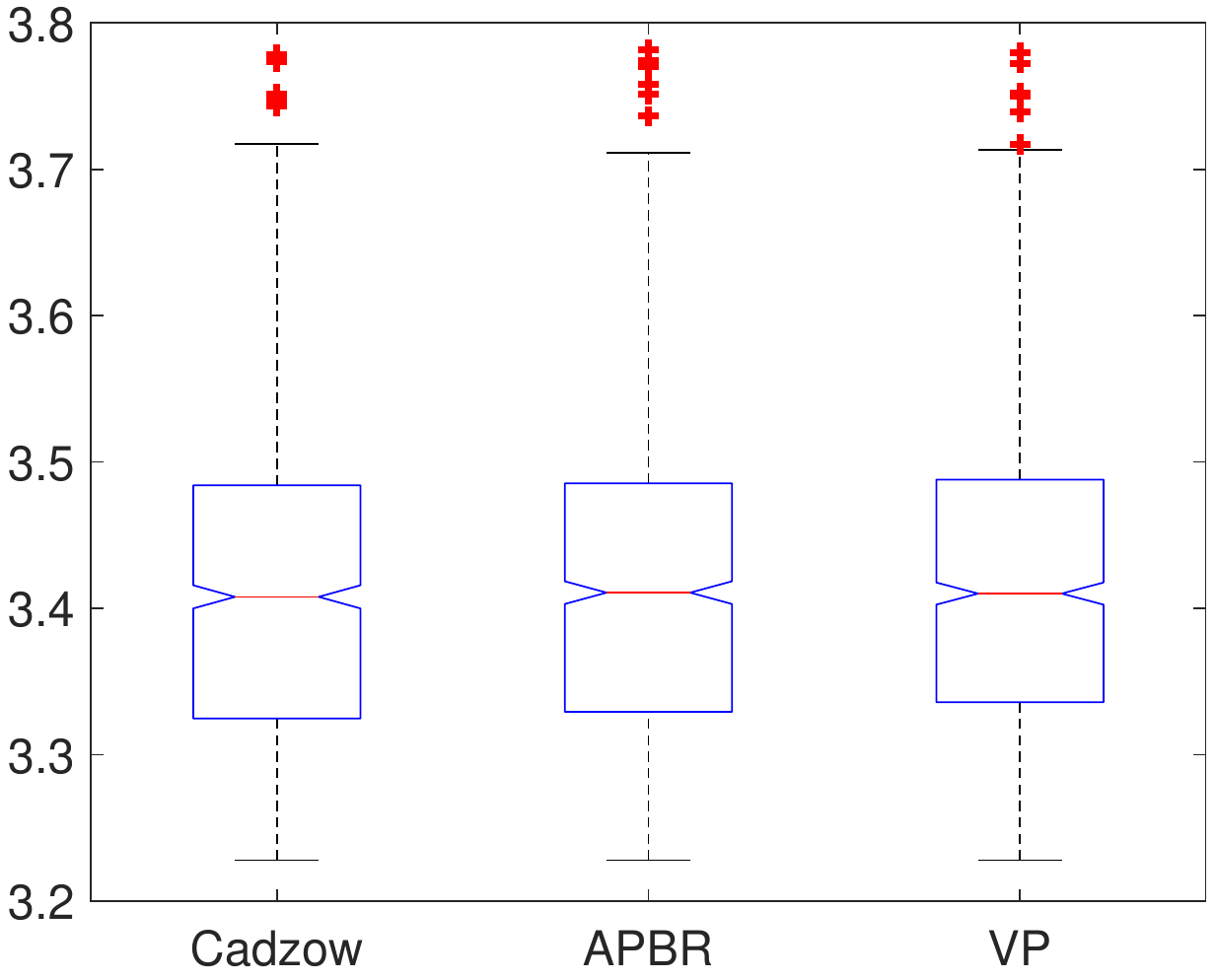}}
  \caption{Boxplots of $\log \|\bp-\bp_{approx} \|_{2}$ taken over 1000 simulations for the case of red noise.}
\label{fig:rednorm}
\end{figure}

We make the following remarks. We note that for the examples considered, there are only marginal differences between each of the methods described. Improvements are found when using the multistart APBR method, which as noted, aims to approximate the global minimum of the problem. One may postulate that in bigger, more complex examples, perhaps this difference would be more pronounced. However it is a pleasing validation of each method that in this example each of the approximations they yield are close. More generally, each method seems relatively robust to the type and amount of additive noise it is exposed to.
Finally, the performance of the VP approach can be improved by choosing different (random) starting points.

\subsection{Example 2: Structured low-rank matrix completion }
To illustrate structured low-rank matrix completion and its application to forecasting we consider a small classical example ``Cowtemp" which records the daily morning temperature of a cow. This time-series is available in several on-line repositories. The time-series consists of 75 observations, and the exercise is to forecast the last $m=14$ observations, based on the first $n=61$ observations. This exercise was also performed in \cite{golyandina2018singular}, comparing forecasts of SSA with various parameter choices. As there, we select $L=28$.

\begin{figure}[ht!]
\centering
    \subfigure[$W_{1}$]{\includegraphics[width=0.16\textwidth]{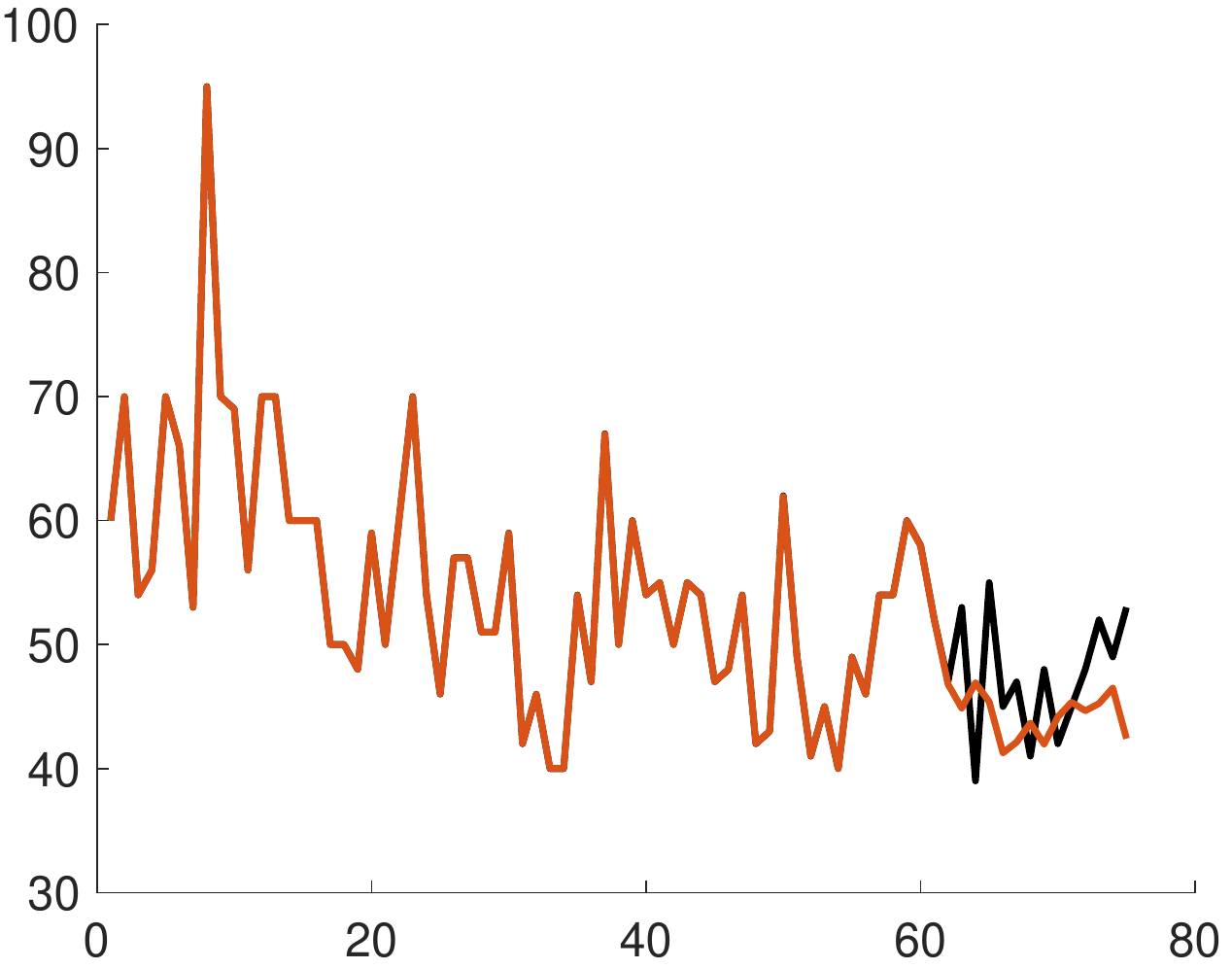}}
    \subfigure[$W_{2}$]{\includegraphics[width=0.16\textwidth]{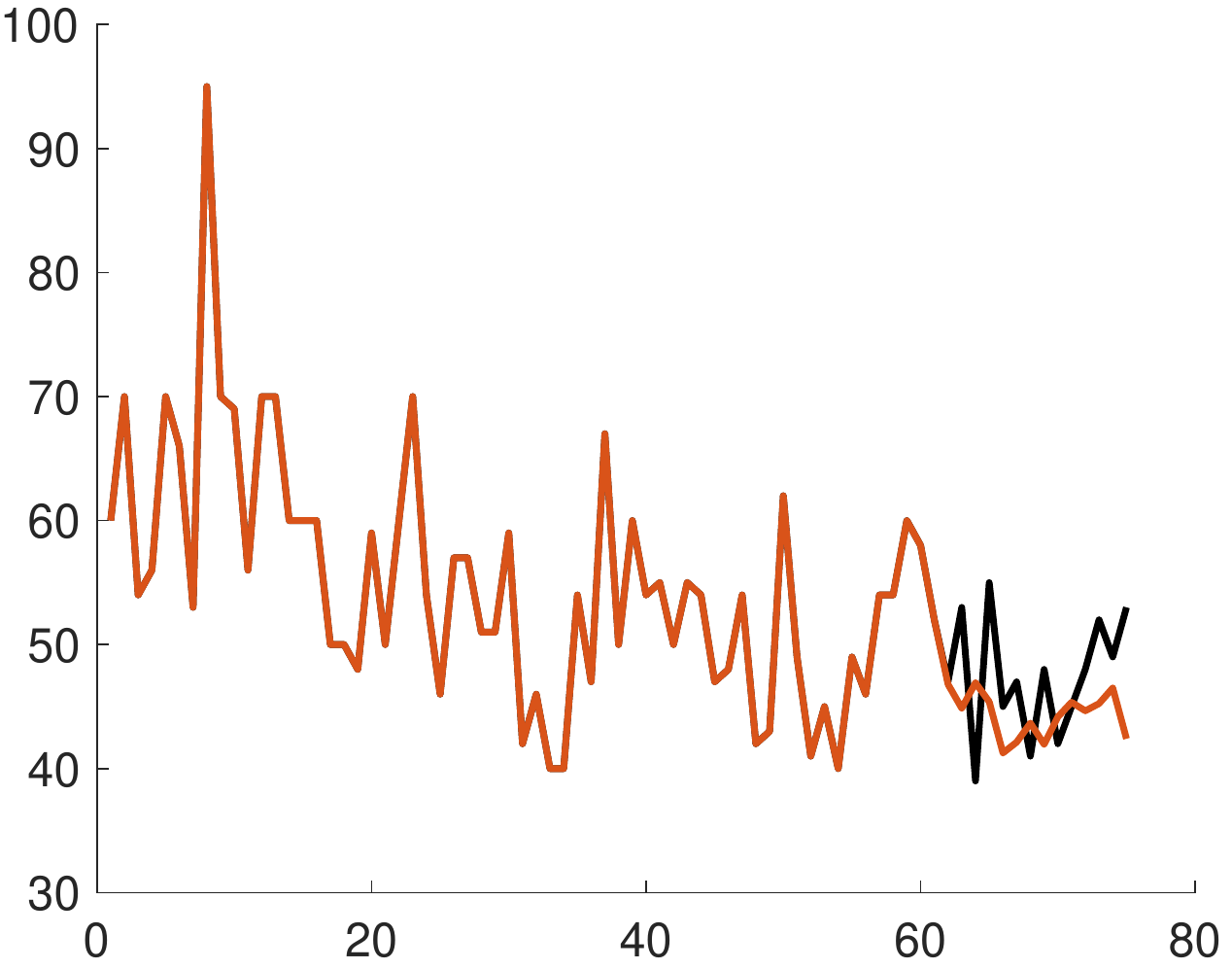}}
    \subfigure[$W_{3}$ with $a=0.01$ and $l=0.01$]{\includegraphics[width=0.16\textwidth]{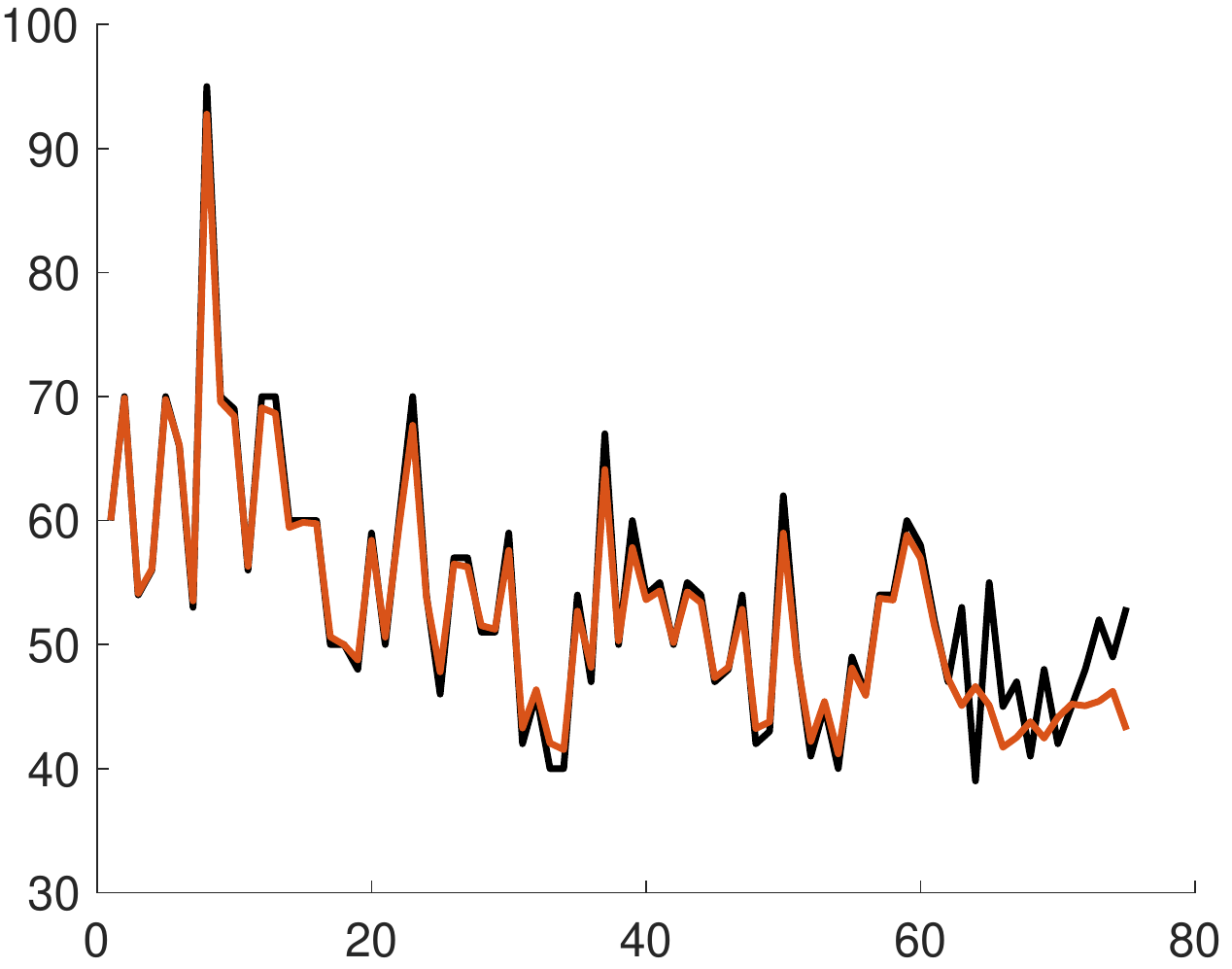}}
  \caption{Forecasts of the Cowtemp time series, with original data (black), for different weighting schemes}
\label{fig:fore}
\end{figure}

We solve \eqref{eq:regularisation3} with $\gamma=100$ and explore different weight vectors $W$. Figure \ref{fig:fore} contains plots of the forecasts of the Cowtemp time series for three weight vectors $W$. First, we give each observation unit weight, $W_{1}=(1,1,\ldots,1)$. Second we introduce the weights such that the weighted vector norm in \eqref{eq:regularisation3} becomes the Frobenius norm $\|\mathcal{S}(\bp_{(1:n)})-\mathcal{S}(\bp_{0})\|_{F} $. Finally we have the exponential weighing, where the $j$th element of $W_{3}$ is given by $a \exp(l j)$ for some scalars $a$ and $l$.

\begin{figure}[ht!]
\centering
    \subfigure[Surface plot]{\includegraphics[width=0.23\textwidth]{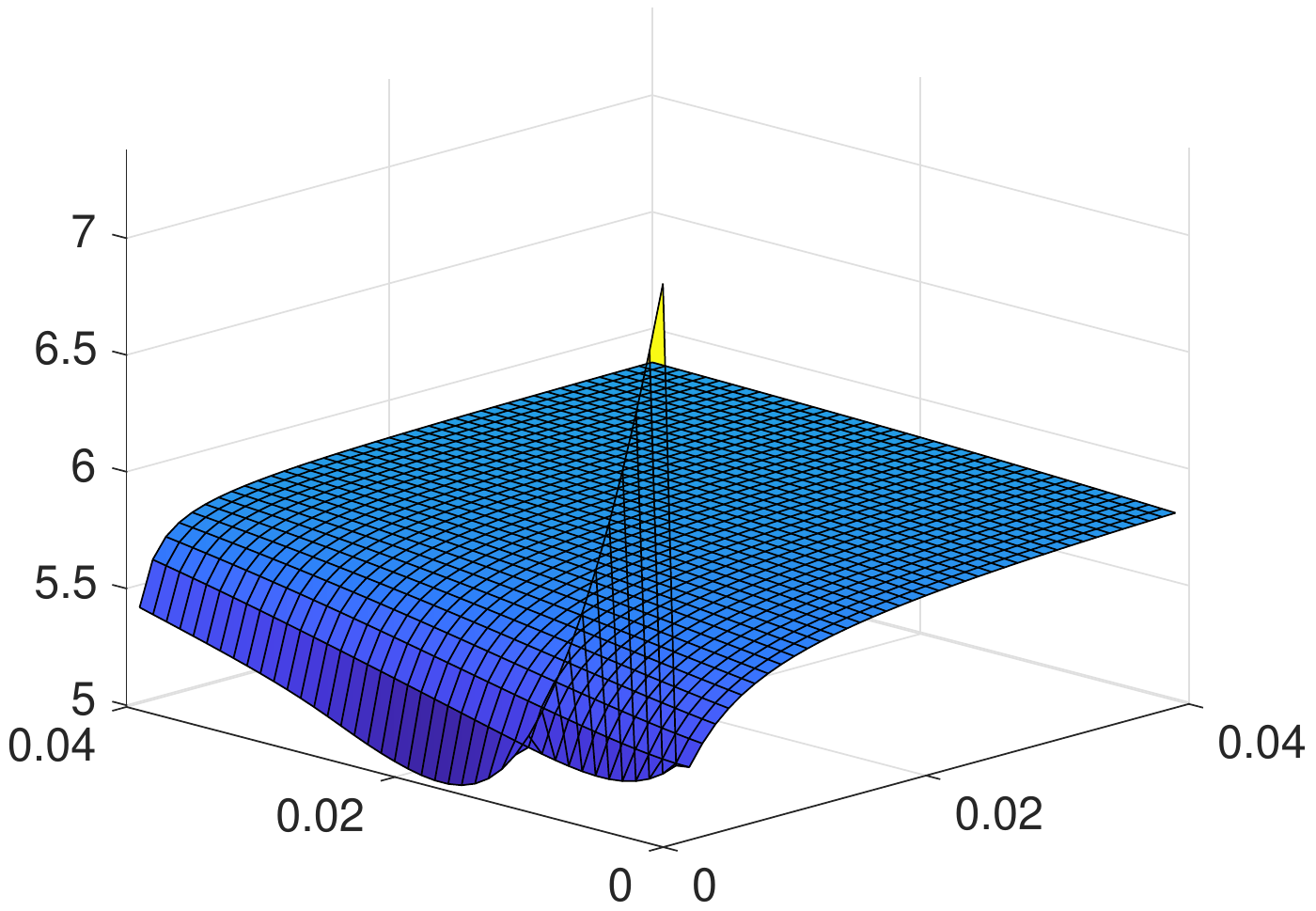}}
    \subfigure[Contour plot]{\includegraphics[width=0.23\textwidth]{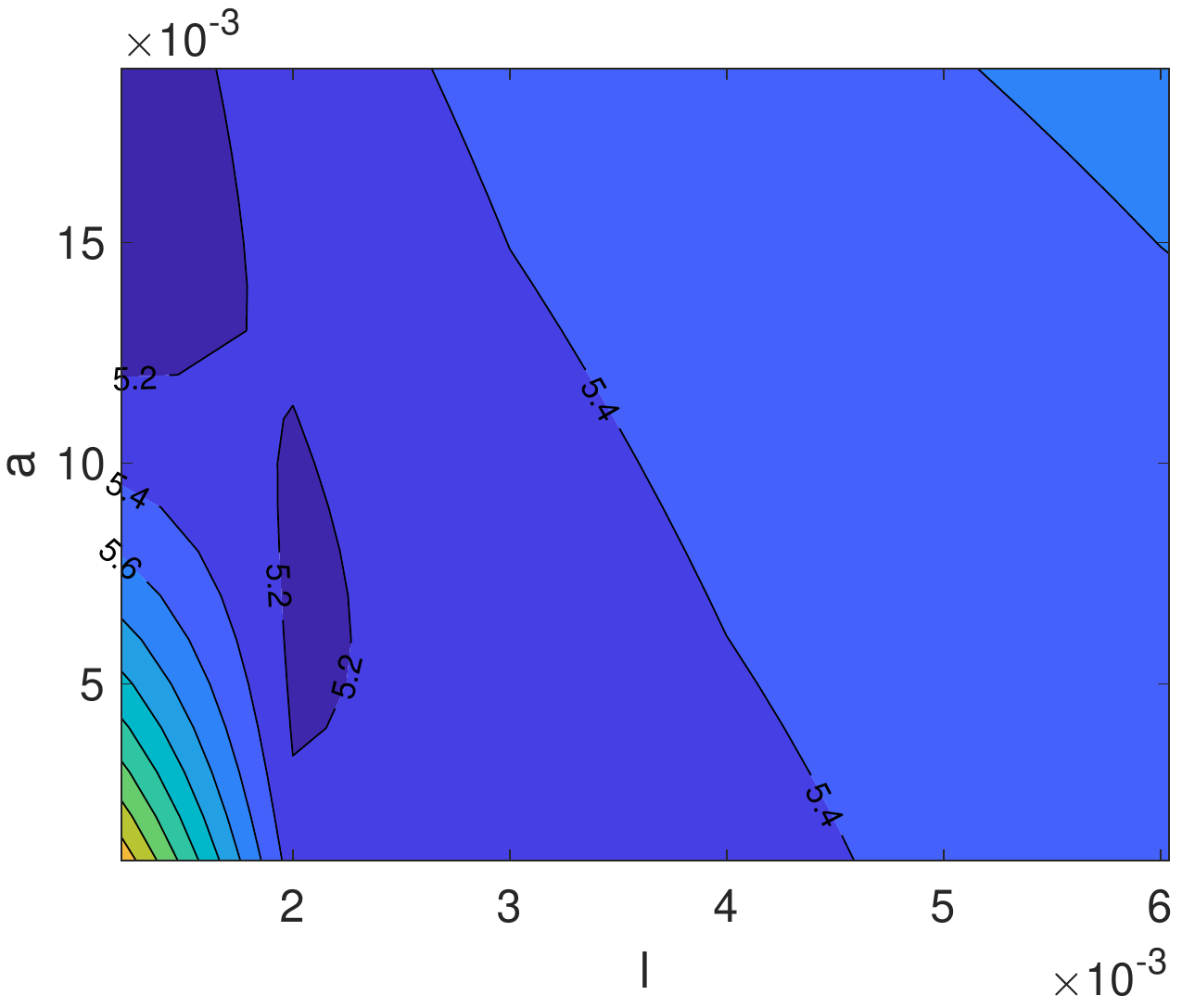}}
  \caption{Root mean square error of the forecast against $a$ and $l$}
\label{fig:rmse}
\end{figure}

\begin{figure}[ht!]
\centering
  \includegraphics[width=0.35\textwidth]{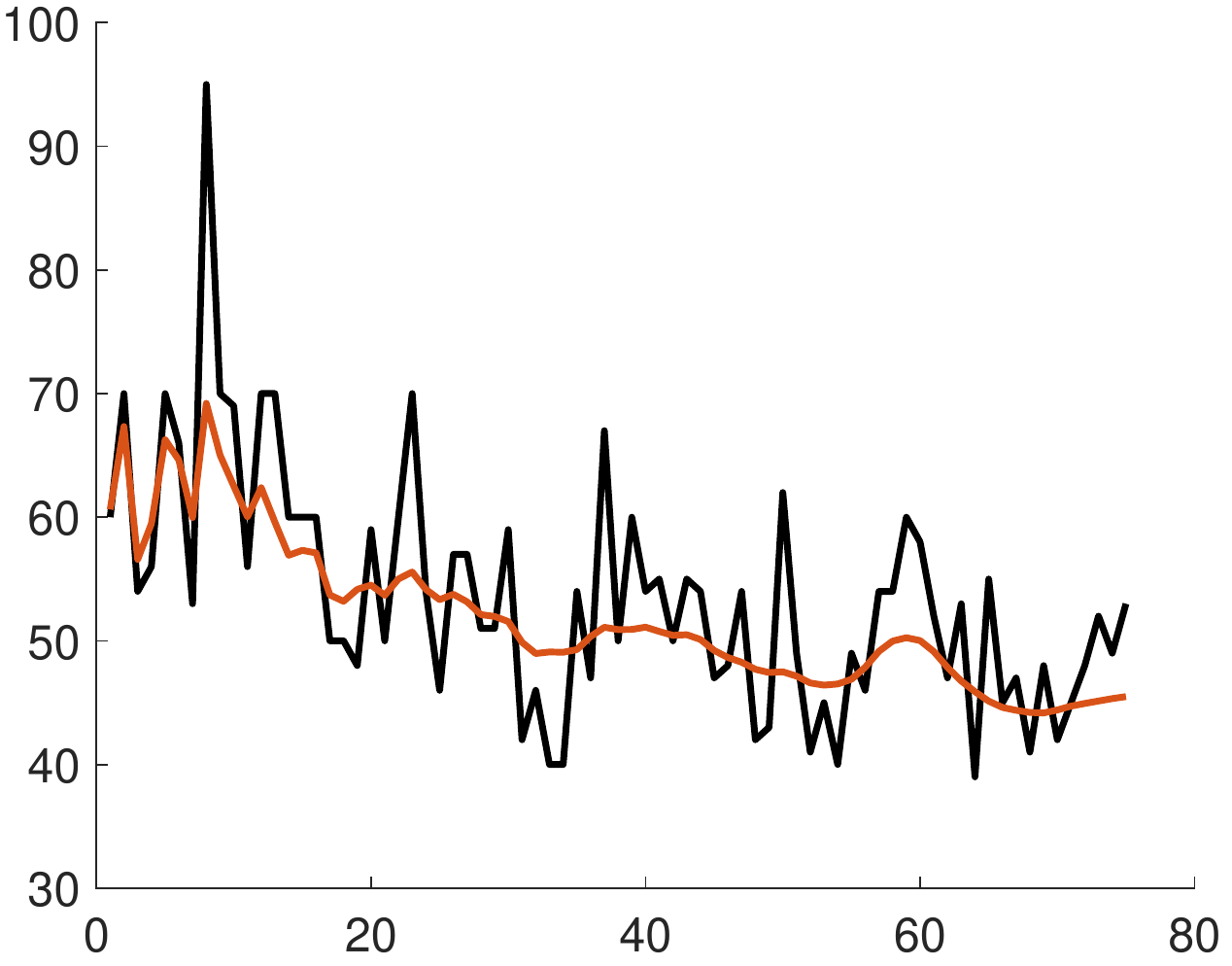}
  \caption{Forecast when $a=0.001$ and $l=0.017$ for weighting scheme $W_{3}$}
\label{fig:bestf}
\end{figure}

Figure \ref{fig:rmse} contains plots of the root mean square error (RMSE) of the forecast obtained using exponential weighting, that is \eqref{eq:regularisation3} with $W=W_3$, for different $a$ and $l$. Note that the approximately optimal RMSE found from \eqref{eq:regularisation3} is 4.9928 obtained when $a=0.001$ and $l=0.017$. A plot of the data with this forecast is given in Figure \ref{fig:bestf}, and it has approximating rank of 28, in this case. The RMSE obtained by SSA as described in \cite{golyandina2018singular} was found to be 5.253602.

\section{Conclusion}
In this paper we have attempted to give a review, commentary and bibliography of papers associated with Hankel low-rank approximation/completion, through the lens of time series modelling and forecasting. It is clear that this is an active area of pursuit, with both rich theoretical work, several algorithms and much possibility for practical application. Open questions remain around when the solution of the nuclear norm relaxation coincides with the `equivalent' rank minimization problem, and this is the motivator behind much present research. It remains to be the case to see such technology embraced by forecasting practitioners, but one would assume, that in time, the work described in this paper will be. Future big data applications on modelling several time series simultaneously are surely on the horizon.
\bibliographystyle{imsart-number}
\bibliography{refs-gillard_usevich}

\begin{thebibliography}{121}

\bibitem{abatzoglou1991tsp}
\begin{barticle}[author]
\bauthor{\bsnm{Abatzoglou},~\bfnm{T.~J.}\binits{T.~J.}},
  \bauthor{\bsnm{Mendel},~\bfnm{J.~M.}\binits{J.~M.}} \AND
  \bauthor{\bsnm{Harada},~\bfnm{G.~A.}\binits{G.~A.}}
(\byear{1991}).
\btitle{The constrained total least squares technique and its applications to
  harmonic superresolution}.
\bjournal{IEEE Transactions on Signal Processing}
\bvolume{39}
\bpages{1070-1087}.
\bdoi{10.1109/78.80955}
\end{barticle}
\endbibitem

\bibitem{adamjan1971analytic}
\begin{barticle}[author]
\bauthor{\bsnm{Adamjan},~\bfnm{Vadym~M}\binits{V.~M.}},
  \bauthor{\bsnm{Arov},~\bfnm{Damir~Z}\binits{D.~Z.}} \AND
  \bauthor{\bsnm{Kre{\u\i}n},~\bfnm{M~G}\binits{M.~G.}}
(\byear{1971}).
\btitle{Analytic properties of Schmidt pairs for a Hankel operator and the
  generalized Schur-Takagi problem}.
\bjournal{Mathematics of the USSR-Sbornik}
\bvolume{15}
\bpages{31}.
\end{barticle}
\endbibitem

\bibitem{akccay2014subspace}
\begin{barticle}[author]
\bauthor{\bsnm{Ak{\c{c}}ay},~\bfnm{H{\"u}seyin}\binits{H.}}
(\byear{2014}).
\btitle{Subspace-based spectrum estimation in frequency-domain by regularized
  nuclear norm minimization}.
\bjournal{Signal Processing}
\bvolume{99}
\bpages{69--85}.
\end{barticle}
\endbibitem

\bibitem{alshabili2019icassp}
\begin{binproceedings}[author]
\bauthor{\bsnm{{Al-Shabili}},~\bfnm{A.}\binits{A.}} \AND
  \bauthor{\bsnm{{Selesnick}},~\bfnm{I.}\binits{I.}}
(\byear{2019}).
\btitle{Sharpening Sparse Regularizers}.
In \bbooktitle{ICASSP 2019 - 2019 IEEE International Conference on Acoustics,
  Speech and Signal Processing (ICASSP)}
\bpages{4908-4912}.
\bdoi{10.1109/ICASSP.2019.8683039}
\end{binproceedings}
\endbibitem

\bibitem{andersson2013alternating}
\begin{barticle}[author]
\bauthor{\bsnm{Andersson},~\bfnm{Fredrik}\binits{F.}} \AND
  \bauthor{\bsnm{Carlsson},~\bfnm{Marcus}\binits{M.}}
(\byear{2013}).
\btitle{Alternating projections on nontangential manifolds}.
\bjournal{Constructive approximation}
\bvolume{38}
\bpages{489--525}.
\end{barticle}
\endbibitem

\bibitem{andersson2017convex}
\begin{barticle}[author]
\bauthor{\bsnm{Andersson},~\bfnm{Fredrik}\binits{F.}},
  \bauthor{\bsnm{Carlsson},~\bfnm{Marcus}\binits{M.}} \AND
  \bauthor{\bsnm{Olsson},~\bfnm{Carl}\binits{C.}}
(\byear{2017}).
\btitle{Convex envelopes for fixed rank approximation}.
\bjournal{Optimization Letters}
\bvolume{11}
\bpages{1783--1795}.
\end{barticle}
\endbibitem

\bibitem{andersson2014admm}
\begin{barticle}[author]
\bauthor{\bsnm{Andersson},~\bfnm{Fredrik}\binits{F.}},
  \bauthor{\bsnm{Carlsson},~\bfnm{Marcus}\binits{M.}},
  \bauthor{\bsnm{Tourneret},~\bfnm{Jean-Yves}\binits{J.-Y.}} \AND
  \bauthor{\bsnm{Wendt},~\bfnm{Herwig}\binits{H.}}
(\byear{2014}).
\btitle{A New Frequency Estimation Method for Equally and Unequally Spaced
  Data}.
\bjournal{IEEE Transactions on Signal Processing}
\bvolume{62}
\bpages{5761-5774}.
\bdoi{10.1109/TSP.2014.2358961}
\end{barticle}
\endbibitem

\bibitem{Antoulas1997}
\begin{binbook}[author]
\bauthor{\bsnm{Antoulas},~\bfnm{A.~C.}\binits{A.~C.}}
(\byear{1997}).
\btitle{On the Approximation of Hankel Matrices}
In \bbooktitle{Operators, Systems and Linear Algebra: Three Decades of
  Algebraic Systems Theory}
\bpages{17--22}.
\bpublisher{Vieweg+Teubner Verlag}, \baddress{Wiesbaden}.
\bdoi{10.1007/978-3-663-09823-2_2}
\end{binbook}
\endbibitem

\bibitem{balle2021optimal}
\begin{barticle}[author]
\bauthor{\bsnm{Balle},~\bfnm{Borja}\binits{B.}},
  \bauthor{\bsnm{Lacroce},~\bfnm{Clara}\binits{C.}},
  \bauthor{\bsnm{Panangaden},~\bfnm{Prakash}\binits{P.}},
  \bauthor{\bsnm{Precup},~\bfnm{Doina}\binits{D.}} \AND
  \bauthor{\bsnm{Rabusseau},~\bfnm{Guillaume}\binits{G.}}
(\byear{2021}).
\btitle{Optimal Spectral-Norm Approximate Minimization of Weighted Finite
  Automata}.
\bjournal{arXiv preprint arXiv:2102.06860}.
\end{barticle}
\endbibitem

\bibitem{beylkin2005exponential}
\begin{barticle}[author]
\bauthor{\bsnm{Beylkin},~\bfnm{Gregory}\binits{G.}} \AND
  \bauthor{\bsnm{Monz\'{o}n},~\bfnm{Lucas}\binits{L.}}
(\byear{2005}).
\btitle{On approximation of functions by exponential sums}.
\bjournal{Applied and Computational Harmonic Analysis}
\bvolume{19}
\bpages{17-48}.
\bdoi{https://doi.org/10.1016/j.acha.2005.01.003}
\end{barticle}
\endbibitem

\bibitem{blomberg2015regpaths}
\begin{barticle}[author]
\bauthor{\bsnm{Blomberg},~\bfnm{Niclas}\binits{N.}},
  \bauthor{\bsnm{Rojas},~\bfnm{Cristian~R.}\binits{C.~R.}} \AND
  \bauthor{\bsnm{Wahlberg},~\bfnm{Bo}\binits{B.}}
(\byear{2015}).
\btitle{Regularization Paths for Re-Weighted Nuclear Norm Minimization}.
\bjournal{IEEE Signal Processing Letters}
\bvolume{22}
\bpages{1980-1984}.
\bdoi{10.1109/LSP.2015.2450505}
\end{barticle}
\endbibitem

\bibitem{borsdorf2012structured}
\begin{bphdthesis}[author]
\bauthor{\bsnm{Borsdorf},~\bfnm{Ruediger}\binits{R.}}
(\byear{2012}).
\btitle{Structured Matrix Nearness Problems: Theory and Algorithms}
\btype{PhD thesis},
\bpublisher{University of Manchester}.
\end{bphdthesis}
\endbibitem

\bibitem{brachat2010tensors}
\begin{barticle}[author]
\bauthor{\bsnm{Brachat},~\bfnm{Jerome}\binits{J.}},
  \bauthor{\bsnm{Comon},~\bfnm{Pierre}\binits{P.}},
  \bauthor{\bsnm{Mourrain},~\bfnm{Bernard}\binits{B.}} \AND
  \bauthor{\bsnm{Tsigaridas},~\bfnm{Elias}\binits{E.}}
(\byear{2010}).
\btitle{Symmetric tensor decomposition}.
\bjournal{Linear Algebra and its Applications}
\bvolume{433}
\bpages{1851--1872}.
\bdoi{http://dx.doi.org/10.1016/j.laa.2010.06.046}
\end{barticle}
\endbibitem

\bibitem{brambilla2008alexander}
\begin{barticle}[author]
\bauthor{\bsnm{Brambilla},~\bfnm{Maria~Chiara}\binits{M.~C.}} \AND
  \bauthor{\bsnm{Ottaviani},~\bfnm{Giorgio}\binits{G.}}
(\byear{2008}).
\btitle{On the Alexander–Hirschowitz theorem}.
\bjournal{Journal of Pure and Applied Algebra}
\bvolume{212}
\bpages{1229-1251}.
\end{barticle}
\endbibitem

\bibitem{bresler1986exact}
\begin{barticle}[author]
\bauthor{\bsnm{Bresler},~\bfnm{Yoram}\binits{Y.}} \AND
  \bauthor{\bsnm{Macovski},~\bfnm{Albert}\binits{A.}}
(\byear{1986}).
\btitle{Exact maximum likelihood parameter estimation of superimposed
  exponential signals in noise}.
\bjournal{IEEE Transactions on Acoustics, Speech, and Signal Processing}
\bvolume{34}
\bpages{1081--1089}.
\end{barticle}
\endbibitem

\bibitem{butcher2016simple}
\begin{barticle}[author]
\bauthor{\bsnm{Butcher},~\bfnm{H.}\binits{H.}} \AND
  \bauthor{\bsnm{Gillard},~\bfnm{J.}\binits{J.}}
(\byear{2017}).
\btitle{Simple Nuclear Norm Based Algorithms for Imputing Missing Data and
  Forecasting in Time Series}.
\bjournal{Statistics and Its Interface}
\bvolume{10}
\bpages{19--25}.
\bdoi{10.4310/SII.2017.v10.n1.a2}
\end{barticle}
\endbibitem

\bibitem{cadzow1988signal}
\begin{barticle}[author]
\bauthor{\bsnm{Cadzow},~\bfnm{James~A}\binits{J.~A.}}
(\byear{1988}).
\btitle{Signal enhancement-a composite property mapping algorithm}.
\bjournal{Acoustics, Speech and Signal Processing, IEEE Transactions on}
\bvolume{36}
\bpages{49--62}.
\end{barticle}
\endbibitem

\bibitem{candes2013super}
\begin{barticle}[author]
\bauthor{\bsnm{Cand{\`e}s},~\bfnm{Emmanuel~J}\binits{E.~J.}} \AND
  \bauthor{\bsnm{Fernandez-Granda},~\bfnm{Carlos}\binits{C.}}
(\byear{2013}).
\btitle{Super-resolution from noisy data}.
\bjournal{Journal of Fourier Analysis and Applications}
\bvolume{19}
\bpages{1229--1254}.
\end{barticle}
\endbibitem

\bibitem{candes2010matrix}
\begin{barticle}[author]
\bauthor{\bsnm{Candes},~\bfnm{Emmanuel~J}\binits{E.~J.}} \AND
  \bauthor{\bsnm{Plan},~\bfnm{Yaniv}\binits{Y.}}
(\byear{2010}).
\btitle{Matrix completion with noise}.
\bjournal{Proceedings of the IEEE}
\bvolume{98}
\bpages{925--936}.
\bdoi{10.1109/JPROC.2009.2035722}
\end{barticle}
\endbibitem

\bibitem{Candes2009}
\begin{barticle}[author]
\bauthor{\bsnm{Cand\`{e}s},~\bfnm{Emmanuel~J.}\binits{E.~J.}} \AND
  \bauthor{\bsnm{Recht},~\bfnm{Benjamin}\binits{B.}}
(\byear{2009}).
\btitle{{Exact Matrix Completion via Convex Optimization}}.
\bjournal{Foundations of Computational Mathematics}
\bvolume{9}
\bpages{717--772}.
\bdoi{10.1007/s10208-009-9045-5}
\end{barticle}
\endbibitem

\bibitem{caprani2019icassp}
\begin{binproceedings}[author]
\bauthor{\bsnm{{Caprani}},~\bfnm{J.}\binits{J.}} \AND
  \bauthor{\bsnm{{Carlsson}},~\bfnm{M.}\binits{M.}}
(\byear{2019}).
\btitle{Quadratic Envelope Regularization for Structured Low Rank
  Approximation}.
In \bbooktitle{ICASSP 2019 - 2019 IEEE International Conference on Acoustics,
  Speech and Signal Processing (ICASSP)}
\bpages{8217-8221}.
\bdoi{10.1109/ICASSP.2019.8683420}
\end{binproceedings}
\endbibitem

\bibitem{chen2014robust}
\begin{barticle}[author]
\bauthor{\bsnm{Chen},~\bfnm{Y.}\binits{Y.}} \AND
  \bauthor{\bsnm{Chi},~\bfnm{Y.}\binits{Y.}}
(\byear{2014}).
\btitle{Robust Spectral Compressed Sensing via Structured Matrix Completion}.
\bjournal{IEEE Transactions on Information Theory}
\bvolume{60}
\bpages{6576-6601}.
\bdoi{10.1109/TIT.2014.2343623}
\end{barticle}
\endbibitem

\bibitem{chin1998agcd}
\begin{binproceedings}[author]
\bauthor{\bsnm{Chin},~\bfnm{Paulina}\binits{P.}},
  \bauthor{\bsnm{Corless},~\bfnm{Robert~M.}\binits{R.~M.}} \AND
  \bauthor{\bsnm{Corliss},~\bfnm{George~F.}\binits{G.~F.}}
(\byear{1998}).
\btitle{Optimization strategies for the approximate GCD problem}.
In \bbooktitle{Proceedings of ISSAC '98}
\bpages{228--235}.
\bpublisher{ACM}.
\end{binproceedings}
\endbibitem

\bibitem{chu2003slra}
\begin{barticle}[author]
\bauthor{\bsnm{Chu},~\bfnm{M.~T.}\binits{M.~T.}},
  \bauthor{\bsnm{Funderlic},~\bfnm{R.~E.}\binits{R.~E.}} \AND
  \bauthor{\bsnm{Plemmons},~\bfnm{L.~J.}\binits{L.~J.}}
(\byear{2003}).
\btitle{Structured low rank approximation}.
\bjournal{Linear Algebra Its Appl.}
\bvolume{366}
\bpages{157--172}.
\end{barticle}
\endbibitem

\bibitem{cifuentes2021convex}
\begin{barticle}[author]
\bauthor{\bsnm{Cifuentes},~\bfnm{Diego}\binits{D.}}
(\byear{2021}).
\btitle{A Convex Relaxation to Compute the Nearest Structured Rank Deficient
  Matrix}.
\bjournal{SIAM Journal on Matrix Analysis and Applications}
\bvolume{42}
\bpages{708-729}.
\bdoi{10.1137/19M1257640}
\end{barticle}
\endbibitem

\bibitem{comon1996sums}
\begin{barticle}[author]
\bauthor{\bsnm{Comon},~\bfnm{P.}\binits{P.}} \AND
  \bauthor{\bsnm{Mourrain},~\bfnm{B.}\binits{B.}}
(\byear{1996}).
\btitle{Decomposition of quantics in sums of powers of linear forms}.
\bjournal{Signal Processing}
\bvolume{53}
\bpages{93-107}.
\bnote{Higher Order Statistics}.
\end{barticle}
\endbibitem

\bibitem{condat2015cadzow}
\begin{barticle}[author]
\bauthor{\bsnm{Condat},~\bfnm{Laurent}\binits{L.}} \AND
  \bauthor{\bsnm{Hirabayashi},~\bfnm{Akira}\binits{A.}}
(\byear{2015}).
\btitle{Cadzow denoising upgraded: A new projection method for the recovery of
  Dirac pulses from noisy linear measurements}.
\bjournal{Sampling Theory in Signal and Image Processing}
\bvolume{14}
\bpages{17--47}.
\end{barticle}
\endbibitem

\bibitem{dai2015nuclear}
\begin{barticle}[author]
\bauthor{\bsnm{Dai},~\bfnm{Liang}\binits{L.}} \AND
  \bauthor{\bsnm{Pelckmans},~\bfnm{Kristiaan}\binits{K.}}
(\byear{2015}).
\btitle{On the nuclear norm heuristic for a {Hankel} matrix completion
  problem}.
\bjournal{Automatica}
\bvolume{51}
\bpages{268--272}.
\bdoi{10.1016/j.automatica.2014.10.045}
\end{barticle}
\endbibitem

\bibitem{de1993structured}
\begin{barticle}[author]
\bauthor{\bsnm{De~Moor},~\bfnm{Bart}\binits{B.}}
(\byear{1993}).
\btitle{Structured total least squares and L2 approximation problems}.
\bjournal{Linear algebra and its applications}
\bvolume{188}
\bpages{163--205}.
\end{barticle}
\endbibitem

\bibitem{dragotti2007fri}
\begin{barticle}[author]
\bauthor{\bsnm{Dragotti},~\bfnm{Pier~Luigi}\binits{P.~L.}},
  \bauthor{\bsnm{Vetterli},~\bfnm{Martin}\binits{M.}} \AND
  \bauthor{\bsnm{Blu},~\bfnm{Thierry}\binits{T.}}
(\byear{2007}).
\btitle{Sampling Moments and Reconstructing Signals of Finite Rate of
  Innovation: Shannon Meets Strang–Fix}.
\bjournal{IEEE Transactions on Signal Processing}
\bvolume{55}
\bpages{1741-1757}.
\bdoi{10.1109/TSP.2006.890907}
\end{barticle}
\endbibitem

\bibitem{elad2004shape}
\begin{barticle}[author]
\bauthor{\bsnm{Elad},~\bfnm{M.}\binits{M.}},
  \bauthor{\bsnm{Milanfar},~\bfnm{P.}\binits{P.}} \AND
  \bauthor{\bsnm{Golub},~\bfnm{G.~H.}\binits{G.~H.}}
(\byear{2004}).
\btitle{Shape from moments - an estimation theory perspective}.
\bjournal{IEEE Transactions on Signal Processing}
\bvolume{52}
\bpages{1814-1829}.
\bdoi{10.1109/TSP.2004.828919}
\end{barticle}
\endbibitem

\bibitem{fazel2002matrix}
\begin{bphdthesis}[author]
\bauthor{\bsnm{Fazel},~\bfnm{Maryam}\binits{M.}}
(\byear{2002}).
\btitle{Matrix rank minimization with applications}
\btype{PhD thesis},
\bpublisher{PhD thesis, Stanford University}.
\end{bphdthesis}
\endbibitem

\bibitem{fazel2001rank}
\begin{binproceedings}[author]
\bauthor{\bsnm{Fazel},~\bfnm{Maryam}\binits{M.}},
  \bauthor{\bsnm{Hindi},~\bfnm{Haitham}\binits{H.}} \AND
  \bauthor{\bsnm{Boyd},~\bfnm{Stephen~P}\binits{S.~P.}}
(\byear{2001}).
\btitle{A rank minimization heuristic with application to minimum order system
  approximation}.
In \bbooktitle{Proceedings of the 2001 American Control Conference.(Cat. No.
  01CH37148)}
\bvolume{6}
\bpages{4734--4739}.
\bpublisher{IEEE}.
\end{binproceedings}
\endbibitem

\bibitem{fazel2013hankel}
\begin{barticle}[author]
\bauthor{\bsnm{Fazel},~\bfnm{Maryam}\binits{M.}},
  \bauthor{\bsnm{Pong},~\bfnm{Ting~Kei}\binits{T.~K.}},
  \bauthor{\bsnm{Sun},~\bfnm{Defeng}\binits{D.}} \AND
  \bauthor{\bsnm{Tseng},~\bfnm{Paul}\binits{P.}}
(\byear{2013}).
\btitle{Hankel Matrix Rank Minimization with Applications to System
  Identification and Realization}.
\bjournal{SIAM Journal on Matrix Analysis and Applications}
\bvolume{34}
\bpages{946-977}.
\bdoi{10.1137/110853996}
\end{barticle}
\endbibitem

\bibitem{fazzi2018agcd}
\begin{barticle}[author]
\bauthor{\bsnm{Fazzi},~\bfnm{A.}\binits{A.}},
  \bauthor{\bsnm{Guglielmi},~\bfnm{N.}\binits{N.}} \AND
  \bauthor{\bsnm{Markovsky},~\bfnm{I.}\binits{I.}}
(\byear{2018}).
\btitle{An {ODE} based method for computing the Approximate Greatest Common
  Divisor of polynomials}.
\bjournal{Numerical algorithms}
\bvolume{81}
\bpages{719--740}.
\bdoi{10.1007/s11075-018-0569-0}
\end{barticle}
\endbibitem

\bibitem{fazzi2021ode}
\begin{barticle}[author]
\bauthor{\bsnm{Fazzi},~\bfnm{A.}\binits{A.}},
  \bauthor{\bsnm{Guglielmi},~\bfnm{N.}\binits{N.}} \AND
  \bauthor{\bsnm{Markovsky},~\bfnm{I.}\binits{I.}}
(\byear{2021}).
\btitle{A gradient system approach for {H}ankel structured low-rank
  approximation}.
\bjournal{Linear Algebra Appl.}
\bdoi{10.1016/j.laa.2020.11.016}
\end{barticle}
\endbibitem

\bibitem{feldmann1999parameterization}
\begin{barticle}[author]
\bauthor{\bsnm{Feldmann},~\bfnm{Sven}\binits{S.}} \AND
  \bauthor{\bsnm{Heinig},~\bfnm{Georg}\binits{G.}}
(\byear{1999}).
\btitle{Parametrization of minimal rank block Hankel matrix extensions and
  minimal partial realizations}.
\bjournal{Integral Equations and Operator Theory}
\bvolume{33}
\bpages{151--171}.
\end{barticle}
\endbibitem

\bibitem{floater2021best}
\begin{barticle}[author]
\bauthor{\bsnm{Floater},~\bfnm{Michael~S}\binits{M.~S.}},
  \bauthor{\bsnm{Manni},~\bfnm{Carla}\binits{C.}},
  \bauthor{\bsnm{Sande},~\bfnm{Espen}\binits{E.}} \AND
  \bauthor{\bsnm{Speleers},~\bfnm{Hendrik}\binits{H.}}
(\byear{2021}).
\btitle{Best low-rank approximations and {Kolmogorov} n-widths}.
\bjournal{SIAM Journal on Matrix Analysis and Applications}
\bvolume{42}
\bpages{330--350}.
\end{barticle}
\endbibitem

\bibitem{gillard2010cadzow}
\begin{barticle}[author]
\bauthor{\bsnm{Gillard},~\bfnm{Jonathan}\binits{J.}}
(\byear{2010}).
\btitle{Cadzow’s basic algorithm, alternating projections and singular
  spectrum analysis}.
\bjournal{Statistics and its Interface}
\bvolume{3}
\bpages{335--343}.
\end{barticle}
\endbibitem

\bibitem{gillard2018structured}
\begin{barticle}[author]
\bauthor{\bsnm{Gillard},~\bfnm{Jonathan}\binits{J.}} \AND
  \bauthor{\bsnm{Usevich},~\bfnm{Konstantin}\binits{K.}}
(\byear{2018}).
\btitle{Structured low-rank matrix completion for forecasting in time series
  analysis}.
\bjournal{International Journal of Forecasting}
\bvolume{34}
\bpages{582--597}.
\end{barticle}
\endbibitem

\bibitem{gillard2013optimization}
\begin{barticle}[author]
\bauthor{\bsnm{Gillard},~\bfnm{Jonathan}\binits{J.}} \AND
  \bauthor{\bsnm{Zhigljavsky},~\bfnm{Anatoly}\binits{A.}}
(\byear{2013}).
\btitle{Optimization challenges in the structured low rank approximation
  problem}.
\bjournal{Journal of Global Optimization}
\bvolume{57}
\bpages{733--751}.
\bdoi{10.1007/s10898-012-9962-8}
\end{barticle}
\endbibitem

\bibitem{gillard2015stochastic}
\begin{barticle}[author]
\bauthor{\bsnm{Gillard},~\bfnm{Jonathan}\binits{J.}} \AND
  \bauthor{\bsnm{Zhigljavsky},~\bfnm{Anatoly}\binits{A.}}
(\byear{2015}).
\btitle{Stochastic algorithms for solving structured low-rank matrix
  approximation problems}.
\bjournal{Communications in Nonlinear Science and Numerical Simulation}
\bvolume{21}
\bpages{70--88}.
\end{barticle}
\endbibitem

\bibitem{gillard2016weighted}
\begin{barticle}[author]
\bauthor{\bsnm{Gillard},~\bfnm{JW}\binits{J.}} \AND
  \bauthor{\bsnm{Zhigljavsky},~\bfnm{Anatoly~A}\binits{A.~A.}}
(\byear{2016}).
\btitle{Weighted norms in subspace-based methods for time series analysis}.
\bjournal{Numerical Linear Algebra with Applications}
\bvolume{23}
\bpages{947--967}.
\end{barticle}
\endbibitem

\bibitem{gillis2011low}
\begin{barticle}[author]
\bauthor{\bsnm{Gillis},~\bfnm{Nicolas}\binits{N.}} \AND
  \bauthor{\bsnm{Glineur},~\bfnm{Fran{\c{c}}ois}\binits{F.}}
(\byear{2011}).
\btitle{Low-rank matrix approximation with weights or missing data is
  {NP}-hard}.
\bjournal{SIAM Journal on Matrix Analysis and Applications}
\bvolume{32}
\bpages{1149--1165}.
\bdoi{10.1137/110820361}
\end{barticle}
\endbibitem

\bibitem{golub2012matrix}
\begin{bbook}[author]
\bauthor{\bsnm{Golub},~\bfnm{Gene~H}\binits{G.~H.}} \AND
  \bauthor{\bsnm{Van~Loan},~\bfnm{Charles~F}\binits{C.~F.}}
(\byear{2012}).
\btitle{Matrix computations}
\bvolume{3}.
\bpublisher{JHU Press}.
\end{bbook}
\endbibitem

\bibitem{golyandina2010choice}
\begin{barticle}[author]
\bauthor{\bsnm{Golyandina},~\bfnm{Nina}\binits{N.}}
(\byear{2010}).
\btitle{On the choice of parameters in Singular Spectrum Analysis and related
  subspace-based methods}.
\bjournal{Statistics and its Interface}
\bvolume{3}.
\end{barticle}
\endbibitem

\bibitem{golyandina2015jss}
\begin{barticle}[author]
\bauthor{\bsnm{Golyandina},~\bfnm{Nina}\binits{N.}},
  \bauthor{\bsnm{Korobeynikov},~\bfnm{Anton}\binits{A.}},
  \bauthor{\bsnm{Shlemov},~\bfnm{Alex}\binits{A.}} \AND
  \bauthor{\bsnm{Usevich},~\bfnm{Konstantin}\binits{K.}}
(\byear{2015}).
\btitle{Multivariate and 2D Extensions of Singular Spectrum Analysis with the
  Rssa Package}.
\bjournal{Journal of Statistical Software, Articles}
\bvolume{67}
\bpages{1--78}.
\bdoi{10.18637/jss.v067.i02}
\end{barticle}
\endbibitem

\bibitem{golyandina2018singular}
\begin{bbook}[author]
\bauthor{\bsnm{Golyandina},~\bfnm{Nina}\binits{N.}},
  \bauthor{\bsnm{Korobeynikov},~\bfnm{Anton}\binits{A.}} \AND
  \bauthor{\bsnm{Zhigljavsky},~\bfnm{Anatoly}\binits{A.}}
(\byear{2018}).
\btitle{Singular spectrum analysis with R}.
\bpublisher{Springer}.
\end{bbook}
\endbibitem

\bibitem{golyandina2001analysis}
\begin{bbook}[author]
\bauthor{\bsnm{Golyandina},~\bfnm{Nina}\binits{N.}},
  \bauthor{\bsnm{Nekrutkin},~\bfnm{Vladimir}\binits{V.}} \AND
  \bauthor{\bsnm{Zhigljavsky},~\bfnm{Anatoly~A}\binits{A.~A.}}
(\byear{2001}).
\btitle{Analysis of time series structure: SSA and related techniques}.
\bpublisher{CRC press}.
\end{bbook}
\endbibitem

\bibitem{golyandina20102dssa}
\begin{bincollection}[author]
\bauthor{\bsnm{Golyandina},~\bfnm{N.}\binits{N.}} \AND
  \bauthor{\bsnm{Usevich},~\bfnm{K.}\binits{K.}}
(\byear{2010}).
\btitle{{2D-extension of Singular Spectrum Analysis}: algorithm and elements of
  theory}.
In \bbooktitle{Matrix Methods: Theory, Algorithms and Applications}
(\beditor{\bfnm{V.}\binits{V.}~\bsnm{Olshevsky}} \AND
  \beditor{\bfnm{E.}\binits{E.}~\bsnm{Tyrtyshnikov}}, eds.)
\bpages{449--474}.
\bpublisher{World Scientific}.
\end{bincollection}
\endbibitem

\bibitem{golyandina2020singular}
\begin{bbook}[author]
\bauthor{\bsnm{Golyandina},~\bfnm{Nina}\binits{N.}} \AND
  \bauthor{\bsnm{Zhigljavsky},~\bfnm{Anatoly}\binits{A.}}
(\byear{2020}).
\btitle{Singular Spectrum Analysis for time series},
\bedition{2nd} ed.
\bpublisher{Springer}.
\end{bbook}
\endbibitem

\bibitem{gonnet2013robust}
\begin{barticle}[author]
\bauthor{\bsnm{{Gonnet}},~\bfnm{Pedro}\binits{P.}},
  \bauthor{\bsnm{{G\"{u}ttel}},~\bfnm{Stefan}\binits{S.}} \AND
  \bauthor{\bsnm{{Trefethen}},~\bfnm{Lloyd~N.}\binits{L.~N.}}
(\byear{2013}).
\btitle{Robust Pad\'{e} Approximation via SVD}.
\bjournal{Siam Review}
\bvolume{55}
\bpages{101--117}.
\end{barticle}
\endbibitem

\bibitem{grussler2018convex}
\begin{barticle}[author]
\bauthor{\bsnm{Grussler},~\bfnm{Christian}\binits{C.}},
  \bauthor{\bsnm{Rantzer},~\bfnm{Anders}\binits{A.}} \AND
  \bauthor{\bsnm{Giselsson},~\bfnm{Pontus}\binits{P.}}
(\byear{2018}).
\btitle{Low-Rank Optimization With Convex Constraints}.
\bjournal{IEEE Transactions on Automatic Control}
\bvolume{63}
\bpages{4000-4007}.
\bdoi{10.1109/TAC.2018.2813009}
\end{barticle}
\endbibitem

\bibitem{guglielmi2017ode}
\begin{barticle}[author]
\bauthor{\bsnm{Guglielmi},~\bfnm{N.}\binits{N.}} \AND
  \bauthor{\bsnm{Markovsky},~\bfnm{I.}\binits{I.}}
(\byear{2017}).
\btitle{An {ODE} based method for computing the distance of co-prime
  polynomials to common divisibility}.
\bjournal{SIAM Journal on Numerical Analysis}
\bvolume{55}
\bpages{1456--1482}.
\bdoi{10.1137/15M1018265}
\end{barticle}
\endbibitem

\bibitem{hage2015robust}
\begin{binproceedings}[author]
\bauthor{\bsnm{Hage},~\bfnm{Clemens}\binits{C.}} \AND
  \bauthor{\bsnm{Kleinsteuber},~\bfnm{Martin}\binits{M.}}
(\byear{2015}).
\btitle{Robust Structured Low-Rank Approximation on the Grassmannian}.
In \bbooktitle{International Conference on Latent Variable Analysis and Signal
  Separation}
\bpages{295--303}.
\bpublisher{Springer}.
\end{binproceedings}
\endbibitem

\bibitem{hall1967combinatorial}
\begin{bbook}[author]
\bauthor{\bsnm{Hall},~\bfnm{M.}\binits{M.}}
(\byear{1967}).
\btitle{Combinatorial Theory}.
\bpublisher{Blaisdell Publishing Company}.
\end{bbook}
\endbibitem

\bibitem{heinig1995mosaic}
\begin{barticle}[author]
\bauthor{\bsnm{Heinig},~\bfnm{G.}\binits{G.}}
(\byear{1995}).
\btitle{Generalized inverses of {Hankel} and {Toeplitz} mosaic matrices}.
\bjournal{Linear Algebra Appl.}
\bvolume{216}
\bpages{43--59}.
\end{barticle}
\endbibitem

\bibitem{heinig1992block}
\begin{barticle}[author]
\bauthor{\bsnm{Heinig},~\bfnm{Georg}\binits{G.}} \AND
  \bauthor{\bsnm{Jankowski},~\bfnm{Peter}\binits{P.}}
(\byear{1992}).
\btitle{Kernel structure of block Hankel and Toeplitz matrices and partial
  realization}.
\bjournal{Linear Algebra and its Applications}
\bvolume{175}
\bpages{1 - 30}.
\bdoi{DOI: 10.1016/0024-3795(92)90299-P}
\end{barticle}
\endbibitem

\bibitem{heinig1984algebraic}
\begin{bbook}[author]
\bauthor{\bsnm{Heinig},~\bfnm{G.}\binits{G.}} \AND
  \bauthor{\bsnm{Rost},~\bfnm{K.}\binits{K.}}
(\byear{1984}).
\btitle{Algebraic methods for Toeplitz-like matrices and operators}.
\bpublisher{Birkh{\"a}user, Boston}.
\end{bbook}
\endbibitem

\bibitem{ho1966}
\begin{barticle}[author]
\bauthor{\bsnm{Ho},~\bfnm{B.~L.}\binits{B.~L.}} \AND
  \bauthor{\bsnm{Kalman},~\bfnm{R.~E.}\binits{R.~E.}}
(\byear{1966}).
\btitle{Editorial: Effective construction of linear state-variable models from
  input/output functions}.
\bjournal{Regelungstechnik}
\bvolume{14}
\bpages{545--548}.
\end{barticle}
\endbibitem

\bibitem{iarobbino1999power}
\begin{bbook}[author]
\bauthor{\bsnm{Iarobbino},~\bfnm{A.}\binits{A.}} \AND
  \bauthor{\bsnm{Kanev},~\bfnm{V.}\binits{V.}}
(\byear{1999}).
\btitle{Power sums, Gorenstein Algebras and Determinantal Loci}.
\bseries{Lecture Notes in Mathematics}
\bvolume{1721}.
\bpublisher{Springer}.
\end{bbook}
\endbibitem

\bibitem{iohvidov1982toeplitz}
\begin{bbook}[author]
\bauthor{\bsnm{Iohvidov},~\bfnm{IS}\binits{I.}}
(\byear{1982}).
\btitle{{Toeplitz} and {Hankel} matrices and forms. Algebraic theory}.
\bpublisher{Birkh{\"a}user, Basel}.
\end{bbook}
\endbibitem

\bibitem{ishteva2014factorization}
\begin{barticle}[author]
\bauthor{\bsnm{Ishteva},~\bfnm{Mariya}\binits{M.}},
  \bauthor{\bsnm{Usevich},~\bfnm{Konstantin}\binits{K.}} \AND
  \bauthor{\bsnm{Markovsky},~\bfnm{Ivan}\binits{I.}}
(\byear{2014}).
\btitle{Factorization approach to structured low-rank approximation with
  applications}.
\bjournal{SIAM Journal on Matrix Analysis and Applications}
\bvolume{35}
\bpages{1180--1204}.
\end{barticle}
\endbibitem

\bibitem{jia2016geophysics}
\begin{barticle}[author]
\bauthor{\bsnm{Jia},~\bfnm{Yongna}\binits{Y.}},
  \bauthor{\bsnm{Yu},~\bfnm{Siwei}\binits{S.}},
  \bauthor{\bsnm{Liu},~\bfnm{Lina}\binits{L.}} \AND
  \bauthor{\bsnm{Ma},~\bfnm{Jianwei}\binits{J.}}
(\byear{2016}).
\btitle{A fast rank-reduction algorithm for three-dimensional seismic data
  interpolation}.
\bjournal{Journal of Applied Geophysics}
\bvolume{132}
\bpages{137-145}.
\bdoi{https://doi.org/10.1016/j.jappgeo.2016.06.010}
\end{barticle}
\endbibitem

\bibitem{kang2016top}
\begin{binproceedings}[author]
\bauthor{\bsnm{Kang},~\bfnm{Zhao}\binits{Z.}},
  \bauthor{\bsnm{Peng},~\bfnm{Chong}\binits{C.}} \AND
  \bauthor{\bsnm{Cheng},~\bfnm{Qiang}\binits{Q.}}
(\byear{2016}).
\btitle{Top-n recommender system via matrix completion}.
In \bbooktitle{Proceedings of the AAAI Conference on Artificial Intelligence}
\bvolume{30}.
\end{binproceedings}
\endbibitem

\bibitem{atikur2013note}
\begin{barticle}[author]
\bauthor{\bsnm{Khan},~\bfnm{M.~A.~R.}\binits{M.~A.~R.}} \AND
  \bauthor{\bsnm{Poskitt},~\bfnm{D.~S.}\binits{D.~S.}}
(\byear{2013}).
\btitle{A note on window length selection in singular spectrum analysis}.
\bjournal{Australian \& New Zealand Journal of Statistics}
\bvolume{55}
\bpages{87--108}.
\end{barticle}
\endbibitem

\bibitem{khrulkov2018desingularization}
\begin{barticle}[author]
\bauthor{\bsnm{Khrulkov},~\bfnm{Valentin}\binits{V.}} \AND
  \bauthor{\bsnm{Oseledets},~\bfnm{Ivan}\binits{I.}}
(\byear{2018}).
\btitle{Desingularization of Bounded-Rank Matrix Sets}.
\bjournal{SIAM Journal on Matrix Analysis and Applications}
\bvolume{39}
\bpages{451-471}.
\bdoi{10.1137/16M1108194}
\end{barticle}
\endbibitem

\bibitem{knirsch2021optimal}
\begin{bmisc}[author]
\bauthor{\bsnm{Knirsch},~\bfnm{Hanna}\binits{H.}},
  \bauthor{\bsnm{Petz},~\bfnm{Markus}\binits{M.}} \AND
  \bauthor{\bsnm{Plonka},~\bfnm{Gerlind}\binits{G.}}
(\byear{2021}).
\btitle{Optimal Rank-1 Hankel Approximation of Matrices: Frobenius Norm,
  Spectral Norm, and Cadzow's Algorithm}.
\end{bmisc}
\endbibitem

\bibitem{korobeynikov2010computation}
\begin{barticle}[author]
\bauthor{\bsnm{Korobeynikov},~\bfnm{Anton}\binits{A.}}
(\byear{2010}).
\btitle{Computation-and space-efficient implementation of SSA}.
\bjournal{Statistics and Its Interface}
\bvolume{3}
\bpages{357--368}.
\end{barticle}
\endbibitem

\bibitem{649704}
\begin{barticle}[author]
\bauthor{\bsnm{Kuijper},~\bfnm{M.}\binits{M.}} \AND
  \bauthor{\bsnm{Willems},~\bfnm{J.~C.}\binits{J.~C.}}
(\byear{1997}).
\btitle{On constructing a shortest linear recurrence relation}.
\bjournal{IEEE Transactions on Automatic Control}
\bvolume{42}
\bpages{1554-1558}.
\bdoi{10.1109/9.649704}
\end{barticle}
\endbibitem

\bibitem{kulis2006learning}
\begin{binproceedings}[author]
\bauthor{\bsnm{Kulis},~\bfnm{Brian}\binits{B.}},
  \bauthor{\bsnm{Sustik},~\bfnm{M{\'a}ty{\'a}s}\binits{M.}} \AND
  \bauthor{\bsnm{Dhillon},~\bfnm{Inderjit}\binits{I.}}
(\byear{2006}).
\btitle{Learning low-rank kernel matrices}.
In \bbooktitle{Proceedings of the 23rd international conference on Machine
  learning}
\bpages{505--512}.
\end{binproceedings}
\endbibitem

\bibitem{kung1978identification}
\begin{binproceedings}[author]
\bauthor{\bsnm{Kung},~\bfnm{S.~Y.}\binits{S.~Y.}}
(\byear{1978}).
\btitle{A new identification and model reduction algorithm via singular value
  decomposition}.
In \bbooktitle{Proceedings of 12th Asilomar Conference on Circuits, Systems and
  Computers}
\bpages{705--714}.
\end{binproceedings}
\endbibitem

\bibitem{kurakin1995linear}
\begin{barticle}[author]
\bauthor{\bsnm{Kurakin},~\bfnm{V~L}\binits{V.~L.}},
  \bauthor{\bsnm{Kuzmin},~\bfnm{A~S}\binits{A.~S.}},
  \bauthor{\bsnm{Mikhalev},~\bfnm{A~V}\binits{A.~V.}} \AND
  \bauthor{\bsnm{Nechaev},~\bfnm{A~A}\binits{A.~A.}}
(\byear{1995}).
\btitle{Linear recurring sequences over rings and modules}.
\bjournal{Journal of mathematical sciences}
\bvolume{76}
\bpages{2793--2915}.
\end{barticle}
\endbibitem

\bibitem{kummerle2019phd}
\begin{bphdthesis}[author]
\bauthor{\bsnm{Kümmerle},~\bfnm{Christian}\binits{C.}}
(\byear{2019}).
\btitle{Understanding and Enhancing Data Recovery Algorithms}
\btype{Dissertation},
\bpublisher{Technische Universität München},
\baddress{München}.
\end{bphdthesis}
\endbibitem

\bibitem{kummerle2019irls}
\begin{binproceedings}[author]
\bauthor{\bsnm{Kümmerle},~\bfnm{Christian}\binits{C.}} \AND
  \bauthor{\bsnm{Verdun},~\bfnm{Claudio~M.}\binits{C.~M.}}
(\byear{2019}).
\btitle{Completion of Structured Low-Rank Matrices via Iteratively Reweighted
  Least Squares}.
In \bbooktitle{2019 13th International conference on Sampling Theory and
  Applications (SampTA)}
\bpages{1-5}.
\bdoi{10.1109/SampTA45681.2019.9030959}
\end{binproceedings}
\endbibitem

\bibitem{lemmerling2001analysis}
\begin{barticle}[author]
\bauthor{\bsnm{Lemmerling},~\bfnm{Philippe}\binits{P.}} \AND
  \bauthor{\bsnm{Van~Huffel},~\bfnm{Sabine}\binits{S.}}
(\byear{2001}).
\btitle{Analysis of the structured total least squares problem for
  Hankel/Toeplitz matrices}.
\bjournal{Numerical algorithms}
\bvolume{27}
\bpages{89--114}.
\end{barticle}
\endbibitem

\bibitem{lewis2008ap}
\begin{barticle}[author]
\bauthor{\bsnm{Lewis},~\bfnm{Adrian~S.}\binits{A.~S.}} \AND
  \bauthor{\bsnm{Malick},~\bfnm{Jérôme}\binits{J.}}
(\byear{2008}).
\btitle{Alternating Projections on Manifolds}.
\bjournal{Mathematics of Operations Research}
\bvolume{33}
\bpages{216--234}.
\end{barticle}
\endbibitem

\bibitem{liu2020simax}
\begin{barticle}[author]
\bauthor{\bsnm{Liu},~\bfnm{Tianxiang}\binits{T.}},
  \bauthor{\bsnm{Markovsky},~\bfnm{Ivan}\binits{I.}},
  \bauthor{\bsnm{Pong},~\bfnm{Ting~Kei}\binits{T.~K.}} \AND
  \bauthor{\bsnm{Takeda},~\bfnm{Akiko}\binits{A.}}
(\byear{2020}).
\btitle{A Hybrid Penalty Method for a Class of Optimization Problems with
  Multiple Rank Constraints}.
\bjournal{SIAM Journal on Matrix Analysis and Applications}
\bvolume{41}
\bpages{1260-1283}.
\bdoi{10.1137/19M1269919}
\end{barticle}
\endbibitem

\bibitem{liu2010interior}
\begin{barticle}[author]
\bauthor{\bsnm{Liu},~\bfnm{Zhang}\binits{Z.}} \AND
  \bauthor{\bsnm{Vandenberghe},~\bfnm{Lieven}\binits{L.}}
(\byear{2010}).
\btitle{Interior-point method for nuclear norm approximation with application
  to system identification}.
\bjournal{SIAM Journal on Matrix Analysis and Applications}
\bvolume{31}
\bpages{1235--1256}.
\end{barticle}
\endbibitem

\bibitem{markovsky2008structured}
\begin{barticle}[author]
\bauthor{\bsnm{Markovsky},~\bfnm{Ivan}\binits{I.}}
(\byear{2008}).
\btitle{Structured low-rank approximation and its applications}.
\bjournal{Automatica}
\bvolume{44}
\bpages{891--909}.
\end{barticle}
\endbibitem

\bibitem{markovsky2010bibliography}
\begin{barticle}[author]
\bauthor{\bsnm{Markovsky},~\bfnm{Ivan}\binits{I.}}
(\byear{2010}).
\btitle{Bibliography on total least squares and related methods}.
\bjournal{Statistics and its interface}
\bvolume{3}
\bpages{329--334}.
\end{barticle}
\endbibitem

\bibitem{markovsky2012sysid}
\begin{binproceedings}[author]
\bauthor{\bsnm{Markovsky},~\bfnm{I.}\binits{I.}}
(\byear{2012}).
\btitle{How effective is the nuclear norm heuristic in solving data
  approximation problems?}
In \bbooktitle{Proc. of the 16th IFAC Symposium on System Identification}
\bpages{316--321}.
\bdoi{10.3182/20120711-3-BE-2027.00125}
\end{binproceedings}
\endbibitem

\bibitem{markovsky2019book}
\begin{bbook}[author]
\bauthor{\bsnm{Markovsky},~\bfnm{I.}\binits{I.}}
(\byear{2019}).
\btitle{Low-Rank Approximation: Algorithms, Implementation, Applications},
\bedition{2nd edition} ed.
\bpublisher{Springer}.
\bdoi{10.1007/978-3-319-89620-5}
\end{bbook}
\endbibitem

\bibitem{Markovsky:2018}
\begin{binproceedings}[author]
\bauthor{\bsnm{Markovsky},~\bfnm{I}\binits{I.}} \AND
  \bauthor{\bsnm{Dragotti},~\bfnm{PL}\binits{P.}}
(\byear{2018}).
\btitle{Using Hankel Structured Low-Rank Approximation for Sparse Signal
  Recovery}.
\bpages{479--487}.
\bpublisher{Springer}.
\bdoi{10.1007/978-3-319-93764-9_44}
\end{binproceedings}
\endbibitem

\bibitem{slra-missdata}
\begin{barticle}[author]
\bauthor{\bsnm{Markovsky},~\bfnm{I.}\binits{I.}} \AND
  \bauthor{\bsnm{Usevich},~\bfnm{K.}\binits{K.}}
(\byear{2013}).
\btitle{Structured low-rank approximation with missing data}.
\bjournal{SIAM Journal on Matrix Analysis and Applications}
\bvolume{34}
\bpages{814--830}.
\bdoi{10.1137/120883050}
\end{barticle}
\endbibitem

\bibitem{slra-software}
\begin{barticle}[author]
\bauthor{\bsnm{Markovsky},~\bfnm{I.}\binits{I.}} \AND
  \bauthor{\bsnm{Usevich},~\bfnm{K.}\binits{K.}}
(\byear{2014}).
\btitle{Software for weighted structured low-rank approximation}.
\bjournal{J. Comput. Appl. Math.}
\bvolume{256}
\bpages{278--292}.
\end{barticle}
\endbibitem

\bibitem{markovsky2007overview}
\begin{barticle}[author]
\bauthor{\bsnm{Markovsky},~\bfnm{Ivan}\binits{I.}} \AND
  \bauthor{\bsnm{Van~Huffel},~\bfnm{Sabine}\binits{S.}}
(\byear{2007}).
\btitle{Overview of total least-squares methods}.
\bjournal{Signal processing}
\bvolume{87}
\bpages{2283--2302}.
\end{barticle}
\endbibitem

\bibitem{markovsky2006siam}
\begin{bbook}[author]
\bauthor{\bsnm{Markovsky},~\bfnm{I.}\binits{I.}},
  \bauthor{\bsnm{Willems},~\bfnm{J.~C.}\binits{J.~C.}}, \bauthor{\bsnm{{Van
  Huffel}},~\bfnm{S.}\binits{S.}} \AND \bauthor{\bsnm{{De
  Moor}},~\bfnm{B.}\binits{B.}}
(\byear{2006}).
\btitle{Exact and Approximate Modeling of Linear Systems: A Behavioral
  Approach}.
\bpublisher{SIAM}.
\bdoi{10.1137/1.9780898718263}
\end{bbook}
\endbibitem

\bibitem{mourrain2018polynomial}
\begin{barticle}[author]
\bauthor{\bsnm{Mourrain},~\bfnm{Bernard}\binits{B.}}
(\byear{2018}).
\btitle{Polynomial-exponential decomposition from moments}.
\bjournal{Foundations of Computational Mathematics}
\bvolume{18}
\bpages{1435--1492}.
\end{barticle}
\endbibitem

\bibitem{mourrain2000multivariate}
\begin{barticle}[author]
\bauthor{\bsnm{Mourrain},~\bfnm{Bernard}\binits{B.}} \AND
  \bauthor{\bsnm{Pan},~\bfnm{Victor~Y}\binits{V.~Y.}}
(\byear{2000}).
\btitle{Multivariate polynomials, duality, and structured matrices}.
\bjournal{Journal of complexity}
\bvolume{16}
\bpages{110--180}.
\end{barticle}
\endbibitem

\bibitem{ongie2017irls}
\begin{barticle}[author]
\bauthor{\bsnm{Ongie},~\bfnm{Gregory}\binits{G.}} \AND
  \bauthor{\bsnm{Jacob},~\bfnm{Mathews}\binits{M.}}
(\byear{2017}).
\btitle{A Fast Algorithm for Convolutional Structured Low-Rank Matrix
  Recovery}.
\bjournal{IEEE Transactions on Computational Imaging}
\bvolume{3}
\bpages{535-550}.
\bdoi{10.1109/TCI.2017.2721819}
\end{barticle}
\endbibitem

\bibitem{oropeza2011seismic}
\begin{barticle}[author]
\bauthor{\bsnm{Oropeza},~\bfnm{Vicente}\binits{V.}} \AND
  \bauthor{\bsnm{Sacchi},~\bfnm{Mauricio}\binits{M.}}
(\byear{2011}).
\btitle{Simultaneous seismic data denoising and reconstruction via multichannel
  singular spectrum analysis}.
\bjournal{Geophysics}
\bvolume{76}
\bpages{V25-V32}.
\bdoi{10.1190/1.3552706}
\end{barticle}
\endbibitem

\bibitem{ottaviani2014exact}
\begin{barticle}[author]
\bauthor{\bsnm{Ottaviani},~\bfnm{Giorgio}\binits{G.}},
  \bauthor{\bsnm{Spaenlehauer},~\bfnm{Pierre-Jean}\binits{P.-J.}} \AND
  \bauthor{\bsnm{Sturmfels},~\bfnm{Bernd}\binits{B.}}
(\byear{2014}).
\btitle{Exact Solutions in Structured Low-Rank Approximation}.
\bjournal{SIAM Journal on Matrix Analysis and Applications}
\bvolume{35}
\bpages{1521-1542}.
\bdoi{10.1137/13094520X}
\end{barticle}
\endbibitem

\bibitem{papy2006common}
\begin{barticle}[author]
\bauthor{\bsnm{Papy},~\bfnm{Jean-Michel}\binits{J.-M.}},
  \bauthor{\bsnm{De~Lathauwer},~\bfnm{Lieven}\binits{L.}} \AND
  \bauthor{\bsnm{Van~Huffel},~\bfnm{Sabine}\binits{S.}}
(\byear{2006}).
\btitle{Common pole estimation in multi-channel exponential data modeling}.
\bjournal{Signal Processing}
\bvolume{86}
\bpages{846--858}.
\end{barticle}
\endbibitem

\bibitem{prony1975}
\begin{barticle}[author]
\bauthor{\bsnm{Prony},~\bfnm{G.~R.~B.}\binits{G.~R.~B.}}
(\byear{an III, 1795}).
\btitle{Essai \'{e}xperimental et analytique: sur les lois de la
  dilatabilit\'{e} de fluides \'{e}lastique et sur celles de la force expansive
  de la vapeur de l'alkool, à diff\'{e}rentes temp\'{e}ratures}.
\bjournal{Journal de l'\'{E}cole Polytechnique Flor\'{e}al et Plairial}
\bvolume{1}
\bpages{24--76}.
\bnote{cahier 22}.
\end{barticle}
\endbibitem

\bibitem{razavilar1996spectral}
\begin{binproceedings}[author]
\bauthor{\bsnm{Razavilar},~\bfnm{Javad}\binits{J.}},
  \bauthor{\bsnm{Li},~\bfnm{Ye}\binits{Y.}} \AND
  \bauthor{\bsnm{Liu},~\bfnm{KJ~Ray}\binits{K.~R.}}
(\byear{1996}).
\btitle{Spectral estimation based on structured low rank matrix pencil}.
In \bbooktitle{1996 IEEE International Conference on Acoustics, Speech, and
  Signal Processing Conference Proceedings}
\bvolume{5}
\bpages{2503--2506}.
\bpublisher{IEEE}.
\end{binproceedings}
\endbibitem

\bibitem{recht2010guaranteed}
\begin{barticle}[author]
\bauthor{\bsnm{Recht},~\bfnm{Benjamin}\binits{B.}},
  \bauthor{\bsnm{Fazel},~\bfnm{Maryam}\binits{M.}} \AND
  \bauthor{\bsnm{Parrilo},~\bfnm{Pablo~A}\binits{P.~A.}}
(\byear{2010}).
\btitle{Guaranteed minimum-rank solutions of linear matrix equations via
  nuclear norm minimization}.
\bjournal{SIAM review}
\bvolume{52}
\bpages{471--501}.
\bdoi{10.1137/070697835}
\end{barticle}
\endbibitem

\bibitem{rish2014sparse}
\begin{bbook}[author]
\bauthor{\bsnm{Rish},~\bfnm{Irina}\binits{I.}} \AND
  \bauthor{\bsnm{Grabarnik},~\bfnm{Genady}\binits{G.}}
(\byear{2014}).
\btitle{Sparse modeling: theory, algorithms, and applications}.
\bpublisher{CRC press}.
\end{bbook}
\endbibitem

\bibitem{rohde2011estimation}
\begin{barticle}[author]
\bauthor{\bsnm{Rohde},~\bfnm{Angelika}\binits{A.}} \AND
  \bauthor{\bsnm{Tsybakov},~\bfnm{Alexandre}\binits{A.}}
(\byear{2011}).
\btitle{Estimation of high-dimensional low-rank matrices}.
\bjournal{Annals of Statistics}
\bvolume{39}
\bpages{887--930}.
\end{barticle}
\endbibitem

\bibitem{rosen1996total}
\begin{barticle}[author]
\bauthor{\bsnm{Rosen},~\bfnm{J~Ben}\binits{J.~B.}},
  \bauthor{\bsnm{Park},~\bfnm{Haesun}\binits{H.}} \AND
  \bauthor{\bsnm{Glick},~\bfnm{John}\binits{J.}}
(\byear{1996}).
\btitle{Total least norm formulation and solution for structured problems}.
\bjournal{SIAM Journal on matrix analysis and applications}
\bvolume{17}
\bpages{110--126}.
\end{barticle}
\endbibitem

\bibitem{rump2003structured}
\begin{barticle}[author]
\bauthor{\bsnm{Rump},~\bfnm{S.}\binits{S.}}
(\byear{2003}).
\btitle{Structured Perturbations Part {I}: Normwise Distances.}
\bjournal{SIAM J. Matrix Anal. Appl.}
\bvolume{25}
\bpages{1--30}.
\end{barticle}
\endbibitem

\bibitem{schost2016quadratic}
\begin{barticle}[author]
\bauthor{\bsnm{Schost},~\bfnm{\'{E}ric}\binits{E.}} \AND
  \bauthor{\bsnm{Spaenlehauer},~\bfnm{Pierre-Jean}\binits{P.-J.}}
(\byear{2016}).
\btitle{A Quadratically Convergent Algorithm for Structured Low-Rank
  Approximation}.
\bjournal{Found. Comput. Math.}
\bvolume{16}
\bpages{457–492}.
\bdoi{10.1007/s10208-015-9256-x}
\end{barticle}
\endbibitem

\bibitem{shen2008sparse}
\begin{barticle}[author]
\bauthor{\bsnm{Shen},~\bfnm{Haipeng}\binits{H.}} \AND
  \bauthor{\bsnm{Huang},~\bfnm{Jianhua~Z}\binits{J.~Z.}}
(\byear{2008}).
\btitle{Sparse principal component analysis via regularized low rank matrix
  approximation}.
\bjournal{Journal of multivariate analysis}
\bvolume{99}
\bpages{1015--1034}.
\end{barticle}
\endbibitem

\bibitem{shin2014mri}
\begin{barticle}[author]
\bauthor{\bsnm{Shin},~\bfnm{Peter~J.}\binits{P.~J.}},
  \bauthor{\bsnm{Larson},~\bfnm{Peder E.~Z.}\binits{P.~E.~Z.}},
  \bauthor{\bsnm{Ohliger},~\bfnm{Michael~A.}\binits{M.~A.}},
  \bauthor{\bsnm{Elad},~\bfnm{Michael}\binits{M.}},
  \bauthor{\bsnm{Pauly},~\bfnm{John~M.}\binits{J.~M.}},
  \bauthor{\bsnm{Vigneron},~\bfnm{Daniel~B.}\binits{D.~B.}} \AND
  \bauthor{\bsnm{Lustig},~\bfnm{Michael}\binits{M.}}
(\byear{2014}).
\btitle{Calibrationless parallel imaging reconstruction based on structured
  low-rank matrix completion}.
\bjournal{Magnetic Resonance in Medicine}
\bvolume{72}
\bpages{959--970}.
\bdoi{10.1002/mrm.24997}
\end{barticle}
\endbibitem

\bibitem{shlemov2014shaped}
\begin{binproceedings}[author]
\bauthor{\bsnm{Shlemov},~\bfnm{A.}\binits{A.}} \AND
  \bauthor{\bsnm{Golyandina},~\bfnm{N.}\binits{N.}}
(\byear{2014}).
\btitle{Shaped extension of singular spectrum analysis}.
In \bbooktitle{Proceedings of the 21st International Symposium on Mathematical
  Theory of Networks and Systems (MTNS 2014)}.
\end{binproceedings}
\endbibitem

\bibitem{srebro2003weighted}
\begin{binproceedings}[author]
\bauthor{\bsnm{Srebro},~\bfnm{Nathan}\binits{N.}} \AND
  \bauthor{\bsnm{Jaakkola},~\bfnm{Tommi}\binits{T.}}
(\byear{2003}).
\btitle{Weighted low-rank approximations}.
In \bbooktitle{Proceedings of the 20th International Conference on Machine
  Learning (ICML-03)}
\bpages{720--727}.
\end{binproceedings}
\endbibitem

\bibitem{stoica2005spectral}
\begin{barticle}[author]
\bauthor{\bsnm{Stoica},~\bfnm{Petre}\binits{P.}},
  \bauthor{\bsnm{Moses},~\bfnm{Randolph~L}\binits{R.~L.}} \betal{et~al.}
(\byear{2005}).
\btitle{Spectral analysis of signals}.
\end{barticle}
\endbibitem

\bibitem{sylvester1886binary}
\begin{barticle}[author]
\bauthor{\bsnm{Sylvester},~\bfnm{J.~J.}\binits{J.~J.}}
(\byear{1886}).
\btitle{Sur une extension d’un th\'{e}or\`{e}me de Clebsch relatif aux
  courbes du quatri\`{e}me degr\'{e}}.
\bjournal{Comptes Rendus, Math. Acad. Sci. Paris}
\bvolume{102}
\bpages{1532--1534}.
\end{barticle}
\endbibitem

\bibitem{udell2019big}
\begin{barticle}[author]
\bauthor{\bsnm{Udell},~\bfnm{Madeleine}\binits{M.}} \AND
  \bauthor{\bsnm{Townsend},~\bfnm{Alex}\binits{A.}}
(\byear{2019}).
\btitle{Why are big data matrices approximately low rank?}
\bjournal{SIAM Journal on Mathematics of Data Science}
\bvolume{1}
\bpages{144--160}.
\end{barticle}
\endbibitem

\bibitem{usevich2010signal}
\begin{barticle}[author]
\bauthor{\bsnm{Usevich},~\bfnm{K.}\binits{K.}}
(\byear{2010}).
\btitle{On signal and extraneous roots in {Singular Spectrum Analysis}}.
\bjournal{Statistics and Its Interface}
\bvolume{3}
\bpages{281--295}.
\bdoi{10.4310/SII.2010.v3.n3.a3}
\end{barticle}
\endbibitem

\bibitem{usevich2016hankel}
\begin{barticle}[author]
\bauthor{\bsnm{Usevich},~\bfnm{Konstantin}\binits{K.}} \AND
  \bauthor{\bsnm{Comon},~\bfnm{Pierre}\binits{P.}}
(\byear{2016}).
\btitle{Hankel Low-Rank Matrix Completion: Performance of the Nuclear Norm
  Relaxation}.
\bjournal{IEEE Journal of Selected Topics in Signal Processing}
\bvolume{10}
\bpages{637--646}.
\end{barticle}
\endbibitem

\bibitem{usevich2018stft}
\begin{binproceedings}[author]
\bauthor{\bsnm{Usevich},~\bfnm{Konstantin}\binits{K.}},
  \bauthor{\bsnm{Emiya},~\bfnm{Valentin}\binits{V.}},
  \bauthor{\bsnm{Brie},~\bfnm{David}\binits{D.}} \AND
  \bauthor{\bsnm{Chaux},~\bfnm{Caroline}\binits{C.}}
(\byear{2018}).
\btitle{Characterization of Finite Signals with Low-Rank {STFT}}.
In \bbooktitle{2018 IEEE Statistical Signal Processing Workshop (SSP)}
\bpages{393-397}.
\bdoi{10.1109/SSP.2018.8450745}
\end{binproceedings}
\endbibitem

\bibitem{usevich2014varpro}
\begin{barticle}[author]
\bauthor{\bsnm{Usevich},~\bfnm{K.}\binits{K.}} \AND
  \bauthor{\bsnm{Markovsky},~\bfnm{I.}\binits{I.}}
(\byear{2014}).
\btitle{Variable projection for affinely structured low-rank approximation in
  weighted 2-norms}.
\bjournal{J. Comput. Appl. Math.}
\bvolume{272}
\bpages{430--448}.
\bdoi{10.1016/j.cam.2013.04.034}
\end{barticle}
\endbibitem

\bibitem{usevich2017agcd}
\begin{barticle}[author]
\bauthor{\bsnm{Usevich},~\bfnm{Konstantin}\binits{K.}} \AND
  \bauthor{\bsnm{Markovsky},~\bfnm{Ivan}\binits{I.}}
(\byear{2017}).
\btitle{Variable projection methods for approximate (greatest) common divisor
  computations}.
\bjournal{Theoretical Computer Science}
\bvolume{681}
\bpages{176-198}.
\bnote{Symbolic Numeric Computation}.
\bdoi{https://doi.org/10.1016/j.tcs.2017.03.028}
\end{barticle}
\endbibitem

\bibitem{van1996formulation}
\begin{barticle}[author]
\bauthor{\bsnm{Van~Huffel},~\bfnm{Sabine}\binits{S.}},
  \bauthor{\bsnm{Park},~\bfnm{Haesun}\binits{H.}} \AND
  \bauthor{\bsnm{Rosen},~\bfnm{J~Ben}\binits{J.~B.}}
(\byear{1996}).
\btitle{Formulation and solution of structured total least norm problems for
  parameter estimation}.
\bjournal{IEEE Transactions on signal processing}
\bvolume{44}
\bpages{2464--2474}.
\end{barticle}
\endbibitem

\bibitem{willems1986linear}
\begin{barticle}[author]
\bauthor{\bsnm{Willems},~\bfnm{Jan~C.}\binits{J.~C.}}
(\byear{1986}).
\btitle{From time series to linear system--- Part I. Finite dimensional linear
  time invariant systems}.
\bjournal{Automatica}
\bvolume{22}
\bpages{561-580}.
\bdoi{https://doi.org/10.1016/0005-1098(86)90066-X}
\end{barticle}
\endbibitem

\bibitem{yang2016vandermonde}
\begin{barticle}[author]
\bauthor{\bsnm{Yang},~\bfnm{Zai}\binits{Z.}},
  \bauthor{\bsnm{Xie},~\bfnm{Lihua}\binits{L.}} \AND
  \bauthor{\bsnm{Stoica},~\bfnm{Petre}\binits{P.}}
(\byear{2016}).
\btitle{Vandermonde decomposition of multilevel {Toeplitz} matrices with
  application to multidimensional super-resolution}.
\bjournal{IEEE Transactions on Information Theory}
\bvolume{62}
\bpages{3685--3701}.
\end{barticle}
\endbibitem

\bibitem{zhang2019optimal}
\begin{barticle}[author]
\bauthor{\bsnm{Zhang},~\bfnm{Ran}\binits{R.}} \AND
  \bauthor{\bsnm{Plonka},~\bfnm{Gerlind}\binits{G.}}
(\byear{2019}).
\btitle{Optimal approximation with exponential sums by a maximum likelihood
  modification of Prony's method.}
\bjournal{Adv. Comput. Math.}
\bvolume{45}
\bpages{1657--1687}.
\end{barticle}
\endbibitem

\bibitem{zvonarev2017iterative}
\begin{barticle}[author]
\bauthor{\bsnm{Zvonarev},~\bfnm{Nikita}\binits{N.}} \AND
  \bauthor{\bsnm{Golyandina},~\bfnm{Nina}\binits{N.}}
(\byear{2017}).
\btitle{Iterative algorithms for weighted and unweighted finite-rank
  time-series approximations}.
\bjournal{Statistics and its Interface}
\bvolume{10}
\bpages{5--18}.
\bdoi{10.4310/SII.2017.v10.n1.a1}
\end{barticle}
\endbibitem

\bibitem{zvonarev2021low}
\begin{barticle}[author]
\bauthor{\bsnm{Zvonarev},~\bfnm{Nikita}\binits{N.}} \AND
  \bauthor{\bsnm{Golyandina},~\bfnm{Nina}\binits{N.}}
(\byear{2021}).
\btitle{Low-rank signal subspace: parameterization, projection and signal
  estimation}.
\bjournal{arXiv preprint arXiv:2101.09779}.
\end{barticle}
\endbibitem

\bibitem{zvonarev2021fast}
\begin{barticle}[author]
\bauthor{\bsnm{Zvonarev},~\bfnm{Nikita}\binits{N.}} \AND
  \bauthor{\bsnm{Golyandina},~\bfnm{Nina}\binits{N.}}
(\byear{2021}).
\btitle{Fast and stable modification of the Gauss-Newton method for low-rank
  signal estimation}.
\bjournal{arXiv preprint arXiv:2106.14215}.
\end{barticle}
\endbibitem

\end{thebibliography}
\end{document}